\documentclass{amsart}

\usepackage{amsfonts}
\usepackage{amsmath}
\usepackage{amsthm}
\usepackage{amssymb}
\usepackage[english]{babel}
\usepackage{graphics}
\usepackage{latexsym}
\usepackage{longtable} 
\usepackage{mathrsfs} 
\usepackage{graphicx}
\usepackage{wrapfig} %\usepackage{txfonts}
\usepackage[pdftex,colorlinks=true,citecolor=cyan]{hyperref}
\usepackage{enumitem}
\usepackage{esint}
\usepackage{subfig}
\usepackage{caption}

\usepackage{tikz}
\usetikzlibrary{arrows,%
                plotmarks}

\usepackage{pgfplots}
\pgfplotsset{compat=1.11}

\newtheorem{theorem}{Theorem}
\newtheorem{corollary}[theorem]{Corollary}
\newtheorem{lemma}[theorem]{Lemma}

\newtheorem{claim}[theorem]{Claim}
\newtheorem{example}[theorem]{Example}

\theoremstyle{definition}
\newtheorem{definition}[theorem]{Definition}

\newcommand{\mL}{\mathcal{L}}
\newcommand{\mH}{\mathcal{H}}
\newcommand{\mF}{\mathcal{F}}

\newcommand{\mO}{\mathcal{O}}

\newcommand{\E}{\mathrm{E}}
\newcommand{\D}{\mathrm{D}}
\renewcommand{\H}{\mathrm{H}}

\newcommand{\R}{\mathbb{R}}

\newcommand{\N}{\mathbb{N}}

\newcommand{\noi}{\noindent}
\newcommand{\ms}{\medskip}

\newcommand{\al}{\alpha}

\newcommand{\Ga}{\Gamma}
\newcommand{\de}{\delta}

\newcommand{\si}{\sigma}

\newcommand{\Om}{\Omega}

\newcommand{\weak }{\, -\!\!\!\!-\!\!\!\rightharpoonup}
\newcommand{\weakstar }{ \overset{\, *_{\phantom{|}}}{{\smash{\weak }}\, } }

\newcommand{\larrow}{\longrightarrow}

\newcommand{\ri}{\rightarrow}
\newcommand{\LL}{\text{\LARGE$\llcorner$}}
\newcommand{\p}{\partial}
\newcommand{\sub}{\subseteq}
\newcommand{\set}{\setminus}
\newcommand{\by}{\times}

\newcommand{\tr}{\mathrm{tr}}
\renewcommand{\div}{\mathrm{div}}

\newcommand{\bt}{\begin{theorem}}\newcommand{\et}{\end{theorem}}
\newcommand{\bd}{\begin{definition}}\newcommand{\ed}{\end{definition}}
\newcommand{\bl}{\begin{lemma}}\newcommand{\el}{\end{lemma}}
\newcommand{\beq}{\begin{equation}}\newcommand{\eeq}{\end{equation}}
\newcommand{\bc}{\begin{claim}}\newcommand{\ec}{\end{claim}}
\newcommand{\bex}{\begin{example}}\newcommand{\eex}{\end{example}}
\newcommand{\bcor}{\begin{corollary}}\newcommand{\ecor}{\end{corollary}}
\newcommand{\bp}{\begin{proof}}\newcommand{\ep}{\end{proof}}

\numberwithin{equation}{section}
\begin{document}

\title[Numerics for vectorial absolute minimisers]{On the numerical approximation of vectorial absolute minimisers in $L^\infty$}

\author{Nikos Katzourakis and Tristan Pryer}

\address{Department of Mathematics and Statistics, University of Reading, Whiteknights, PO Box 220, Reading RG6 6AX, United Kingdom}
\email{n.katzourakis@reading.ac.uk}

\address{Department of Mathematics and Statistics, University of Reading, Whiteknights, PO Box 220, Reading RG6 6AX, United Kingdom}

  \thanks{\!\!\!\!\!\!\!\!\texttt{N.K. has been partially financially supported through the EPSRC grant EP/N017412/1. T.P. has been partially financially supported through the EPSRC grant EP/P000835/1.}}
  
  \email{t.pryer@reading.ac.uk}
%\subjclass[2010]{35D99, 35D40, 35J47, 35J47, 35J92, 35J70, 35J99}

\date{}


\begin{abstract} Let $\Omega$ be an open set. We consider the supremal functional
\[ \tag{1} \label{1}
\ \ \ \ \ \ \mathrm{E}_\infty (u,\mathcal{O})\, :=\, \| \D u \|_{L^\infty( \mathcal{O} )}, \ \ \ \mathcal{O} \subseteq \Omega \text{ open},
\]
applied to locally Lipschitz mappings $u : \R^n \supseteq \Om \larrow
\R^N$, where $n,N\in \N$. This is the model functional of Calculus of
Variations in $L^\infty$. The area is developing rapidly, but the
vectorial case of $N\geq 2$ is still poorly understood. Due to the
non-local nature of \eqref{1}, usual minimisers are not truly
optimal. The concept of so-called absolute minimisers is the primary
contender in the direction of variational concepts. However, these
cannot be obtained by direct minimisation and the question of their
existence under prescribed boundary data is open when $n,N\geq
2$. Herein we present numerical experiments based on a new method
recently proposed by the first author in the papers
\cite{KM2,KS2}.

\end{abstract}

\maketitle

%\tableofcontents

\section{Introduction} \label{section1}

Let $n,N\in \N$ be fixed integers and $\Omega$ an open set in $\R^n$. In this paper we perform numerical experiments and make some theoretical observations regarding appropriately defined minimisers of the functional
\beq
\label{1.1}
\ \ \ \ \ \ \mathrm{E}_\infty (u,\mathcal{O})\, :=\, \| \D u \|_{L^\infty( \mathcal{O} )}, \ \ \ \mathcal{O} \subseteq \Omega \text{ open},
\eeq
which is the archetypal model functionals of the Calculus of Variations in the space $L^\infty$. The appropriate functional setting to place this functional is the Sobolev space $W^{1,\infty}(\Om;\R^N)$, which (by the Rademacher theorem) consists of a.e.\ differentiable continuous mappings $u : \R^n \supseteq \Om \larrow \R^N$ with $L^\infty$ derivative $\D u$. Here the gradient is understood as an $N\by n$ matrix-valued map 
\[
\D u=(\D_i u_\al)_{i=1...n}^{\al=1...N} \ : \ \ \Om \larrow \R^{N\by n}
\]
where $\D_i \equiv \p/\p x_i$ and the $L^\infty$ norm in \eqref{1.1} is understood as the essential supremum of the Euclidean (Frobenius) norm on $\R^{N\by n}$, defined through the associated trace inner product $A:B := \tr (A^\top B)$.

The area of Calculus of Variations in $L^\infty$ has a relatively short history in Analysis, pioneered by G. Aronsson in the 1960s, see e.g.\ \cite{A1, A2, A3}. Nonetheless, since its inception it has attracted the interest of many mathematicians due to both the theoretical importance (see e.g.\ \cite{ACJS, BJW1, BJW2, CEG, J, EY, P} and the expository texts \cite{ACJ, C, K8}) as well as due to the relevance to various and diverse applications (see for instance \cite{BN, GNP, ETT, PSSW, BEJ, CR, K1}).

One of the starting difficulties to the study of \eqref{1.1} is that, quite surprisingly, the standard (global) minimisers are not truly optimal objects any more, in the sense that ``they can be minimised even further" and in general may not solve any kind of Euler-Lagrange equations. Namely, there are many simple examples even in one space dimension where a minimiser does not minimise with respect to its own boundary conditions on subdomains. The latter property is of course automatic for integral functionals. In essence, this owes to the ``non-local" nature of \eqref{1.1}: although one can see \eqref{1.1} as the limit case of the $p$-Dirichlet functional as $p\to \infty$
\beq
\label{1.2}
\ \ \ \ \ \ \mathrm{E}_p (u,\mathcal{O})\, :=\, \left(\int_\mO | \D u |^p\right)^{\!\!1/p}, \ \ \mathcal{O} \subseteq \Omega \text{ open}, \ \ \ u \in W^{1,p}(\Om;\R^N),
\eeq
when passing to the extreme case of $p=\infty$, the $\sigma$-additivity properties of the Lebesgue integral are lost. This was already realised by Aronsson, who introduced the next concept.
\bd[Absolute minimisers]  
\label{definition1}
A mapping $u\in W^{1,\infty}(\Om;\R^N)$ is called an absolute minimiser of \eqref{1.1} on $\Om$ if
\[
\E_\infty (u,\mO)\, \leq \, \E_\infty(u+\phi,\mO)
\]
for any open sets $\mO \sub \Om$ and any variation $\phi \in W^{1,\infty}_0(\mO;\R^N)$.
\ed

Although the area is developing rapidly, the vectorial case of $N\geq 2$ is still poorly understood. In particular, to date neither the existence nor the uniqueness of absolute minimisers with prescribed boundary data on $\p\Om$ are known, unless $\min\{n,N\}=1$, namely either when $n=1$ (curves in $\R^N$) or $N=1$ (scalar functions on $\R^n$), even in more general situations when one has a general Hamiltonian $\H(.,u,\D u)$ instead of $|\D u|$. The difficulty lies in the fact that for the ``localised" variational concept of Definition \ref{definition1}, the standard approach of direct minimisation (\cite{D,FL,GM}) of the functional \eqref{1.1} in the affine space 
\[
W^{1,\infty}_{u_0}(\Om;\R^N) \, := \, u_0 + W^{1,\infty}_{0}(\Om;\R^N)
\]
for a fixed $u_0 \in W^{1,\infty}(\Om;\R^N)$ does {\it not} in general yield absolute minimisers, but merely global minimisers which minimise on $\Om$ only. Hence, in the absence of alternative methods, the main vectorial tool is to resort to $L^p$ approximations and then let $p\to \infty$. To this end, let $(u_p)_p$ be the family of $p$-minimisers in $W^{1,\infty}_{u_0}(\Om;\R^N)$. The problem that one encounters in the attempt to pass to the limit in
\[
\E_p (u_p,\mO)\, \leq \, \E_p(u_p+\phi,\mO),
\]
is that the weak* sequential compactness $u_p \weakstar u_\infty$ one easily has is not strong enough to obtain that the limit is an absolute minimiser. This obstacle can be bypassed when $\min\{n,N\}=1$ because of either the scalar nature of the competing functions and comparison methods, or because of the one-dimensionality of their domain of definition (see for instance \cite{BJW1, K9, AbK}), although the mode of convergence itself cannot in general be strengthened. In particular, the functional \eqref{1.1} is only weakly* lower semi-continuous and {\it not} weakly* continuous.
 
Herein we present and perform numerical experiments based on a new
method recently proposed by the first author in the joint papers
\cite{KM2, KS2}, where various existence and uniqueness have been
obtained. This approach is motivated by the paper \cite{EY} and
relates to developments in mass transportation (see
e.g.\ \cite{DEP,DP,S}, underpinning also the approach followed in the
recent paper \cite{KM}, despite in different guises since the higher
order case treated therein is quite special. We refrain from giving
any details now, and instead expound on the main ideas of this method
in some detail in Section \ref{section2}.

The numerical method we employ is a finite element approximation,
based on an earlier work of Barrett and Liu on numerical methods for
elliptic systems \cite{BL}. Therein the authors prove that, for a
fixed exponent $p$, the method converges to the respective
$p$-harmonic mapping under certain regularity assumptions on the
solution.  We would like to stress that significant care must be taken
with numerical computations using this approach because the underlying
nonlinear system is ill-conditioned. This owes to the nonlinearity of
the problem which grows exponentially with $p$. Work to overcome this
issue includes, for example, the work of Huang, Li and Liu
\cite{HuangLiLiu:2007} where preconditioners based on gradient descent
algorithms are designed and shown to work well for $p$ up to
$1000$. 

Our mechanism for designing approximations of absolute minimisers is
based on the approximation of the $L^\infty$ minimisation problem by a
respective $L^p$ one. We then utilise the approach described in
\cite{Pry,KP,KP2} where it was shown that by forming an appropriate
limit we are able to select candidates for numerical approximation
along a ``good'' sequence of solutions, the $p$-harmonic mappings.

The purpose of this work is to demonstrate some key properties of
vectorial absolute minimisers using an analytically justifiable
numerical scheme which currently is the \emph{only} technique
available to give insight into the limiting vector-valued problem. We
note that our goal is \emph{not} to construct an efficient
approximation method for vectorial absolute minimisers; indeed this
indirect approximation of the limiting problem is not computationally
efficient.

We conclude by noting that, albeit the concept of absolute minimisers
is the primary contender in the direction of variational concepts, is
not the unique candidate in the vector case of $N\geq 2$. This is
connected to that fact that, when $N\geq2$, the associated
Euler-Lagrange $\infty$-Laplace system admits smooth non-minimising
solutions which are characterised by two distinct variational
concepts, see e.g.\ \cite{AK, K2, KS}.

\ms

\noi \textbf{Conflict of interest.} The authors declare that there is conflict of interest.

\ms

\section{Approximating absolute minimisers through $p$-concentration measures} \label{section2}

In this section we present the recently proposed method to construct absolute minimisers. The ideas herein are taken from the papers in preparation \cite{KS2,KM2}. For the clarity of this discussion suppose that $u \in C^1(\Om;\R^N)$. The central point of this approach is to bypass the difficulties caused by the lack of continuity of the essential-supremum with respect to weak convergence by attempting to write the supremum as an integral but for a different quirky measure, namely
\beq
\label{2.1}
\E_\infty(u,\mO)\, =\, \int_{\overline{\mO}} |\D u| \, \mathrm{d} \si
\eeq
where $\si$ is a Radon probability measure on the closure $\overline{\mO}$, which in principle is $\mO$-dependent as well as $u$-dependent. Of course one can define an infinity of different families of measures $\{\si^{\mO}: \mO\Subset \Om\}$ performing this role: for instance, for any $\mO \Subset \Om$ choose any point $x_{\mO}$ in the arg-max set of $|\D u|$, namely any point $x_{\mO}  \in \overline{\mO}$ such that
\[
|\D u|(x_{\mO}) \, =\, \E_\infty(u,\mO)
\]
and then choose the Dirac mass $\si^{\mO}:= \de_{x_{\mO}}$. However, this condition by itself does not suffice for our purposes and has to be coupled by another condition. Firstly, in order to utilise duality arguments with ease, let us modify \eqref{2.1} to its $L^2$ variant
\beq
\label{2.2}
\E_\infty(u,\mO)^2\, =\, \int_{\overline{\mO}} |\D u|^2 \, \mathrm{d} \si
\eeq
and suppose in addition that the matrix-valued Radon measure $\D u \, \si$ is divergence-free on $\overline{\mO}$, namely that
 \beq
\label{2.3}
\div\big( \D u \, \si \big)\, =\, 0.
\eeq
The PDE \eqref{2.3} is to be understood in a sense stronger than the usual distributional sense, due to the emergence of the boundary. In fact, \eqref{2.3} has to be interpreted as
\[
\int_{\overline{\mO}} \D u : \D \phi \, \mathrm{d} \si \, =\, 0, \ \ \text{ for all }\phi \in C^1_0(\overline{\mO};\R^N),
\]
namely in the space of those $C^1$ test maps $\phi : \mO \larrow \R^n$ which vanish on $\p\mO$ and extend continuously together with their derivative on $\p\mO$ (but whose derivative may not vanish on $\p\mO$). 

Although the pair of conditions \eqref{2.2}-\eqref{2.3} seem to be ad-hoc and unmotivated, the payoff is that it is quite easy to see that if one has a $C^1$ map $u : \Om \larrow \R^N$ such that for any $\mO \Subset \Om$ there is a Radon probability measure $\si^\mO$ on $\overline{\mO}$ such that the pair $u,\si^\mO$ satisfies \eqref{2.2}-\eqref{2.3}, then in fact $u$ is an absolute minimiser in the sense of Aronsson (Definition \ref{1.1}). Indeed, fix a test function $\phi \in C^1_0(\overline{\mO};\R^N)$ and for simplicity suppose that $\p\mO$ has zero Lebesgue $n$-measure. Then, we have
\[
\begin{split}
\E_\infty(u,\mO)^2 &\overset{\eqref{2.2}}{=}\, \int_{\overline{\mO}} |\D u|^2 \, \mathrm{d} \si^\mO 
\\
& \overset{\eqref{2.3}}{=}\, \int_{\overline{\mO}} |\D u|^2 \, \mathrm{d} \si \, +\, 2 \int_{\overline{\mO}} \D u : \D \phi \, \mathrm{d} \si^\mO 
\\
& \ \leq \, \int_{\overline{\mO}} |\D u|^2 \, \mathrm{d} \si^\mO \, +\, 2 \int_{\overline{\mO}} \D u : \D \phi \, \mathrm{d}\si^\mO \, +\, \int_{\overline{\mO}} |\D \phi|^2 \, \mathrm{d} \si^\mO
\\
& \ =\, \int_{\overline{\mO}} |\D u + \D \phi|^2 \, \mathrm{d} \si^\mO 
\\
&\ \leq \, \si^\mO(\overline{\mO}) \, \sup_{\overline{\mO}} |\D u + \D \phi|^2
\\
&\ =\, \E_\infty(u+\phi,\mO)^2.
\end{split}
\]
Hence, $u$ is indeed minimising on $\mO$, at least among smooth variations. The general case then needs some approximation arguments, appropriately adapted to $L^\infty$. What is less clear is how one can actually prove the existence of such objects in $W^{1,\infty}_{u_0}(\Om;\R^N)$ with given boundary conditions. The idea is to use $L^p$ approximations in the following way: for each $p>n$, consider the minimisation problem
\[
\E_p(u_p,\Om)\, =\, \inf\Big\{\E_p(u,\Om)\ : \ u \in W^{1,p}_{u_0}(\Om;\R^N) \Big\}.
\]
By standard variational arguments (see e.g.\ \cite{D,FL}), this problem has a solution $u_p$, which is a weak solution to the celebrated $p$-Laplace system:
\beq
\div \big(|\D u|^{p-2}\D u \big)\, =\, 0, \ \ \text{ in }\Om.
\eeq
In particular, this means that for any $\mO \Subset \Om$ and any $\phi \in W^{1,p}_0(\mO;\R^N)$, $u_p$ satisfies
\beq 
\label{2.5}
\int_\mO \D u_p :\D \phi \, |\D u_p|^{p-2}\, =\, 0
\eeq
By defining the (absolutely continuous) probability measure 
\beq
\si^{\mO}_p(A) \,:=\, \frac{ \ \displaystyle \int_{A\cap \mO} |\D u_p|^{p-2} \ }{ \ \displaystyle \int_\mO |\D u_p|^{p-2} \ },  \ \ \ A\sub \Om \text{ Borel},
\eeq
we may rewrite \eqref{2.5} as
\beq 
\label{2.7}
\int_\mO \D u_p :\D \phi \, \mathrm{d} \si^{\mO}_p\, =\, 0.
\eeq
By juxtaposing \eqref{2.7} with \eqref{2.3}, one sees that upon devising appropriate analytic tools in order to pass in an appropriate weak* sense to some limit
\[
(\D u_p,\si^{\mO}_p) \weakstar (\D u_\infty, \si^{\mO}_\infty) 
\]
subsequentially as $p\to \infty$, one would in principle obtain
\eqref{2.3}. This is particularly challenging, as one the one hand
mass might be lost towards the boundary $\p\mO$ and one has to
consider measures on $\overline{\mO}$ rather than on $\mO$, even
though by definition $\si^{\mO}_p (\p\mO)=0$. On the other hand, the
main cause of additional ramifications is the possible ``oscillations"
of the pair of weakly* converging objects $(\D u_p,\si^{\mO}_p)$ and
one needs particular ``compensated compactness" mechanisms to pass to
the limit and show that the oscillations occur in such a way that
cancel each other. Despite being highly non-obvious, the above
approach can indeed work and this is done in \cite{KS2, KM2}.

It remains to show that, at least formally, one also obtains $\eqref{2.2}$ in the limit as $p\to \infty$. To see this, note that for any $\al \in (0,1)$ and for 
\[
\mathrm{L}_p(\mO)\,:=\, \left( \int_\mO |\D u_p|^{p-2}\right)^{\!\!1/(p-2)},
\]
directly from the definition of $\si^{\mO}_p$ we have the estimate
\[
\begin{split}
\si^{\mO}_p\big( \big\{|\D u_p|\leq \al \mathrm{L}_p(\mO) \big\} \big)\, & =\, \frac{ \ \displaystyle \int_{ \{|\D u_p|\leq \al \mathrm{L}_p(\mO) \} \cap \mO} |\D u_p|^{p-2} \ }{ \ \displaystyle \int_\mO |\D u_p|^{p-2} \ }
\\
& \leq \, \al^{p-2}\frac{ \ \displaystyle \int_{ \mO} (\mathrm{L}_p)^{p-2} \ }{ \ \displaystyle \int_\mO |\D u_p|^{p-2} \ }
\\
& = \, \al^{p-2}.
\end{split}
\]
Since one can check with simple arguments that 
\[
\E_\infty(u_\infty,\mO) \, \leq \, \underset{p\to \infty}{\liminf} \, \mathrm{L}_p(\mO)
\]
(with semi-continuity methods similar to e.g.\ in \cite[Ch.\ 8]{K8}), we see that, at least formally
\[
\si^{\mO}_\infty\big( \big\{|\D u_\infty|\leq \al E_\infty(u_\infty,\mO) \big\} \big)\, =
\, 0, \ \ \text{ for all }\al \in (0,1),
\]
which implies
\beq
\label{2.8}
\si^{\mO}_\infty\big( \big\{|\D u_\infty| < E_\infty(u_\infty,\mO) \big\} \big)\, =
\, 0.
\eeq
Equality \eqref{2.8} implies that the limit probability measure $\si^{\mO}_\infty$ is supported on the arg-max set 
\[
\overline{\mO} \cap \big\{|\D u_\infty| = E_\infty(u_\infty,\mO) \big\}, 
\]
since its complement is a nullset with respect to $\si^{\mO}_\infty$. Then, we have
\[
\begin{split}
\int_{\overline{\mO}} |\D u_\infty|^2\, \mathrm{d}\si^{\mO}_\infty\, &=\, \int_{\overline{\mO} \cap  \{|\D u_\infty| = E_\infty(u_\infty,\mO) \} } |\D u_\infty|^2\, \mathrm{d}\si^{\mO}_\infty
\\
&=\, E_\infty(u_\infty,\mO)^2 \, \si^{\mO}_\infty \big(\overline{\mO} \set \big\{|\D u_\infty| < E_\infty(u_\infty,\mO) \big\}\big)
\\
&=\, E_\infty(u_\infty,\mO)^2 \, \Big(\si^{\mO}_\infty \big(\overline{\mO}\big) - \si^{\mO}_\infty \big(\big\{|\D u_\infty| < E_\infty(u_\infty,\mO) \big\}\big)\Big)
\\
&=\, E_\infty(u_\infty,\mO)^2 .
\end{split}
\]
Hence, \eqref{2.2} indeed follows.

\newcommand{\T}{\mathcal T}

\newcommand{\poly}[1]{\ensuremath{\mathbb P}^{#1}}

\newcommand{\fes}{\mathbb V}

\newcommand{\sob}[2]{\ensuremath{W^{#1,#2}}}

\newcommand{\leb}[1]{\ensuremath{L^{#1}}}

\ms

\section{Numerical Approximations of $p$-Concentration Measures}
\label{section3}

In this section we perform several numerical experiments in both the
scalar and the vectorial case for various appropriately selected
boundary conditions. These experiments demonstrate in a concrete
fashion the plethora of possible behaviours of the limit measures
$\si^{\mO}_\infty$, depending on the shape of $\mO$ and its position
in relation to the level sets of the modulus of the gradient $|\D
u_\infty|$.

The method we use is a conforming finite element discretisation of the
$p$-Laplace system analysed in \cite{BL} for fixed $p$. We will
describe the discretisation and summarise extensive numerical
experiments aimed at quantifying the behaviour of the limit measure
$\si^{\mO}_\infty$. We refer to \cite{Pry,KP,KP2} for analytic
justification.

We let $\T{}$ be an admissible triangulation of $\Omega$, namely,
$\T{}$ is a finite collection of sets such that

\smallskip
\begin{enumerate}
\item $K\in\T{}$ implies $K$ is an open triangle,
\item for any $K,J\in\T{}$ we have that $\overline K\cap\overline J$ is
  either $\emptyset$, a vertex, an
  edge, or the whole of $\overline K$ and $\overline J$ and
\item $\cup_{K\in\T{}}\overline K=\overline\Omega$.
\end{enumerate}

\smallskip

\noi The shape regularity constant of $\T{}$ is defined as the number
\begin{equation}
  \label{eqn:def:shape-regularity}
  \mu(\T{}) \,:=\, \inf_{K\in\T{}} \frac{\rho_K}{h_K},
\end{equation}
where $\rho_K$ is the radius of the largest ball contained inside
$K$ and $h_K$ is the diameter of $K$. An indexed family of
triangulations $\{\T^n\}_n$ is called \emph{shape regular} if 
\begin{equation}
  \label{eqn:def:family-shape-regularity}
  \mu\, :=\, \inf_n\mu(\T^n)\, >\, 0.
\end{equation}
Further, we define $h:\Omega \to \R$ to be the {piecewise constant}
\emph{meshsize function} of $\T{}$ given by
\begin{equation}
h({x})\, :=\, \max_{\overline K\ni {x}}\, h_K.
\end{equation}
A mesh is called quasiuniform when there exists a positive constant
$C$ such that $\max_{x\in\Omega} h(x) \le C \min_{x\in\Omega} h(x)$. In what
follows we shall assume that all triangulations are shape-regular and
quasiuniform.

We let $\poly 1(\T{})$ denote the space of piecewise linear
polynomials over the triangulation $\T{}$, i.e.,
\begin{equation}
  \poly 1 (\T{}) \, = \, \big\{ \phi \text{ such that } \phi|_K \in \poly 1 (K) \big\}
\end{equation}
 and introduce the \emph{finite element space}
\begin{gather}
  \label{eqn:def:finite-element-space}
  \fes \,:= \, \poly 1(\T{}) \cap C^{0}(\Omega)
\end{gather}
to be the usual space of continuous piecewise linear polynomial
functions.

\smallskip

\subsection{Galerkin discretisation}
We consider the Galerkin discretisation, to find $U \in [\fes]^N$ with
$U\vert_{\partial \Omega} = I_h g$, the piecewise Lagrange
interpolant, such that
\begin{equation}
  \label{eq:plapdis}
  \int_\Omega |\D U |^{p-2} \D U : {\D \Phi} \,=\, 0,\ \ \forall \Phi \in [\fes]^N.
\end{equation}
This is a conforming finite element discretisation of the vectorial
$p$-Laplacian system proposed in \cite{BL}. Existence and uniqueness
of solution to (\ref{eq:plapdis}) follows from examination of the
$p$-functional
  \begin{equation}
    \label{eq:p-energy}
    \E_p(u,\Omega) \, =\, \left(\int_\Om |\text{D}u|^p \right)^{1/p}.
  \end{equation}
  Notice that (\ref{eq:p-energy}) is strictly convex and coercive on
  $\sob{1}{p}_0(\Omega,\R^N)$ so we may apply standard arguments from
  the Calculus of Variations showing that the minimisation problem is
  well posed. Hence, there exists $u\in \sob{1}{p}_g(\Omega,\R^N)$
  such that
  \begin{equation}
    \label{eq:min}
    \E_p(u,\Omega) \,=\, \min_{ v\in \sob{1}{p}_0(\Omega,\R^N) } \E_p(v,\Omega).
  \end{equation}
  Since $[\fes]^N \subset \sob{1}{p}(\Omega,\R^N)$ the same
  argument applies for the Galerkin approximation.

  Further, for fixed $p$, let $\{U_p\}$ be the finite element
  approximation generated by solving (\ref{eq:plapdis}) (indexed by
  $h$) and $u_p$, the weak solution of the $p$-Laplace system, then we
  have that
  \begin{equation}
    U_p \larrow  u_p  \ \text{ in }  C^{0}(\overline{\Omega},\R^N), \text{ as } h \ri  0.
  \end{equation}

\smallskip

\subsection{Numerical experiments}

Our numerical approximations are achieved using Galerkin
approximations to the $p$-Laplacian for sufficiently high values of
$p$.  We focus on studying the behaviour solutions have as $p$
increases which allow us to make various conjectures on the behaviour
of their asymptotic limit as $p\ri\infty$.

The computation of $p$-harmonic mappings is an extremely challenging
problem in its own right. The class of nonlinearity in the problem
results in the algebraic system, which ultimately yields the finite
element solution, being extremely badly conditioned. One method to
tackle this class of problems is by using preconditioners based on
descent algorithms \cite{HuangLiLiu:2007}. For extremely large $p$,
say $p \geq 10000$, this may be required; however for our purposes
we restrict our attention to $p \sim 100$. This yields sufficient
accuracy for the results we are illustrating.

We emphasise that even the case $p\sim 100$ is computationally difficult
to handle. The numerical approximation we are using is based on a
damped Newton solver. As it is well known, Newton solvers require a
sufficiently close initial guess in order to converge. A reasonable
initial guess for the $p$-Laplacian is given by numerically
approximating with the $q$-Laplacian for $q < p$ sufficiently close to
$p$. This leads to an iterative process in the generation of the
initial guess, i.e., we solve the $2$-Laplacian as an initial guess to
the $3$-Laplacian which serves as an initial guess to the
$4$-Laplacian, and so on. In each of our experiments the number of
nonlinear iterations required to achieve a relative tolerance of
$10^{-8}$ was achieved in less that $20$ iterations.

Using the methodology of $p$-approximation which we
advocate here, it has been numerically demonstrated in the scalar
case the rates of convergence both in $p$ and in $h$ that we
expect to achieve \cite[Section 4]{Pry}. It was noticed that these
rates were dependant on the regularity of the underlying
$\infty$-harmonic function and that
\begin{equation}
  \| u_\infty - U_{p^*} \|_{\leb{\infty}} \approx O(h)
\end{equation}
for solutions $u_\infty \in C^{\infty}$ and
\begin{equation}
  \| u_\infty - U_{p^*}\|_{\leb{\infty}} = O(h^{1/3})
\end{equation}
if $u_\infty \in C^{1,1/3}$. Where we use $p^*$ as the $\text{argmin}_p
\|u_\infty - U_p\|_{\leb{\infty}}$.
  
In both cases as $h$ is decreased, an increasing value of $p$ is
required to achieve optimal approximation (in $h$). This suggests a
coupling $p = C h^\alpha$ is necessary to achieve convergence, where
the $\alpha$ is determined by the regularity expected in
$u_\infty$. We found experimentally that coupling $p = h^{-1/2}$
worked well for the singular case ($u_\infty \in C^{1,1/3}$) and $p =
h^{-1}$ for the smooth case ($u_\infty \in C^{\infty}$).

We mention that we do \emph{not} have access to exact solutions in
general. For the experiments in Sections
\ref{sec:aronsson}--\ref{sec:vec} we only provide boundary data
and in principle the solutions to this problem are non-smooth and must
be interpreted in an appropriate weak sense. For the scalar problem
this is through the viscosity solution framework \cite{CIL}.

\providecommand{\figwidth}{\textwidth}
\newcommand{\figscale}{1.1}

\ms

\subsection{Scalar case - Aronsson solution}
\label{sec:aronsson}

\begin{figure}[h!]
  \caption[] {\label{fig:t1} A test to characterise the
    $\sigma_\infty^\mO$ for boundary data given by the Aronsson
    solution in (\ref{eq:aronsson}). We show the domain $\mO$ and the
    measure $\sigma_p^\mO(\mathcal T)$ for increasing values of
    $p$. In each of the subfigures B-F we plot 10 contours that are
    equally spaced between the minimum and maximum values of
    $\sigma_p^\mO$. Notice that as $p$ increases these contours tend
    toward the top right corner of $\mO$ indicating that mass is
    concentrated there. 
  }
  \begin{center}
    \subfloat[{}]
             {
               \hspace{.2cm}
               \begin{tikzpicture}[scale=2]
                 \newcommand\Square[1]{+(-#1,-#1) rectangle +(#1,#1)}
                 \draw [very thin, gray] (-1,-1) grid (1,1);  
                 %\draw [red] (2,3) +(-2pt,-2pt) rectangle +(2pt,2pt) ;
                 \draw [fill=red] (.5,.5) \Square{.375};
                 \node at (0.5,0.5) {$\mO$};
                 \draw [fill=blue] (.875,.875) circle (.025);
                 \node [above] at (.875,.875) {$x_0$};
                 \begin{scope}[thick,font=\scriptsize]
                   \draw [->] (-1.25,0) -- (1.25,0) node [above left]  {$x$};
                   \draw [->] (0,-1.25) -- (0,1.25) node [below right] {$y$};
                   \foreach \n in {-1,1}{%
                     \draw (\n,-3pt) -- (\n,3pt)   node [above] {$\n$};
                     \draw (-3pt,\n) -- (3pt,\n)   node [right] {$\n$};
                   }
                 \end{scope} 
                 %\path [draw=none,fill=gray,semitransparent] (+1,-1) circle (3);
               \end{tikzpicture}
             }
     \hfill
     \subfloat[{
        {
          $p=2$
        }
    }]{
      \includegraphics[scale=\figscale,width=0.45\figwidth]{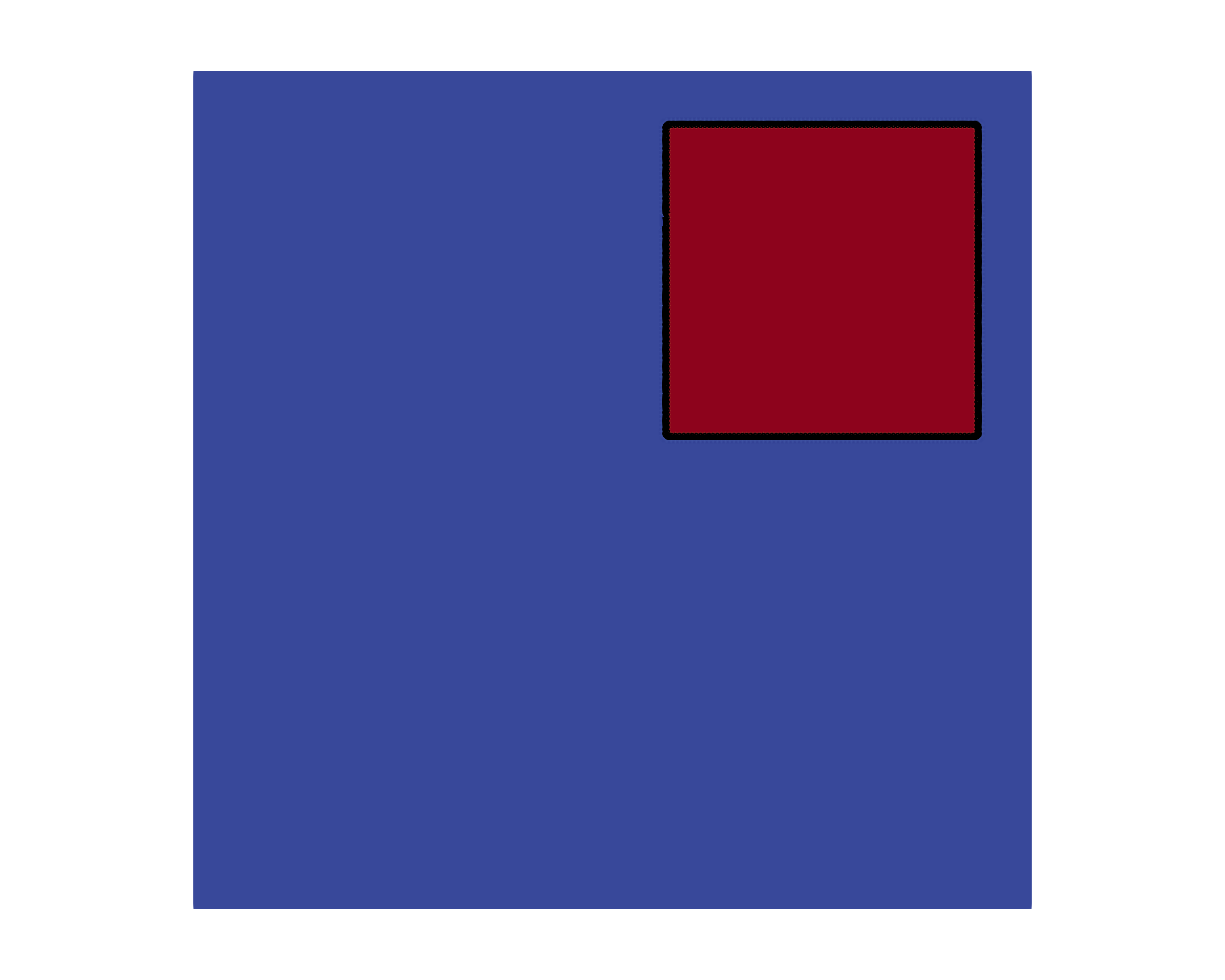}
    }
    \\
    \subfloat[{
        {
          $p=4$
        }  
    }]{
      \includegraphics[scale=\figscale,width=0.45\figwidth]{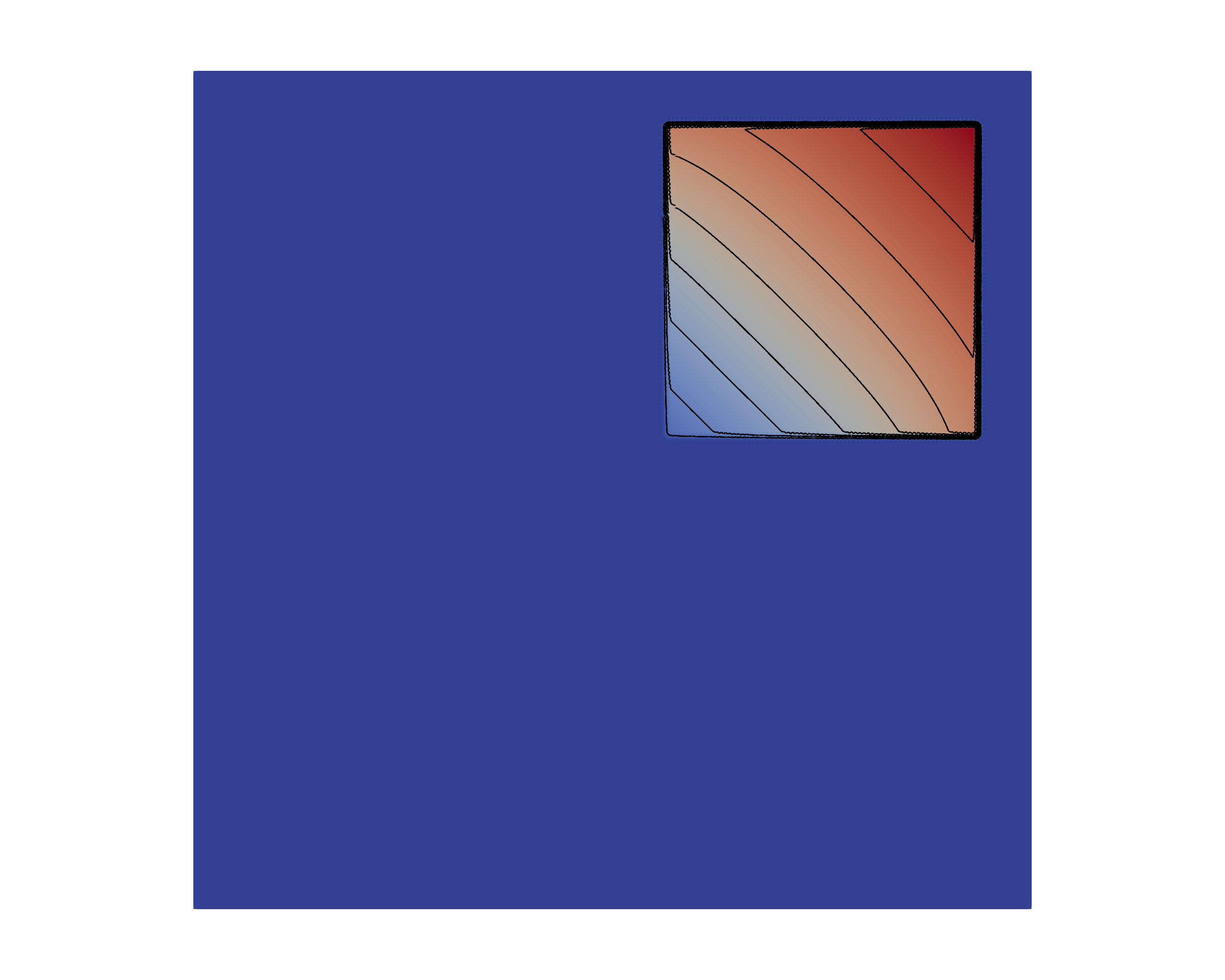}
    }
     \hfill
    \subfloat[{
        %%%%%%%%%%%%%%%%%%%%%%%%%%%%%%%%%%%%%%%%%%%%%%%%%%%%%%%%%%%%%%% 
        {
          $p=10$
        }
        %%%%%%%%%%%%%%%%%%%%%%%%%%%%%%%%%%%%%%%%%%%%%%%%%%%%%%%%%%%%%%%% 
    }]{
      \includegraphics[scale=\figscale,width=0.45\figwidth]{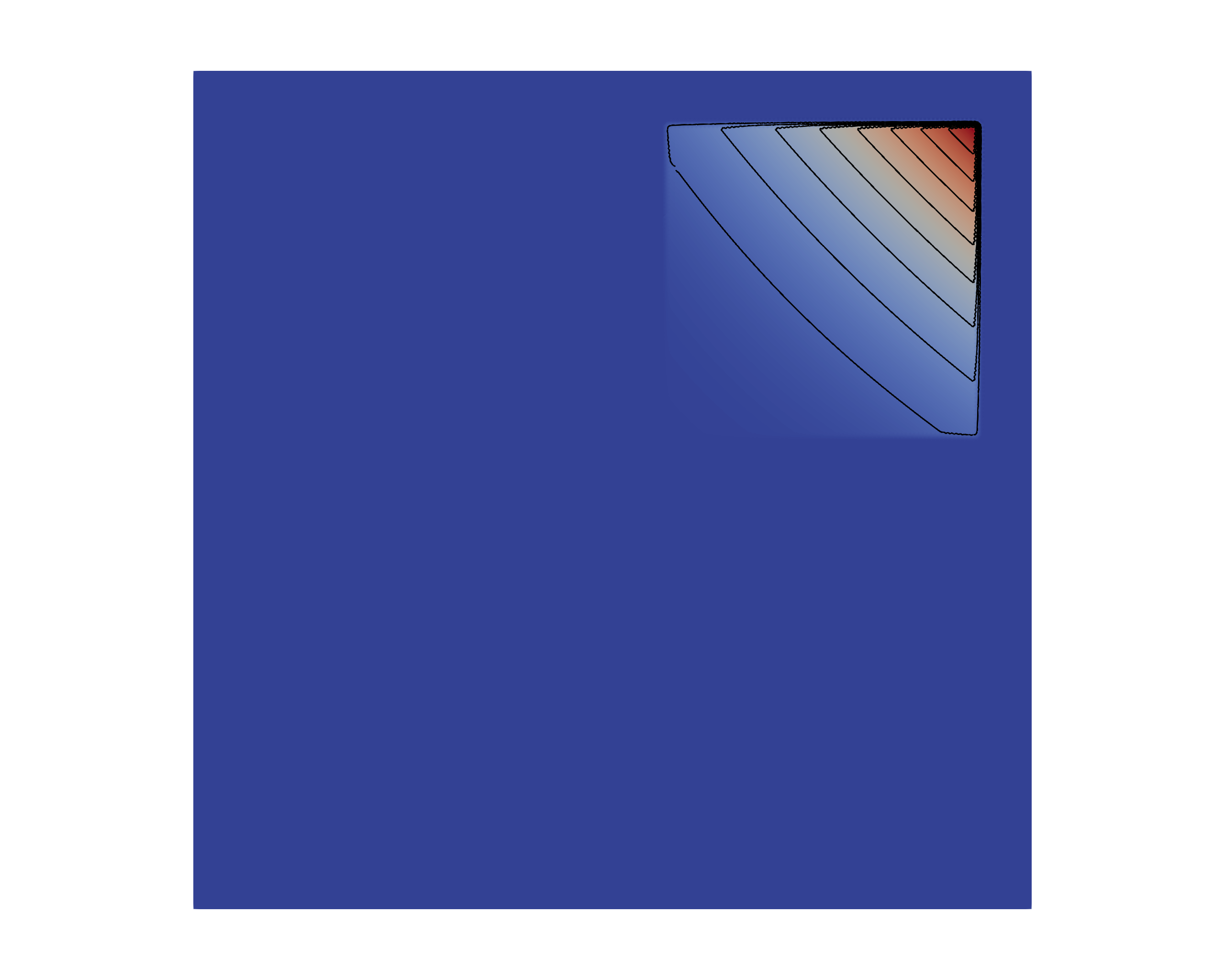}
        }
        \\
      \subfloat[{
      %%%%%%%%%%%%%%%%%%%%%%%%%%%%%%%%%%%%%%%%%%%%%%%%%%%%%%%%%%%%%%% 
        {
          $p=20$
        }
        %%%%%%%%%%%%%%%%%%%%%%%%%%%%%%%%%%%%%%%%%%%%%%%%%%%%%%%%%%%%%%%% 
    }]{
      \includegraphics[scale=\figscale,width=0.45\figwidth]{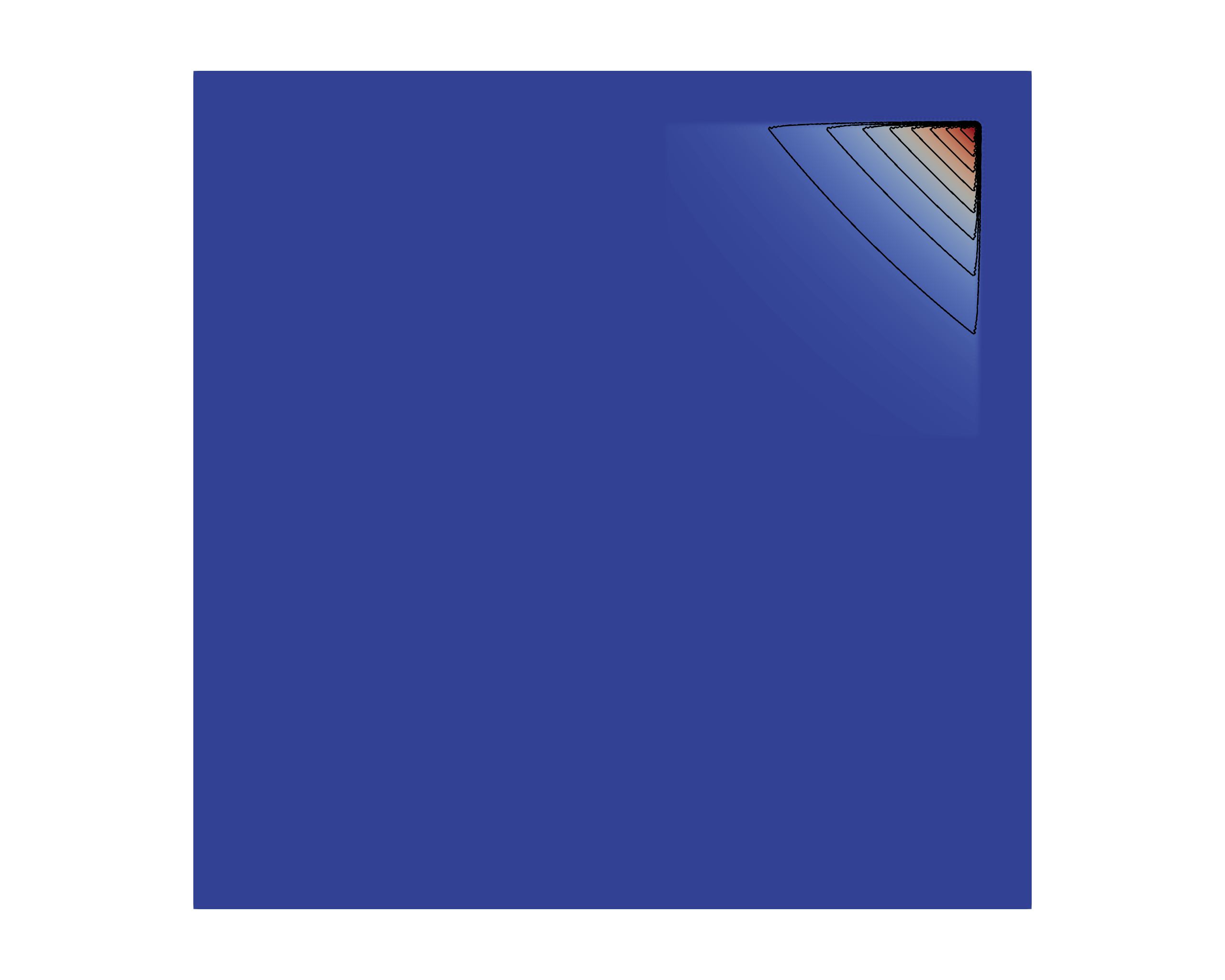}
    }
    %%%%%%%%%%%%%%%%%%%%%%%%%%%%%%%%%%%%%%%%%%%%%%%%%%%%%%%%%%%%%%%%%%%% 
     \hfill
     \subfloat[{
      %%%%%%%%%%%%%%%%%%%%%%%%%%%%%%%%%%%%%%%%%%%%%%%%%%%%%%%%%%%%%%% 
        {
          $p=100$
        }
        %%%%%%%%%%%%%%%%%%%%%%%%%%%%%%%%%%%%%%%%%%%%%%%%%%%%%%%%%%%%%%%% 
    }]{
      \includegraphics[scale=\figscale,width=0.45\figwidth]{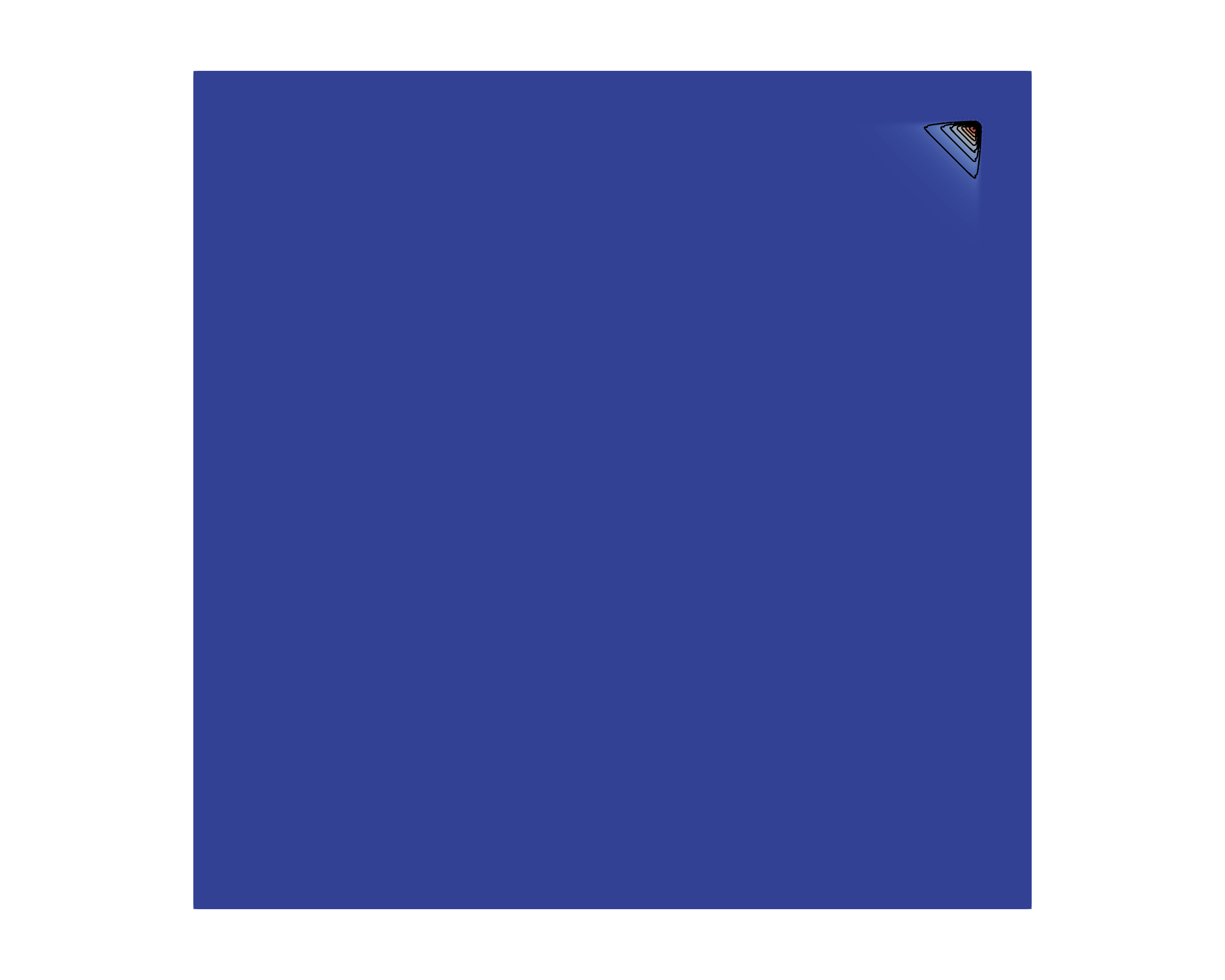}
    }
    %%%%%%%%%%%%%%%%%%%%%%%%%%%%%%%%%%%%%%%%%%%%%%%%%%%%%%%%%%%%%%%%%%%% 
  \end{center}
\end{figure}
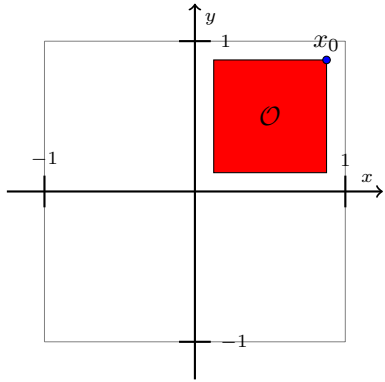

\begin{figure}[h!]
  \caption[] {\label{fig:t2} As Figure \ref{fig:t1} with a different
    $\mO$. Notice that as $p$ increases these contours tend
    toward the corners of $\mO$ indicating that mass is concentrated
    there. We conjecture that $\sigma_\infty^\mO = \frac{1}{4}\sum_{j=0}^3
    \delta_{x_j}$, with $x_j$ denoting the four corners of $\mO$.}
  \begin{center}
    \subfloat[{The domain $\mO$ and the level sets of $|\D u_\infty|$.}]
             {
               \hspace{.2cm}
    \begin{tikzpicture}[scale=2]
      \newcommand\Square[1]{+(-#1,-#1) rectangle +(#1,#1)}
      \draw [very thin, gray] (-1,-1) grid (1,1);  
      \draw plot[smooth, samples=50, domain=0:1]
      (\x, {max((1 - 3*\x^(2/3) + 3*\x^(4/3) - \x^2),0)^0.5}) -| (0,0) -- cycle;
      \draw plot[smooth, samples=50, domain=0:1]
      (-\x, {max((1 - 3*\x^(2/3) + 3*\x^(4/3) - \x^2),0)^0.5}) -| (0,0) -- cycle;
      \draw plot[smooth, samples=50, domain=0:1]
      (-\x, {-max((1 - 3*\x^(2/3) + 3*\x^(4/3) - \x^2),0)^0.5}) -| (0,0) -- cycle;
      \draw plot[smooth, samples=50, domain=0:1]
      (\x, {-max((1 - 3*\x^(2/3) + 3*\x^(4/3) - \x^2),0)^0.5}) -| (0,0) -- cycle;
      \draw [fill=red] (.0,.0) \Square{.25};
      \draw [fill=blue] (-.25,-.25) circle (.025);
      \node [left] at (-.2,-.25) {{$x_0$}};
      \draw [fill=blue] (.25,-.25) circle (.025);
      \node [right] at (.2,-.25) {{$x_1$}};
      \draw [fill=blue] (-.25,.25) circle (.025);
      \node [left] at (-.2,.25) {{$x_2$}};
      \draw [fill=blue] (.25,.25) circle (.025);
      \node [right] at (.2,.25) {{$x_3$}};
      \begin{scope}[thick,font=\scriptsize]
       \draw [->] (-1.25,0) -- (1.25,0) node [above left]  {$x$};
       \draw [->] (0,-1.25) -- (0,1.25) node [below right] {$y$};
       \foreach \n in {-1,1}{%
         \draw (\n,-3pt) -- (\n,3pt)   node [above] {$\n$};
         \draw (-3pt,\n) -- (3pt,\n)   node [right] {$\n$};
       }
          \end{scope}
    \end{tikzpicture}
             }
     \hfill
     \subfloat[{
        {
          $p=2$
        }
    }]{
      \includegraphics[scale=\figscale,width=0.45\figwidth]{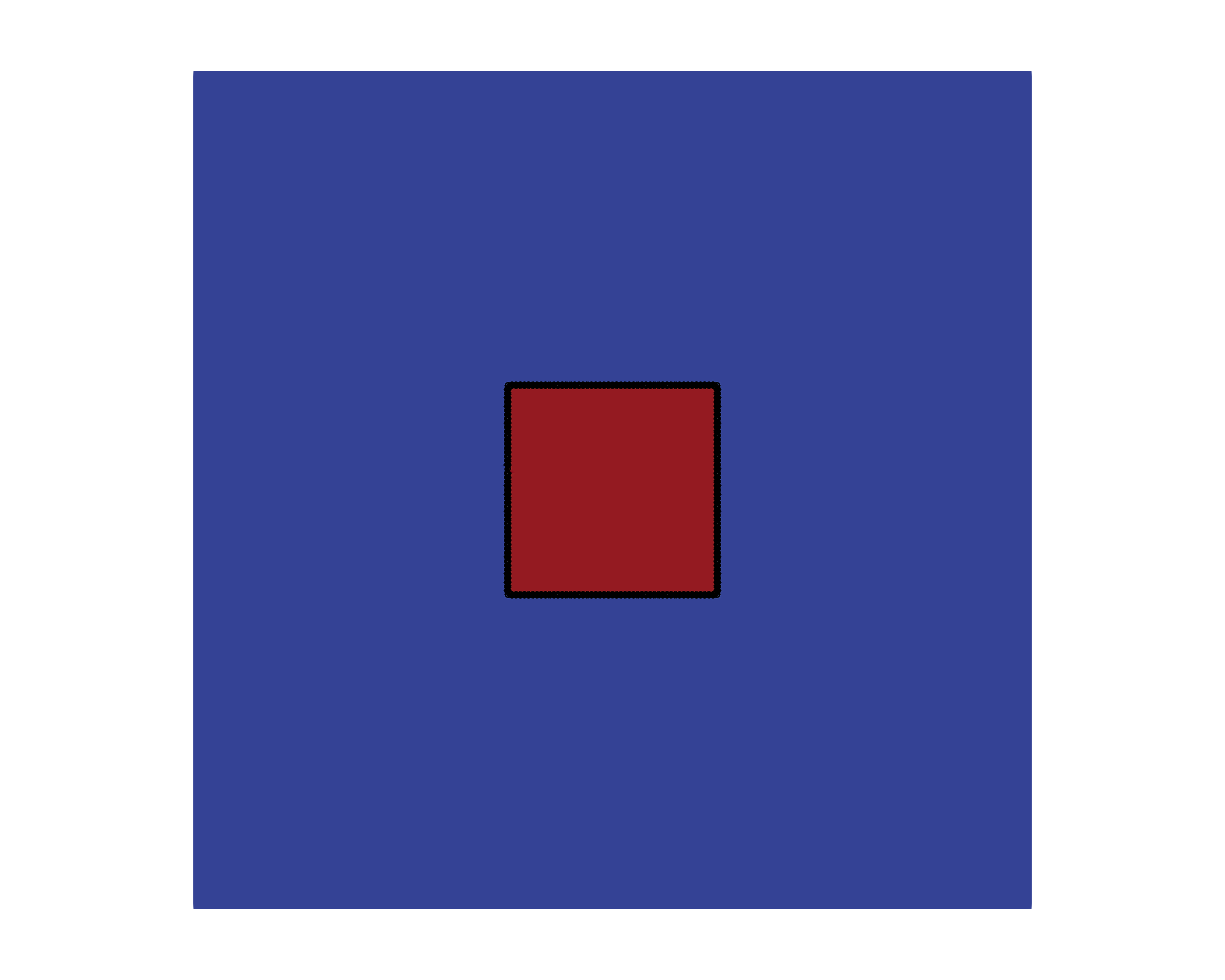}
    }
    \\
    \subfloat[{
        {
          $p=4$
        }  
    }]{
      \includegraphics[scale=\figscale,width=0.45\figwidth]{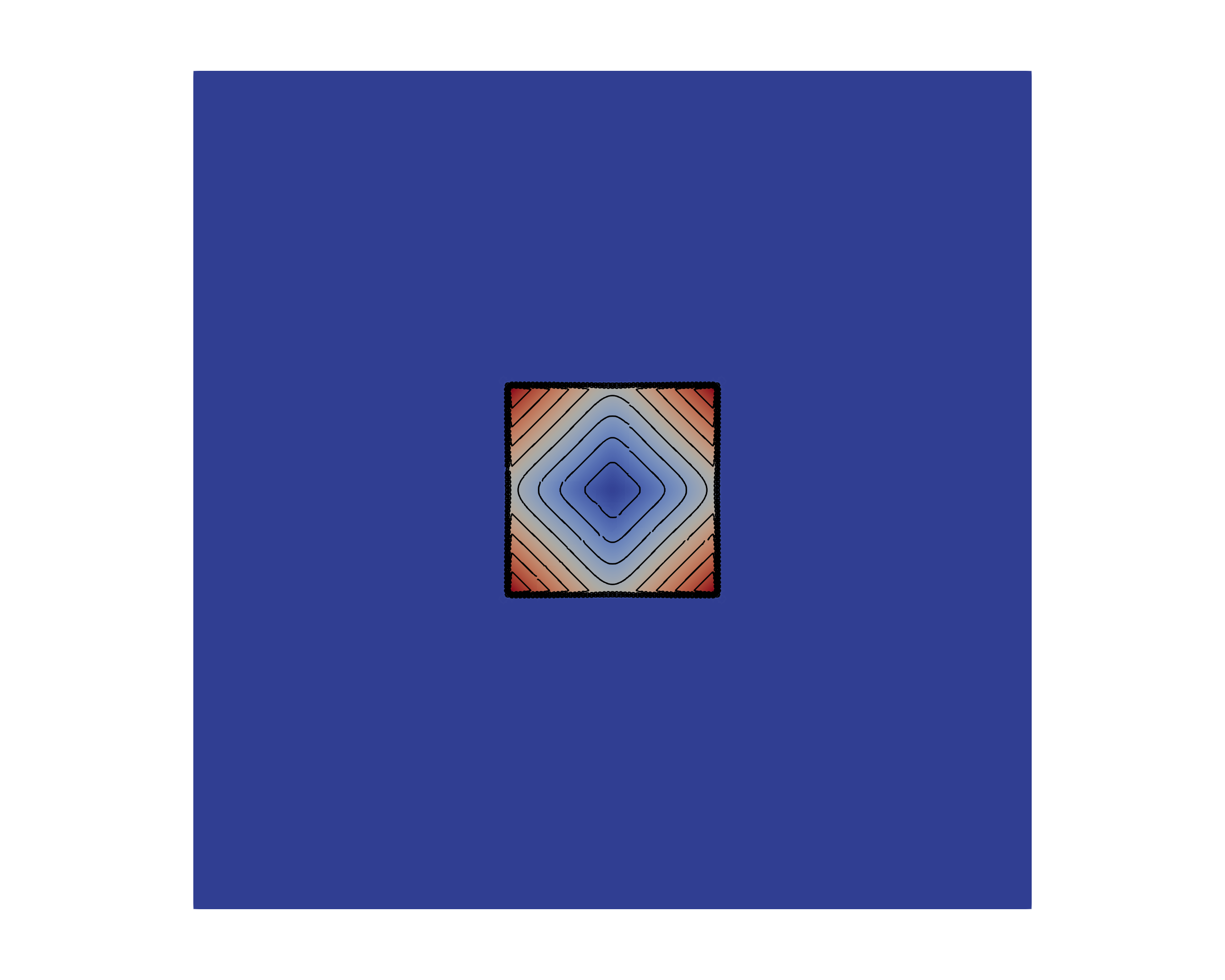}
    }
     \hfill
    \subfloat[{
        %%%%%%%%%%%%%%%%%%%%%%%%%%%%%%%%%%%%%%%%%%%%%%%%%%%%%%%%%%%%%%% 
        {
          $p=10$
        }
        %%%%%%%%%%%%%%%%%%%%%%%%%%%%%%%%%%%%%%%%%%%%%%%%%%%%%%%%%%%%%%%% 
    }]{
      \includegraphics[scale=\figscale,width=0.45\figwidth]{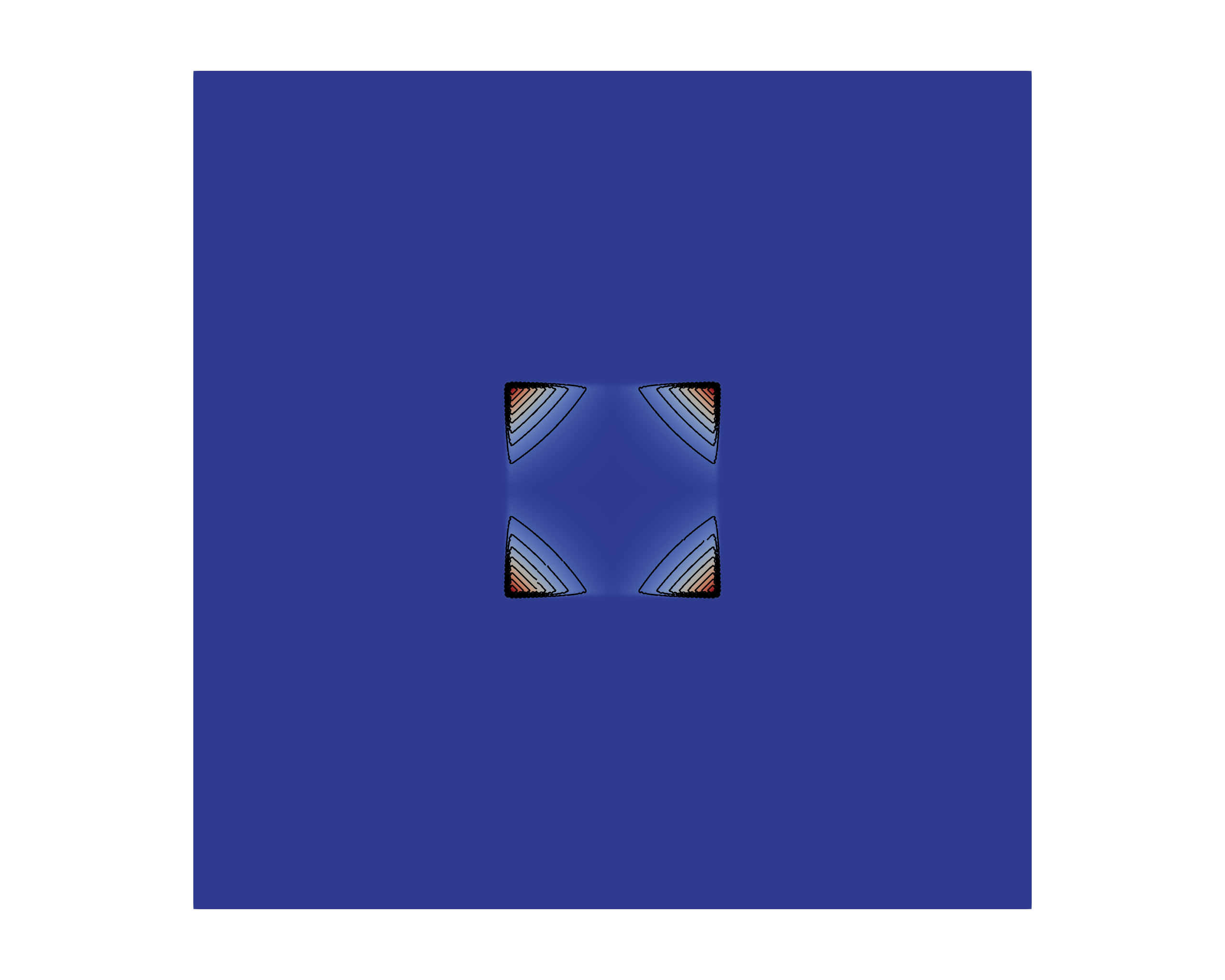}
        }
        \\
      \subfloat[{
      %%%%%%%%%%%%%%%%%%%%%%%%%%%%%%%%%%%%%%%%%%%%%%%%%%%%%%%%%%%%%%% 
        {
          $p=20$
        }
        %%%%%%%%%%%%%%%%%%%%%%%%%%%%%%%%%%%%%%%%%%%%%%%%%%%%%%%%%%%%%%%% 
    }]{
      \includegraphics[scale=\figscale,width=0.45\figwidth]{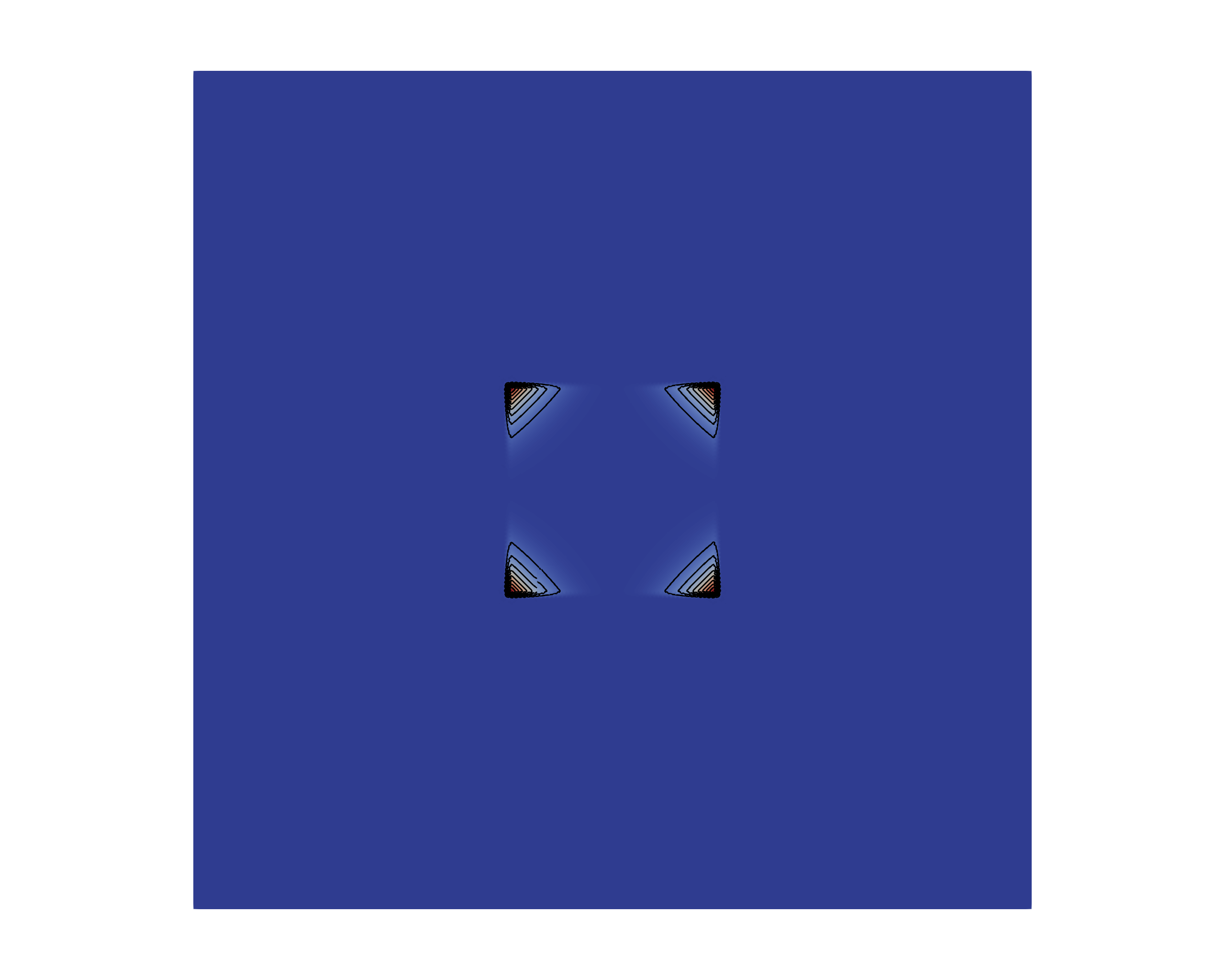}
    }
    %%%%%%%%%%%%%%%%%%%%%%%%%%%%%%%%%%%%%%%%%%%%%%%%%%%%%%%%%%%%%%%%%%%% 
     \hfill
     \subfloat[{
      %%%%%%%%%%%%%%%%%%%%%%%%%%%%%%%%%%%%%%%%%%%%%%%%%%%%%%%%%%%%%%% 
        {
          $p=100$
        }
        %%%%%%%%%%%%%%%%%%%%%%%%%%%%%%%%%%%%%%%%%%%%%%%%%%%%%%%%%%%%%%%% 
    }]{
      \includegraphics[scale=\figscale,width=0.45\figwidth]{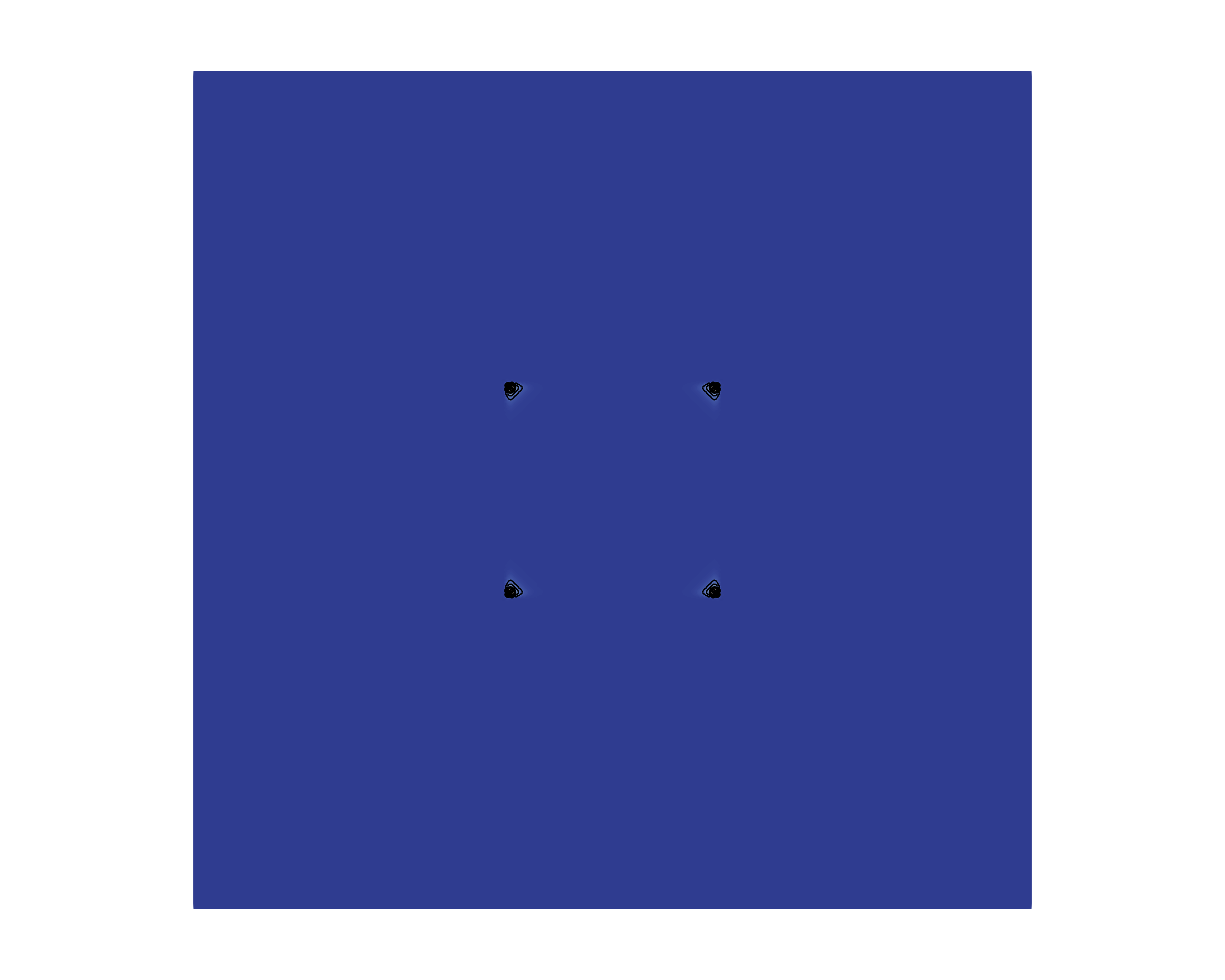}
    }
    %%%%%%%%%%%%%%%%%%%%%%%%%%%%%%%%%%%%%%%%%%%%%%%%%%%%%%%%%%%%%%%%%%%% 
  \end{center}
\end{figure}
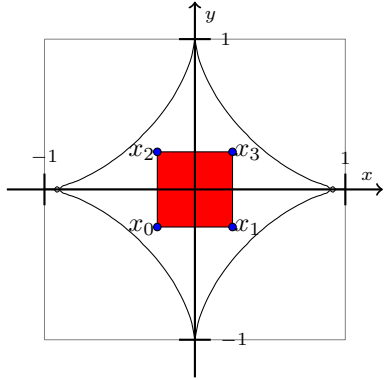

\begin{figure}[h!]
  \caption[] {\label{fig:t3} As Figure \ref{fig:t1} with a different
    $\mO$. The corners of the square defining $\mO$ are touching the
    level sets of $|\D u_\infty|$. Again, we conjecture that
    $\sigma_\infty^\mO = \frac{1}{4} \sum_{j=0}^3 \delta_{x_j}$, with $x_j$
    denoting the four corners of $\mO$.}
  \begin{center}
    \subfloat[{The domain $\mO$ and the level sets of $| \D u_\infty |$.}]
             {
               \hspace{.2cm}
    \begin{tikzpicture}[scale=2]
      \newcommand\Square[1]{+(-#1,-#1) rectangle +(#1,#1)}
      \draw [very thin, gray] (-1,-1) grid (1,1);  
      \draw plot[smooth, samples=50, domain=0:1]
      (\x, {max((1 - 3*\x^(2/3) + 3*\x^(4/3) - \x^2),0)^0.5}) -| (0,0) -- cycle;
      \draw plot[smooth, samples=50, domain=0:1]
      (-\x, {max((1 - 3*\x^(2/3) + 3*\x^(4/3) - \x^2),0)^0.5}) -| (0,0) -- cycle;
      \draw plot[smooth, samples=50, domain=0:1]
      (-\x, {-max((1 - 3*\x^(2/3) + 3*\x^(4/3) - \x^2),0)^0.5}) -| (0,0) -- cycle;
      \draw plot[smooth, samples=50, domain=0:1]
      (\x, {-max((1 - 3*\x^(2/3) + 3*\x^(4/3) - \x^2),0)^0.5}) -| (0,0) -- cycle;
      \draw [fill=red] (0.,0.) \Square{.35};
      \draw [fill=blue] (-.35,-.35) circle (.025);
      \node [left] at (-.35,-.35) {{$x_0$}};
      \draw [fill=blue] (.35,-.35) circle (.025);
      \node [right] at (.35,-.35) {{$x_1$}};
      \draw [fill=blue] (-.35,.35) circle (.025);
      \node [left] at (-.35,.35) {{$x_2$}};
      \draw [fill=blue] (.35,.35) circle (.025);
      \node [right] at (.35,.35) {{$x_3$}};
      \begin{scope}[thick,font=\scriptsize]
       \draw [->] (-1.25,0) -- (1.25,0) node [above left]  {$x$};
       \draw [->] (0,-1.25) -- (0,1.25) node [below right] {$y$};
       \foreach \n in {-1,1}{%
         \draw (\n,-3pt) -- (\n,3pt)   node [above] {$\n$};
         \draw (-3pt,\n) -- (3pt,\n)   node [right] {$\n$};
       }
      \end{scope}
    \end{tikzpicture}
             }
     \hfill
     \subfloat[{
        {
          $p=2$
        }
    }]{
      \includegraphics[scale=\figscale,width=0.45\figwidth]{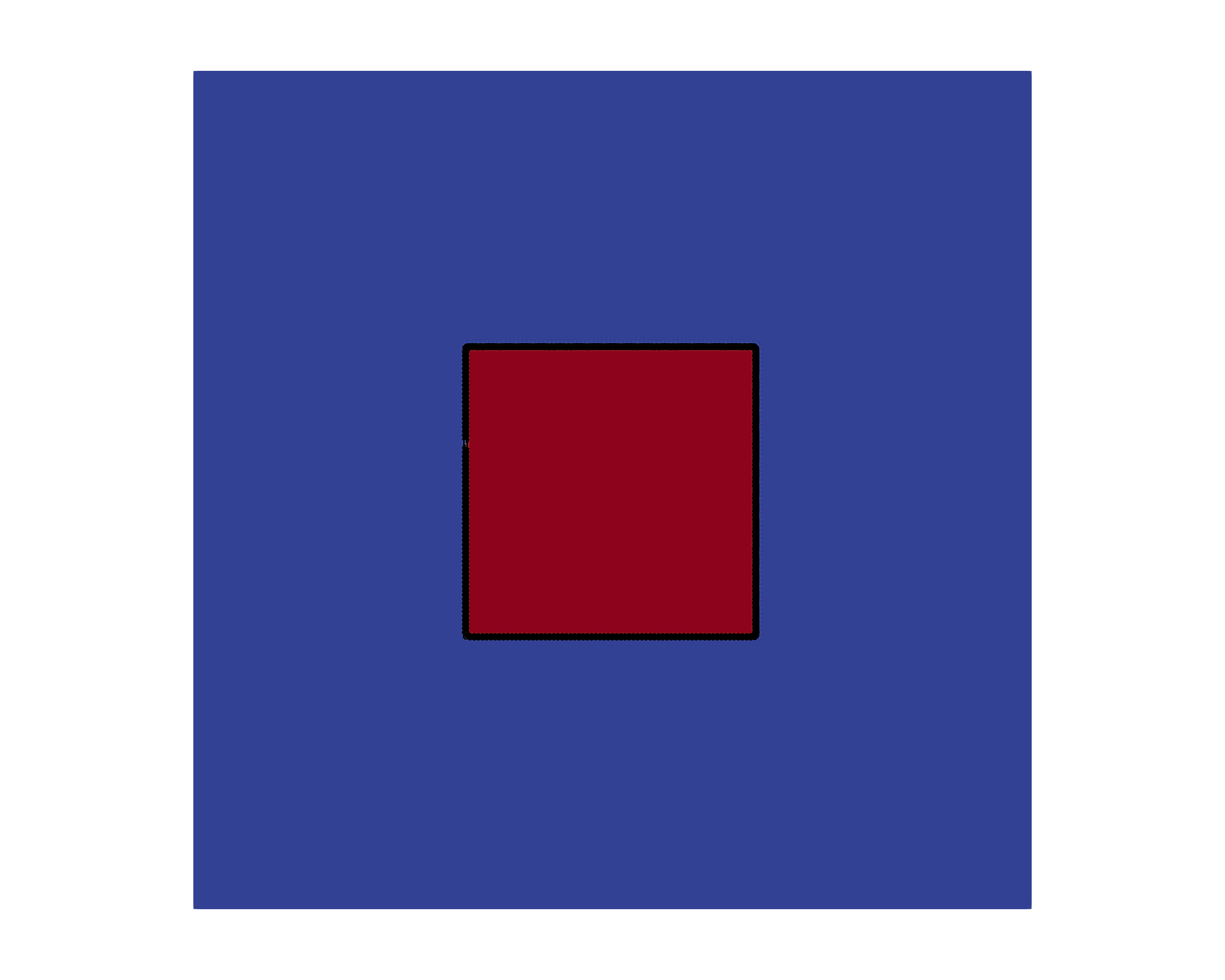}
    }
    \\
    \subfloat[{
        {
          $p=4$
        }  
    }]{
      \includegraphics[scale=\figscale,width=0.45\figwidth]{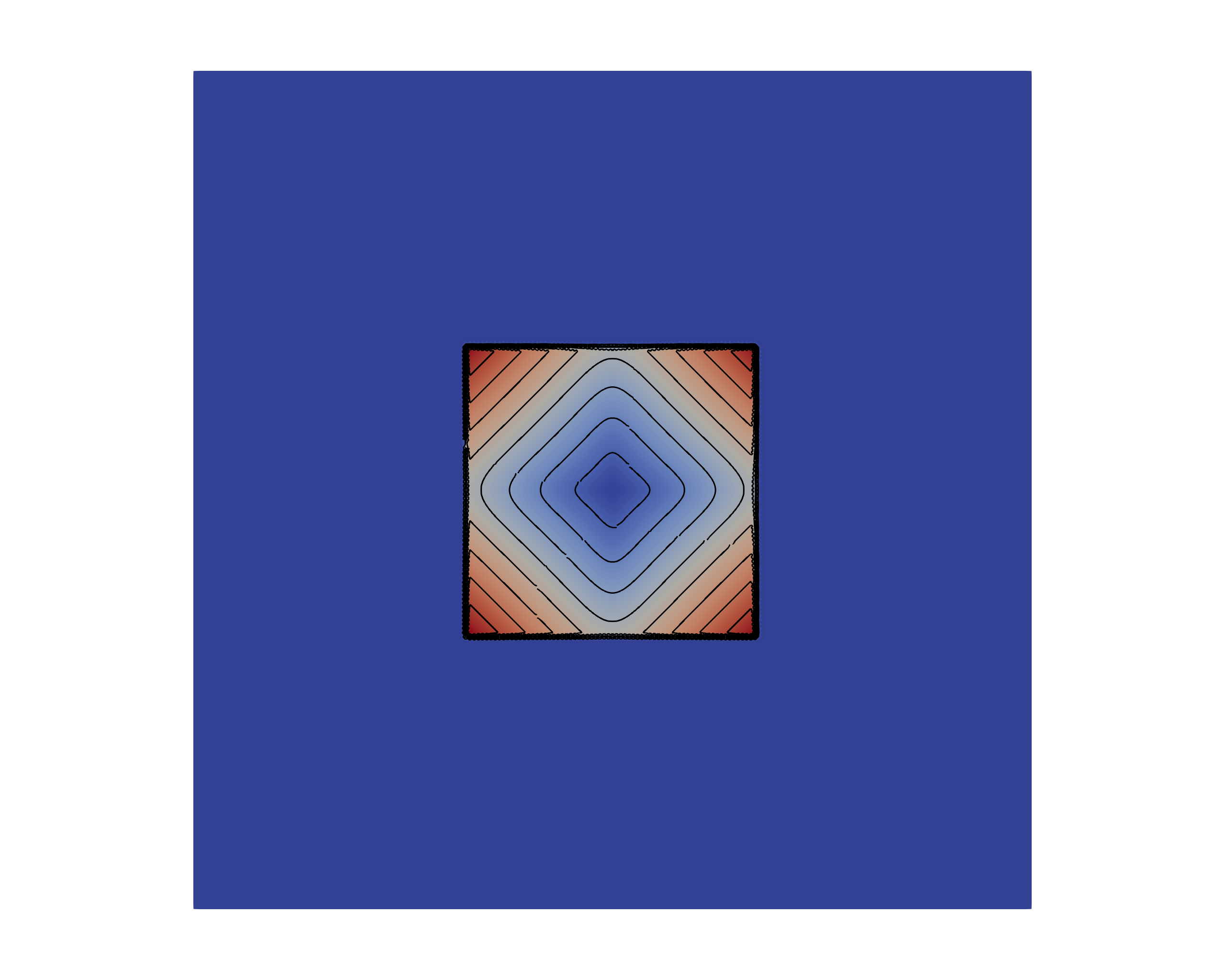}
    }
     \hfill
    \subfloat[{
        %%%%%%%%%%%%%%%%%%%%%%%%%%%%%%%%%%%%%%%%%%%%%%%%%%%%%%%%%%%%%%% 
        {
          $p=10$
        }
        %%%%%%%%%%%%%%%%%%%%%%%%%%%%%%%%%%%%%%%%%%%%%%%%%%%%%%%%%%%%%%%% 
    }]{
      \includegraphics[scale=\figscale,width=0.45\figwidth]{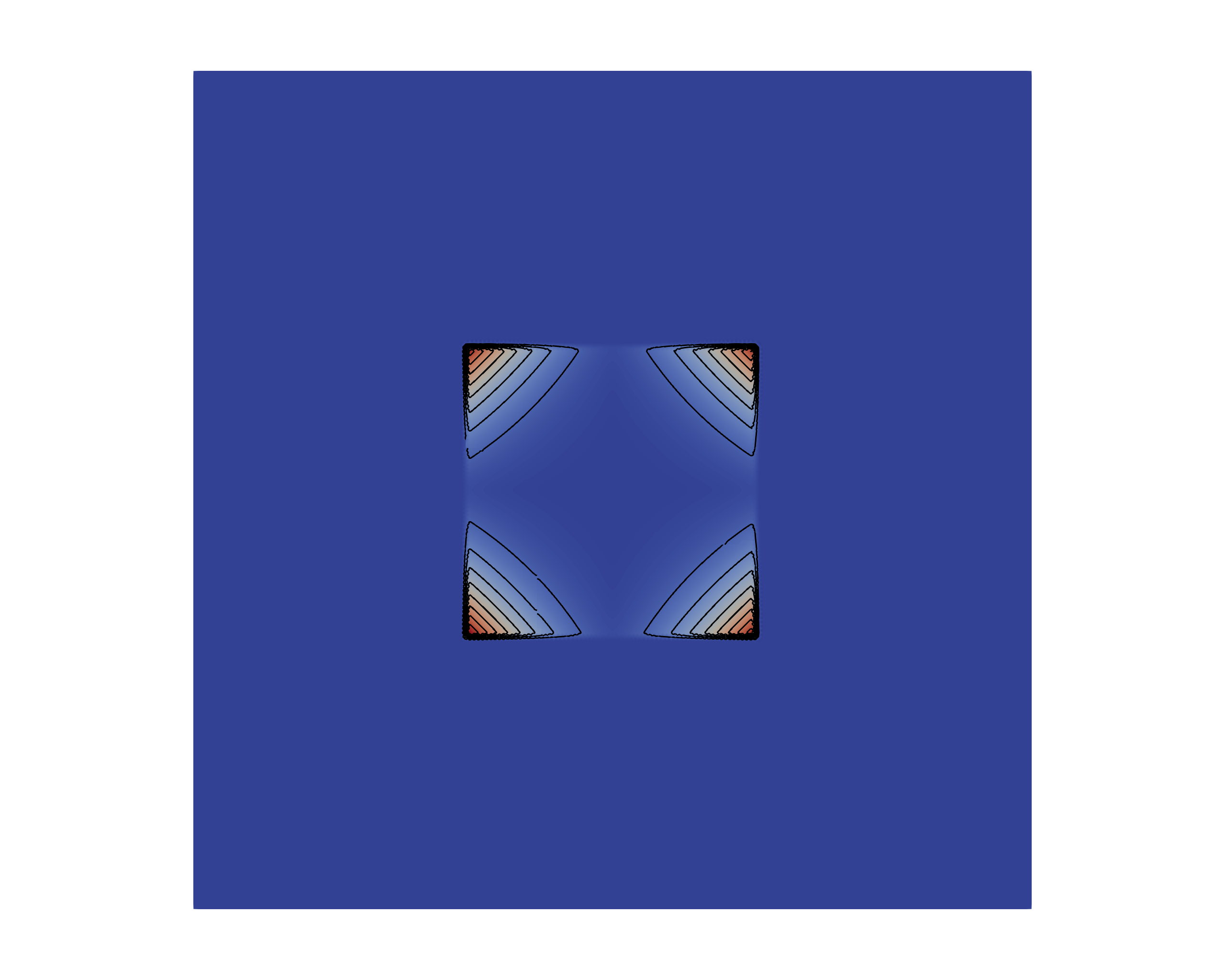}
        }
        \\
      \subfloat[{
      %%%%%%%%%%%%%%%%%%%%%%%%%%%%%%%%%%%%%%%%%%%%%%%%%%%%%%%%%%%%%%% 
        {
          $p=20$
        }
        %%%%%%%%%%%%%%%%%%%%%%%%%%%%%%%%%%%%%%%%%%%%%%%%%%%%%%%%%%%%%%%% 
    }]{
      \includegraphics[scale=\figscale,width=0.45\figwidth]{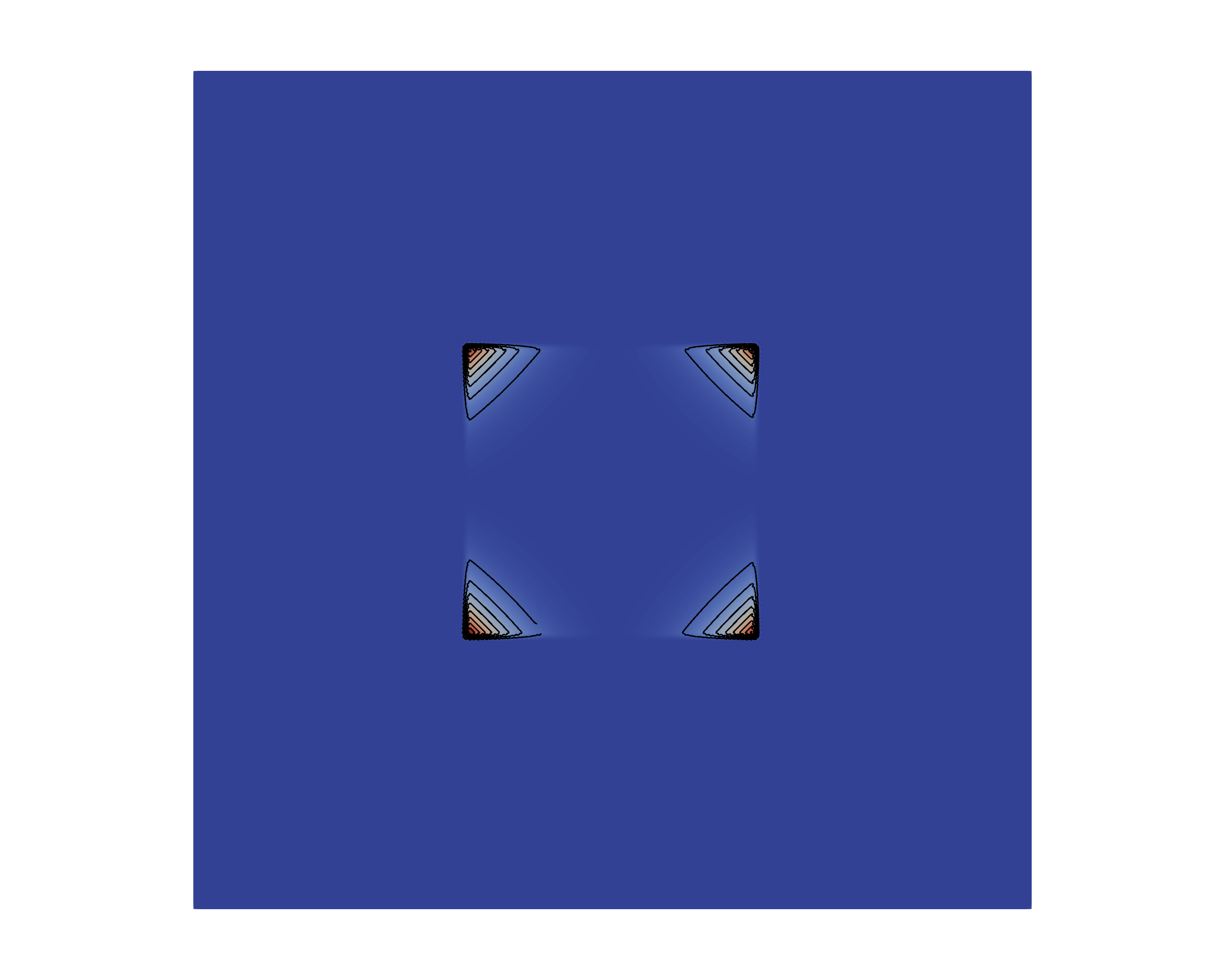}
    }
    %%%%%%%%%%%%%%%%%%%%%%%%%%%%%%%%%%%%%%%%%%%%%%%%%%%%%%%%%%%%%%%%%%%% 
     \hfill
     \subfloat[{
      %%%%%%%%%%%%%%%%%%%%%%%%%%%%%%%%%%%%%%%%%%%%%%%%%%%%%%%%%%%%%%% 
        {
          $p=100$
        }
        %%%%%%%%%%%%%%%%%%%%%%%%%%%%%%%%%%%%%%%%%%%%%%%%%%%%%%%%%%%%%%%% 
    }]{
      \includegraphics[scale=\figscale,width=0.45\figwidth]{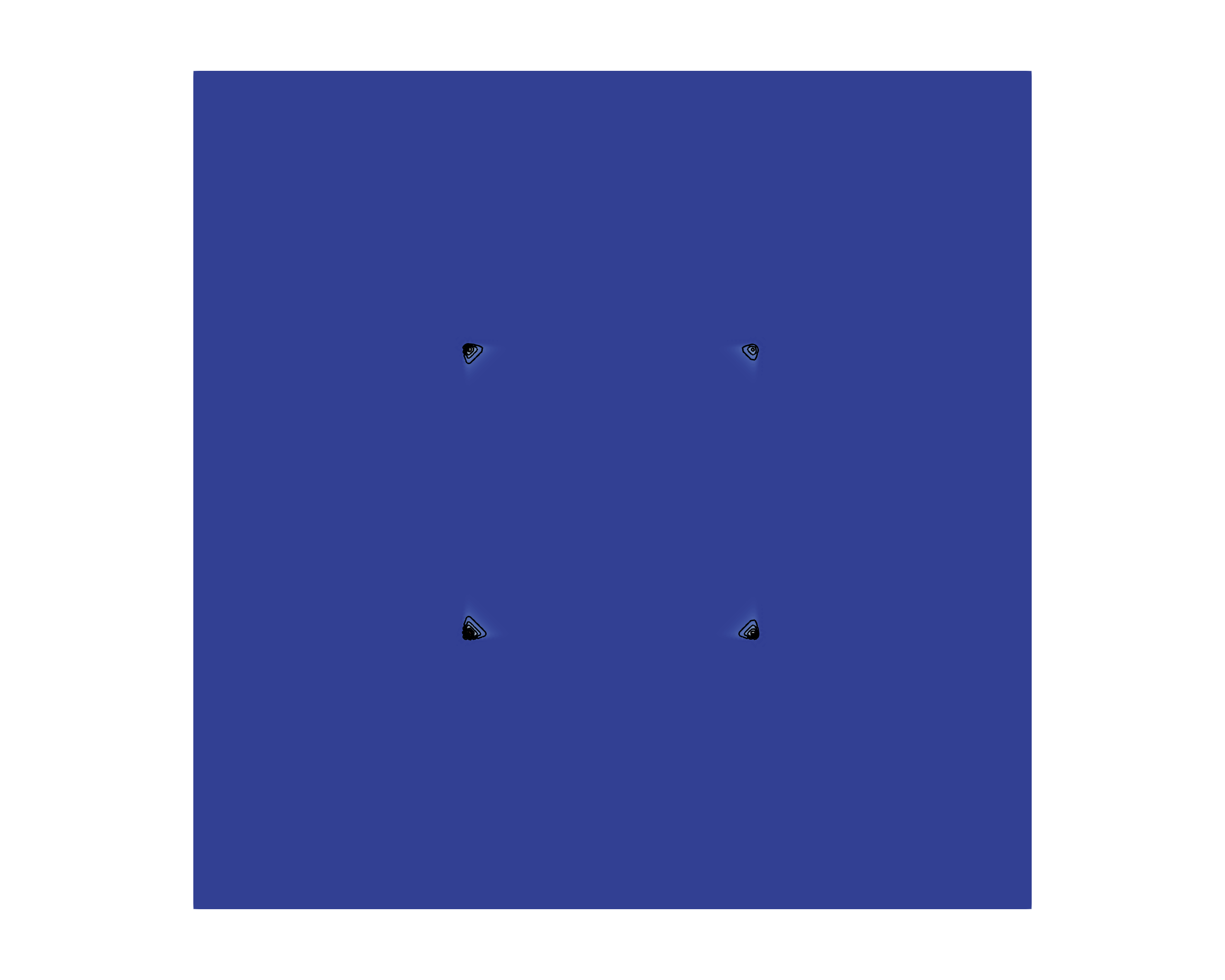}
    }
    %%%%%%%%%%%%%%%%%%%%%%%%%%%%%%%%%%%%%%%%%%%%%%%%%%%%%%%%%%%%%%%%%%%% 
  \end{center}
\end{figure}
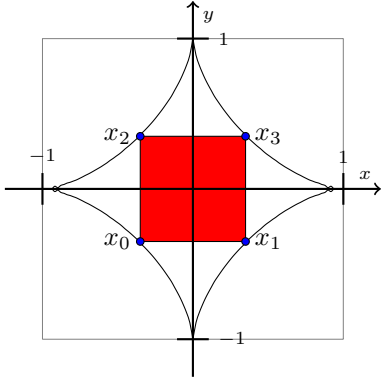

\begin{figure}[h!]
  \caption[] {\label{fig:t4} As Figure \ref{fig:t2} with a different
    $\mO$. The boundary of $\mO$ now contains a portion of the level sets
    of $|\D u_\infty|$. Notice that as $p$ increases, mass is concentrated on
    four portions of the boundary, where the red square intersects with the level sets. We conjecture that $\si^\mO=\frac{1}{C}\sum_{i=0}^3\mH^1\LL_{\Ga_i}$, where $\mH^1$ is the Hausdorff measure, $\Ga_0,...,\Ga_3$ are the four boundary portions and $C:=\sum_{i=0}^3\mH^1(\Ga_i)$.}
  \begin{center}
    \subfloat[{The domain $\mO$ and the level sets of $| \D u_\infty |$.}]
             {
               \hspace{.2cm}
    \begin{tikzpicture}[scale=2]
      \newcommand\Square[1]{+(-#1,-#1) rectangle +(#1,#1)}
      \draw [very thin, gray] (-1,-1) grid (1,1);  
      \draw plot[smooth, samples=50, domain=0:1]
      (\x, {max((1 - 3*\x^(2/3) + 3*\x^(4/3) - \x^2),0)^0.5}) -| (0,0) -- cycle;
      \draw plot[smooth, samples=50, domain=0:1]
      (-\x, {max((1 - 3*\x^(2/3) + 3*\x^(4/3) - \x^2),0)^0.5}) -| (0,0) -- cycle;
      \draw plot[smooth, samples=50, domain=0:1]
      (-\x, {-max((1 - 3*\x^(2/3) + 3*\x^(4/3) - \x^2),0)^0.5}) -| (0,0) -- cycle;
      \draw plot[smooth, samples=50, domain=0:1]
      (\x, {-max((1 - 3*\x^(2/3) + 3*\x^(4/3) - \x^2),0)^0.5}) -| (0,0) -- cycle;
      \draw [fill=red] plot[smooth, samples=50, domain=0.225:0.5]
      (\x, {max((1 - 3*\x^(2/3) + 3*\x^(4/3) - \x^2),0)^0.5}) -|  (0.5,0.) -- (0., 0.) -- (0.,.5) -- cycle;
      \draw [fill=red] plot[smooth, samples=50, domain=0.225:0.5]
      (-\x, {max((1 - 3*\x^(2/3) + 3*\x^(4/3) - \x^2),0)^0.5}) -|  (-0.5,0.) -- (0., 0.) -- (0.,.5) -- cycle;
      \draw [fill=red] plot[smooth, samples=50, domain=0.225:0.5]
      (\x, -{max((1 - 3*\x^(2/3) + 3*\x^(4/3) - \x^2),0)^0.5}) -|  (0.5,0.) -- (0., 0.) -- (0.,-.5) -- cycle;
      \draw [fill=red] plot[smooth, samples=50, domain=0.225:0.5]
      (-\x, -{max((1 - 3*\x^(2/3) + 3*\x^(4/3) - \x^2),0)^0.5}) -|  (-0.5,0.) -- (0., 0.) -- (0.,-.5) -- cycle;

      \draw [blue,ultra thick] plot[smooth, samples=50, domain=0.225:0.5]
      (\x, {max((1 - 3*\x^(2/3) + 3*\x^(4/3) - \x^2),0)^0.5});
      \draw [fill=blue] (.225,.5) circle (.025);
      \draw [fill=blue] (.5,.225) circle (.025);

      \draw [blue,ultra thick] plot[smooth, samples=50, domain=0.225:0.5]
      (-\x, {max((1 - 3*\x^(2/3) + 3*\x^(4/3) - \x^2),0)^0.5});
      \draw [fill=blue] (-.225,.5) circle (.025);
      \draw [fill=blue] (-.5,.225) circle (.025);

      \draw [blue,ultra thick] plot[smooth, samples=50, domain=0.225:0.5]
      (\x, -{max((1 - 3*\x^(2/3) + 3*\x^(4/3) - \x^2),0)^0.5});
      \draw [fill=blue] (.225,-.5) circle (.025);
      \draw [fill=blue] (.5,-.225) circle (.025);

      \draw [blue,ultra thick] plot[smooth, samples=50, domain=0.225:0.5]
      (-\x, -{max((1 - 3*\x^(2/3) + 3*\x^(4/3) - \x^2),0)^0.5});
      \draw [fill=blue] (-.225,-.5) circle (.025);
      \draw [fill=blue] (-.5,-.225) circle (.025);
      
      \begin{scope}[thick,font=\scriptsize]
       \draw [->] (-1.25,0) -- (1.25,0) node [above left]  {$x$};
       \draw [->] (0,-1.25) -- (0,1.25) node [below right] {$y$};
       \foreach \n in {-1,1}{%
         \draw (\n,-3pt) -- (\n,3pt)   node [above] {$\n$};
         \draw (-3pt,\n) -- (3pt,\n)   node [right] {$\n$};
       }
     \end{scope}
    \end{tikzpicture}
             }
     \hfill
     \subfloat[{
        {
          $p=2$
        }
    }]{
      \includegraphics[scale=\figscale,width=0.45\figwidth]{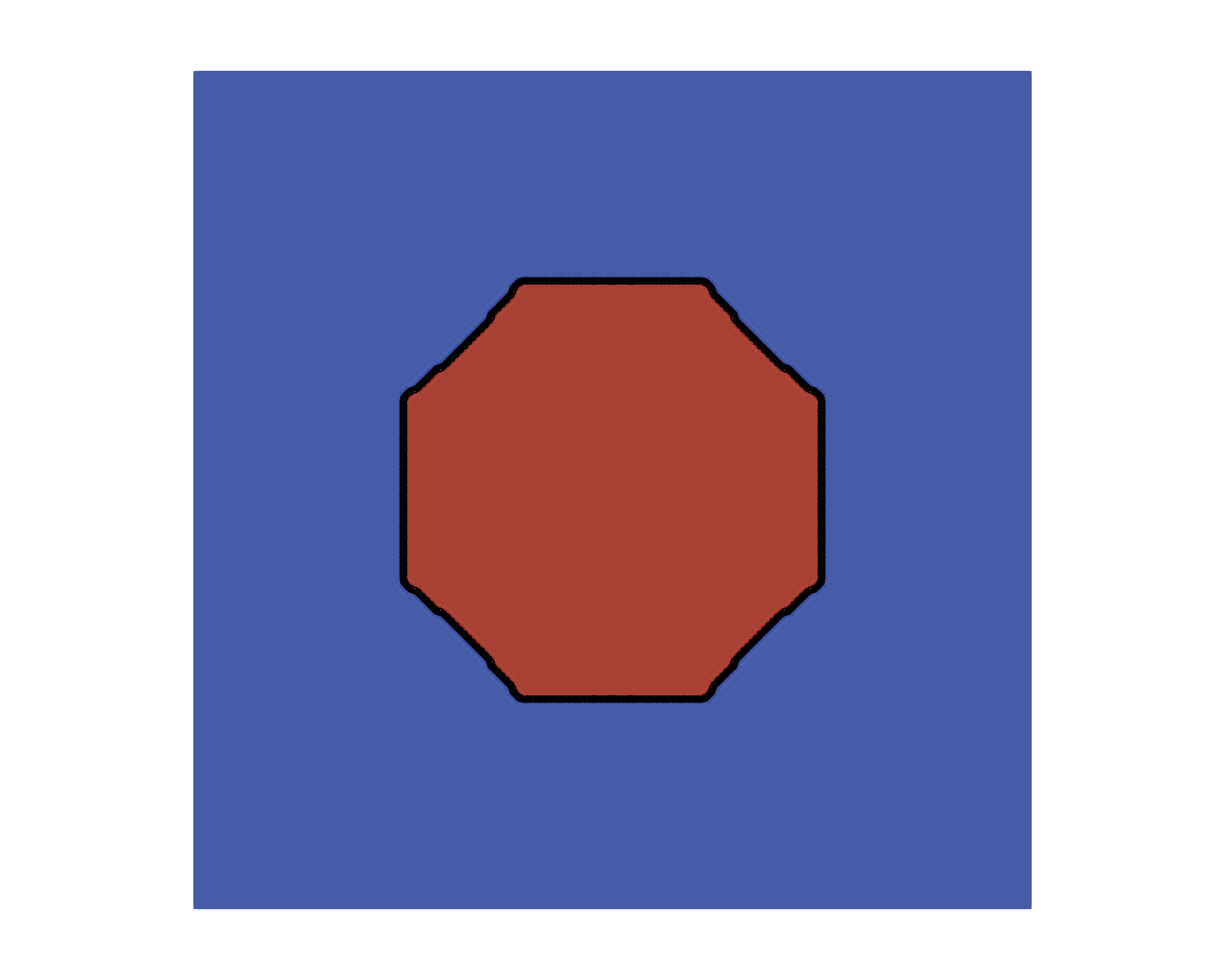}
    }
    \\
    \subfloat[{
        {
          $p=4$
        }  
    }]{
      \includegraphics[scale=\figscale,width=0.45\figwidth]{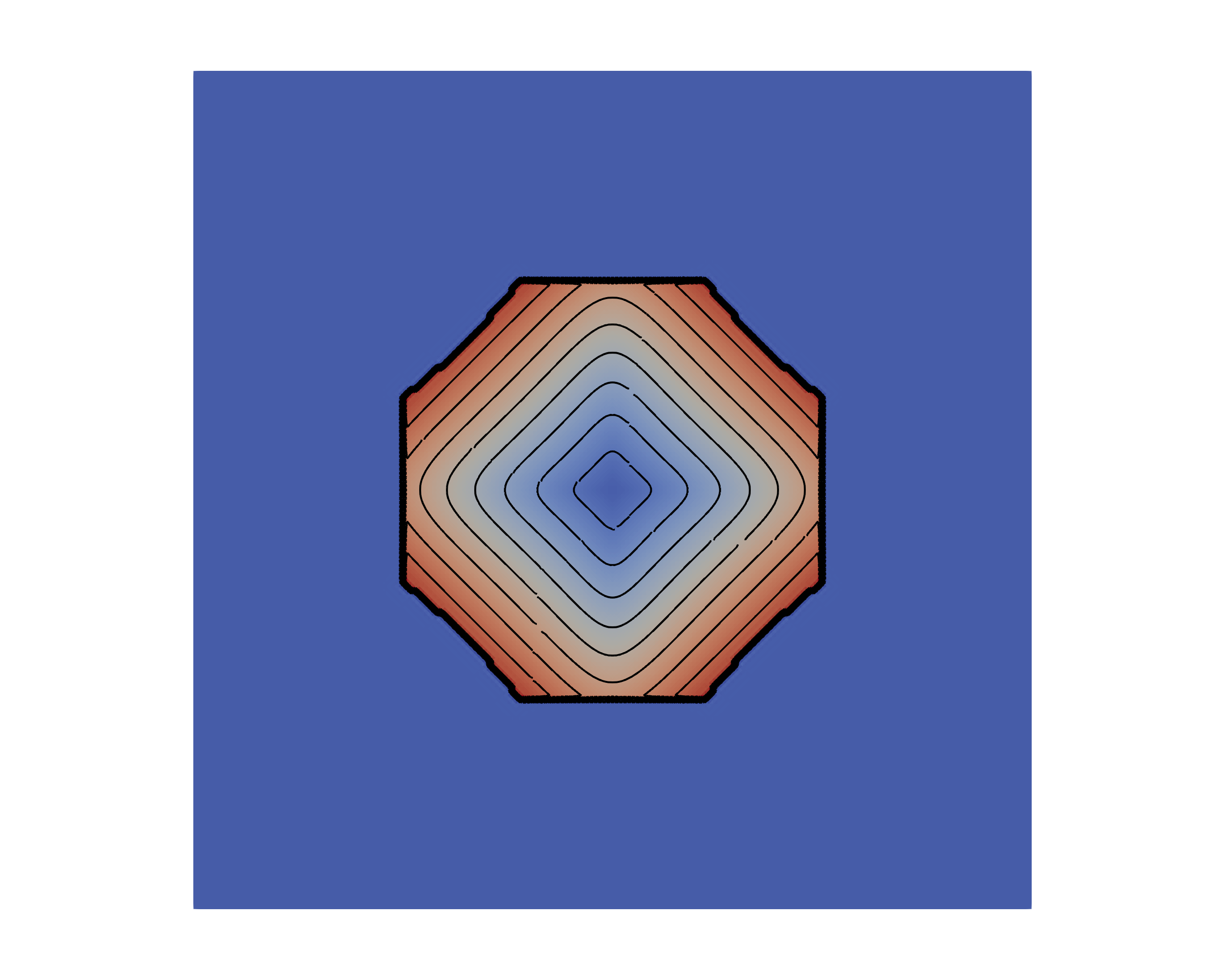}
    }
     \hfill
    \subfloat[{
        %%%%%%%%%%%%%%%%%%%%%%%%%%%%%%%%%%%%%%%%%%%%%%%%%%%%%%%%%%%%%%% 
        {
          $p=10$
        }
        %%%%%%%%%%%%%%%%%%%%%%%%%%%%%%%%%%%%%%%%%%%%%%%%%%%%%%%%%%%%%%%% 
    }]{
      \includegraphics[scale=\figscale,width=0.45\figwidth]{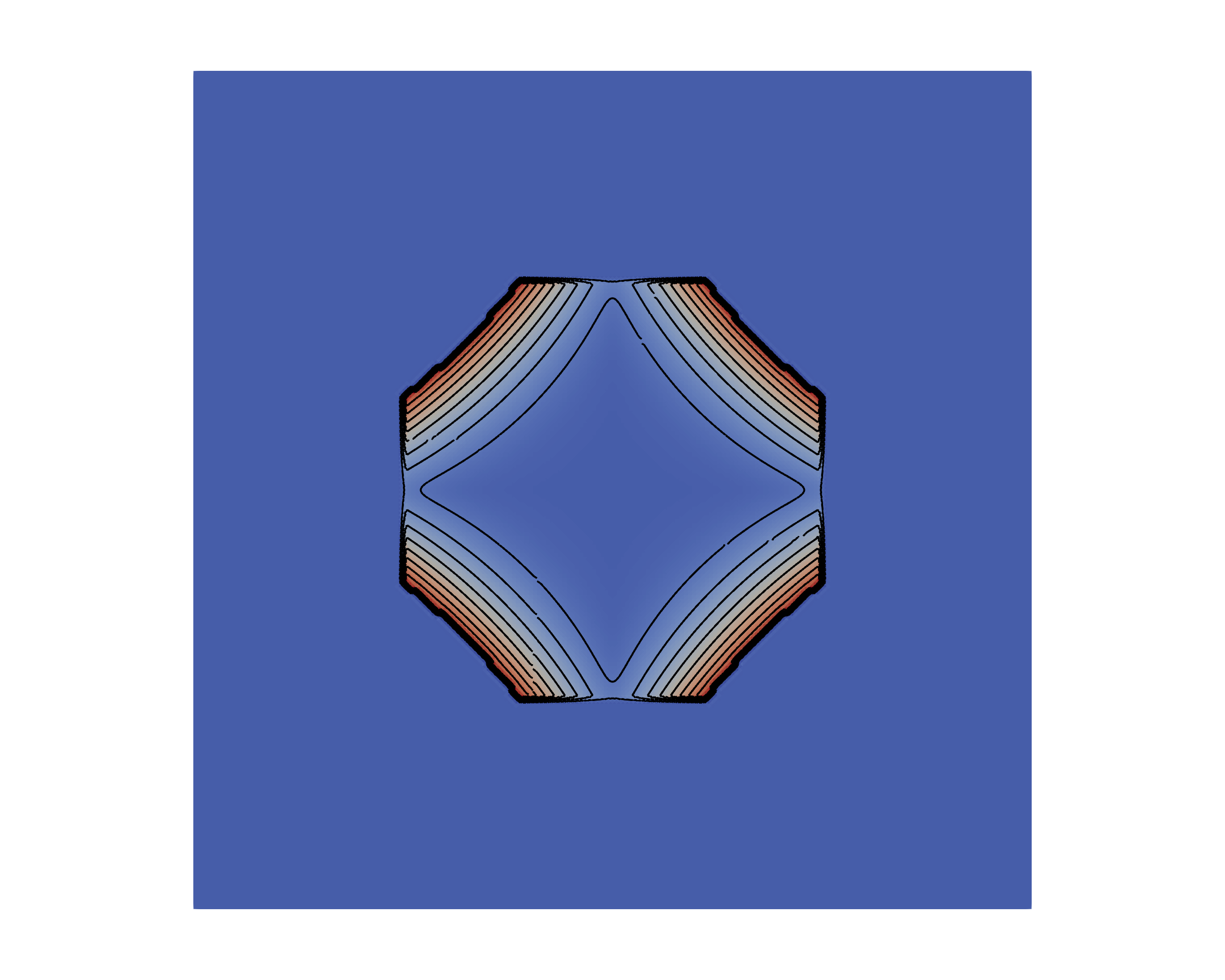}
        }
        \\
      \subfloat[{
      %%%%%%%%%%%%%%%%%%%%%%%%%%%%%%%%%%%%%%%%%%%%%%%%%%%%%%%%%%%%%%% 
        {
          $p=20$
        }
        %%%%%%%%%%%%%%%%%%%%%%%%%%%%%%%%%%%%%%%%%%%%%%%%%%%%%%%%%%%%%%%% 
    }]{
      \includegraphics[scale=\figscale,width=0.45\figwidth]{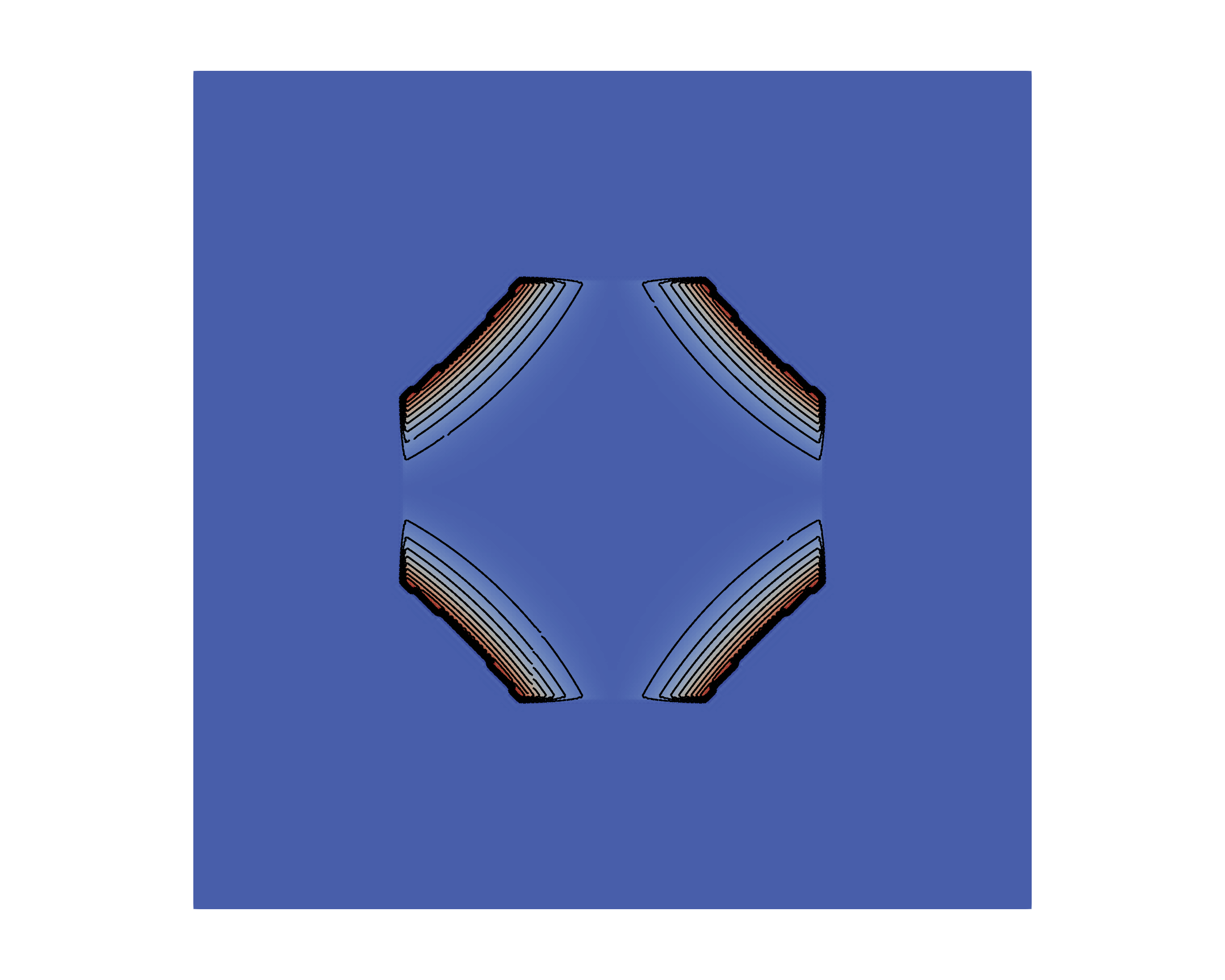}
    }
    %%%%%%%%%%%%%%%%%%%%%%%%%%%%%%%%%%%%%%%%%%%%%%%%%%%%%%%%%%%%%%%%%%%% 
     \hfill
     \subfloat[{
      %%%%%%%%%%%%%%%%%%%%%%%%%%%%%%%%%%%%%%%%%%%%%%%%%%%%%%%%%%%%%%% 
        {
          $p=100$
        }
        %%%%%%%%%%%%%%%%%%%%%%%%%%%%%%%%%%%%%%%%%%%%%%%%%%%%%%%%%%%%%%%% 
    }]{
      \includegraphics[scale=\figscale,width=0.45\figwidth]{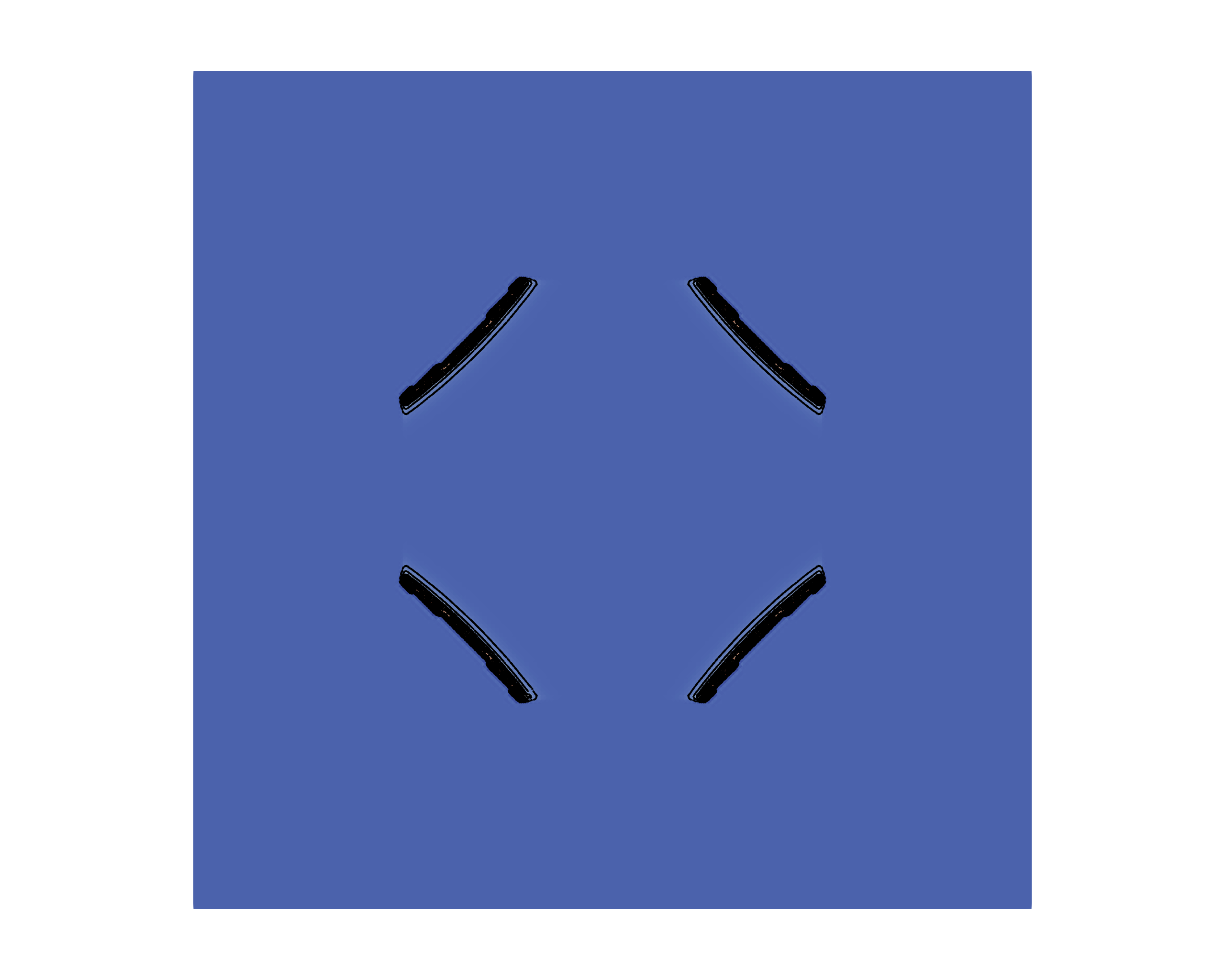}
    }
    %%%%%%%%%%%%%%%%%%%%%%%%%%%%%%%%%%%%%%%%%%%%%%%%%%%%%%%%%%%%%%%%%%%% 
  \end{center}
\end{figure}

Here we examine the scalar problem with prescribed boundary data given
by the celebrated Aronsson solution
\begin{equation}
  \label{eq:aronsson}
  u_\infty(x,y) = |x|^{4/3} - |y|^{4/3}.
\end{equation}
We study the measure $\sigma_p^\mO(\mathcal T)$ with various domains
$\mO$ and increasing values of $p$ to enable conjectures to be made as
to the structure of $\sigma_\infty^\mO$ in some test cases. In
particular, we examine the interplay between $\sigma_\infty^\mO$ and
the level sets of $|\D u_\infty|$. In Figures \ref{fig:t1}--\ref{fig:t4} we
plot the domain $\mO$ and $\sigma_p^\mO(\mathcal T)$ for
$p=2,4,10,20,100$. Specific conjectures regarding the behaviour of the measure $\si^\mO$ is provided at the captions of each figure, since the particular behaviour depends on the shape of the domain and its position in relation to the level sets of the gradient, but generally it is a sum of lower-dimensional concentration measures.

\ms

\subsection{Scalar case - Eikonal solution}

\begin{figure}[h!]
  \caption[] {\label{fig:scalareikonal} A test to characterise the
    $\sigma_\infty^\mO$ for boundary data given by the Eikonal
    function in (\ref{eq:cone}). We show the domain $\mO$ and the
    measure $\sigma_p^\mO(\mathcal T)$ for increasing values of
    $p$. In each of the subfigures B-F we plot 10 contours that are
    equally spaced between the minimum and maximum values of
    $\sigma_p^\mO$. Notice that as $p$ increases these contours remain
    equally spaced indicating mass is not pushed to the boundary.}
  \begin{center}
    \subfloat[{}]
             {
               \hspace{.2cm}
               \begin{tikzpicture}[scale=2]
                 \newcommand\Square[1]{+(-#1,-#1) rectangle +(#1,#1)}
                 \draw [very thin, gray] (-1,-1) grid (1,1);  
                 %\draw [red] (2,3) +(-2pt,-2pt) rectangle +(2pt,2pt) ;
                 \draw [fill=red] (.5,.5) \Square{.375};
                 \node at (0.5,0.5) {$\mO$};
                 \draw [fill=blue] (.875,.875) circle (.025);
                 %\node [above] at (.875,.875) {$\sigma_\infty^\mO = \mathcal L^2$};
                 \begin{scope}[thick,font=\scriptsize]
                   \draw [->] (-1.25,0) -- (1.25,0) node [above left]  {$x$};
                   \draw [->] (0,-1.25) -- (0,1.25) node [below right] {$y$};
                   \foreach \n in {-1,1}{%
                     \draw (\n,-3pt) -- (\n,3pt)   node [above] {$\n$};
                     \draw (-3pt,\n) -- (3pt,\n)   node [right] {$\n$};
                   }
                 \end{scope} 
                 %\path [draw=none,fill=gray,semitransparent] (+1,-1) circle (3);
               \end{tikzpicture}
             }
     \hfill
     \subfloat[{
        {
          $p=2$
        }
    }]{
      \includegraphics[scale=\figscale,width=0.45\figwidth]{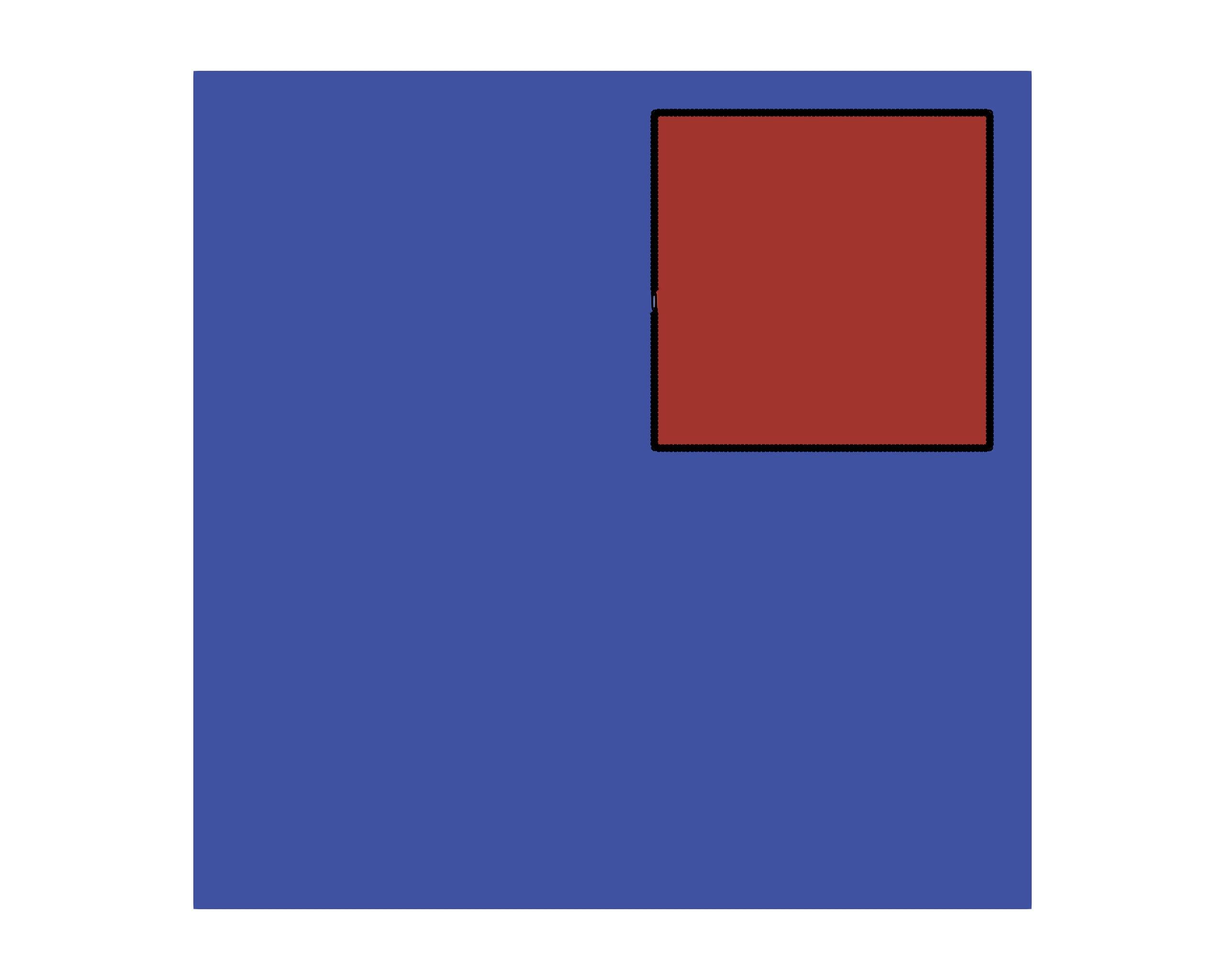}
    }
    \\
    \subfloat[{
        {
          $p=4$
        }  
    }]{
      \includegraphics[scale=\figscale,width=0.45\figwidth]{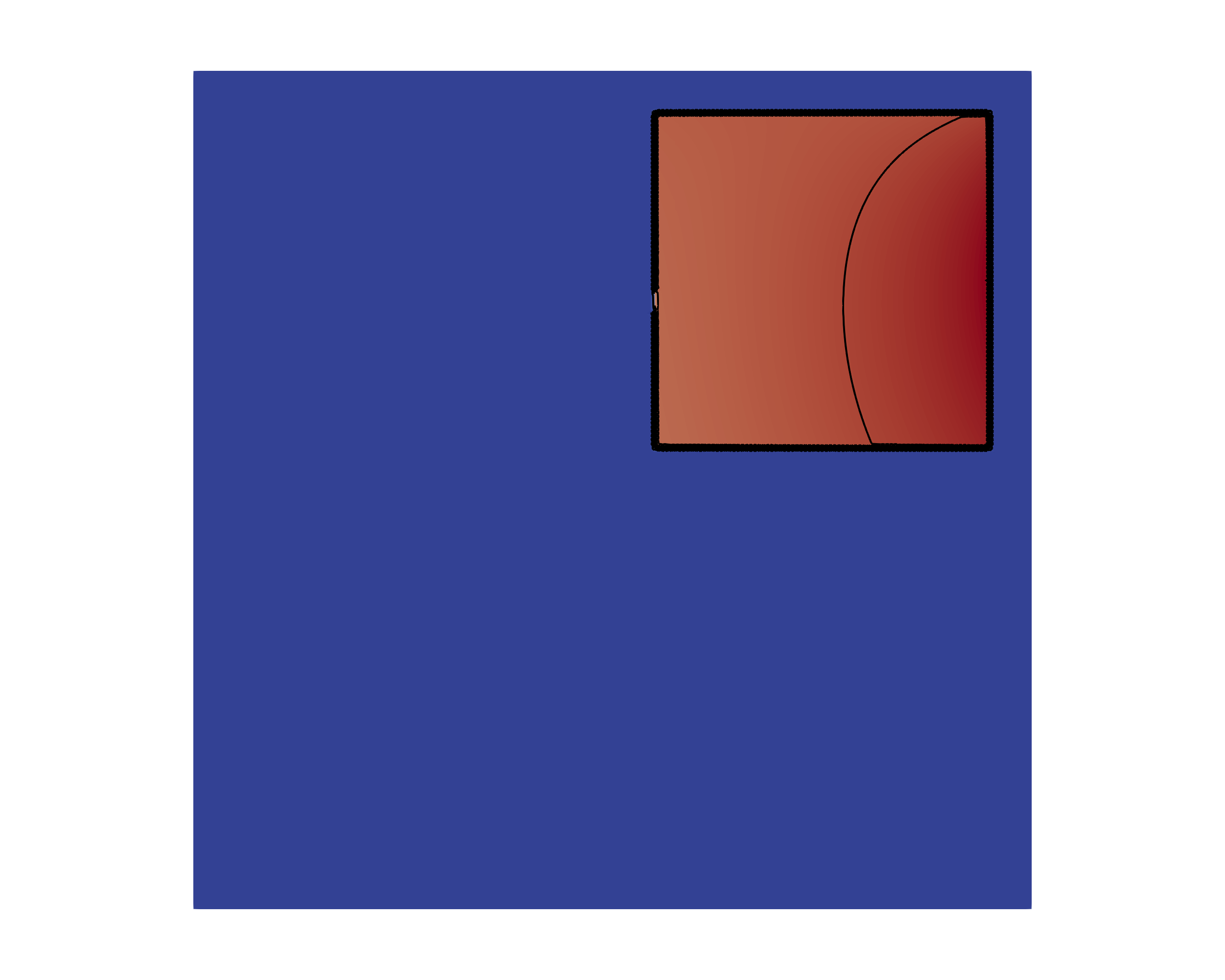}
    }
     \hfill
    \subfloat[{
        %%%%%%%%%%%%%%%%%%%%%%%%%%%%%%%%%%%%%%%%%%%%%%%%%%%%%%%%%%%%%%% 
        {
          $p=10$
        }
        %%%%%%%%%%%%%%%%%%%%%%%%%%%%%%%%%%%%%%%%%%%%%%%%%%%%%%%%%%%%%%%% 
    }]{
      \includegraphics[scale=\figscale,width=0.45\figwidth]{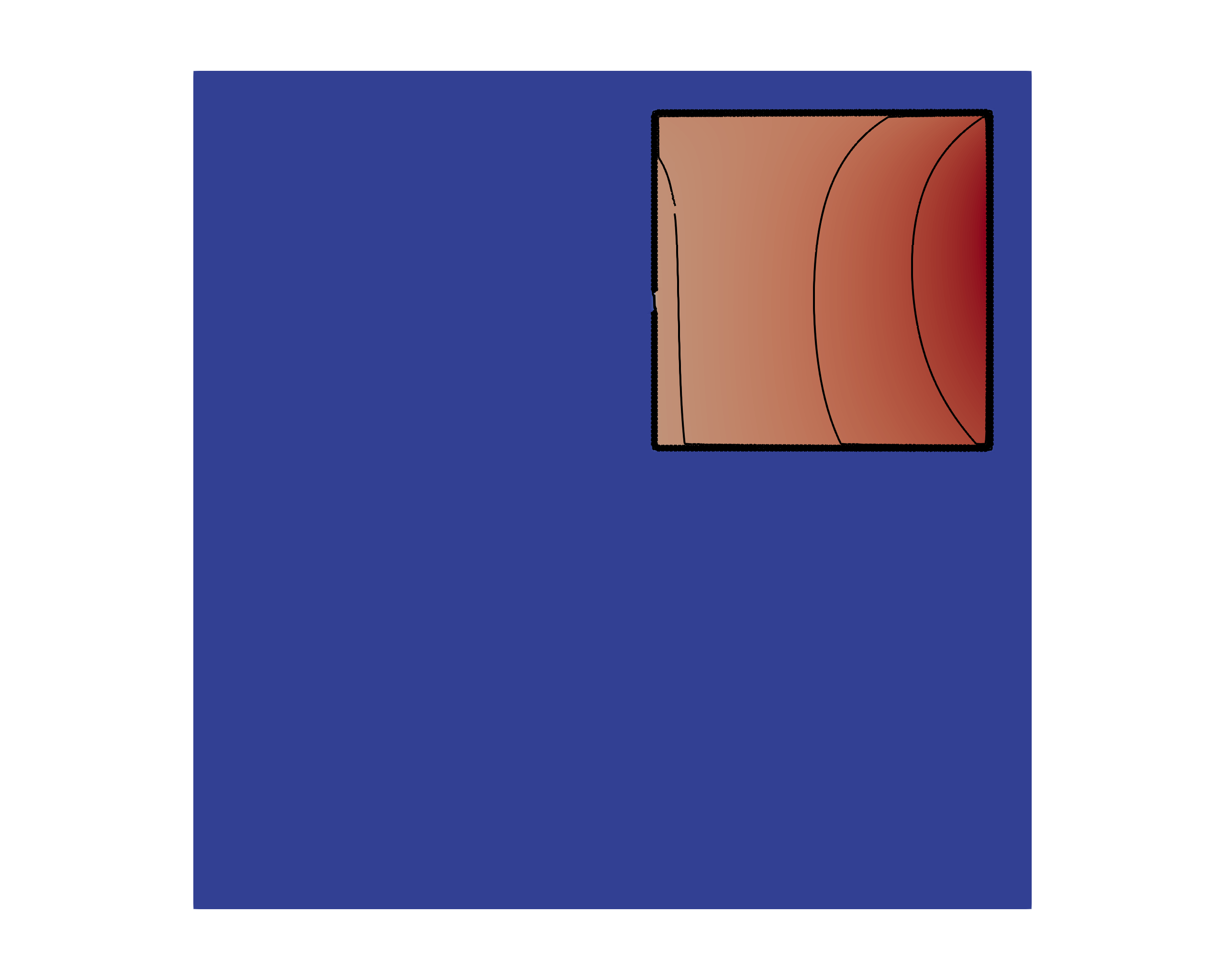}
        }
        \\
      \subfloat[{
      %%%%%%%%%%%%%%%%%%%%%%%%%%%%%%%%%%%%%%%%%%%%%%%%%%%%%%%%%%%%%%% 
        {
          $p=20$
        }
        %%%%%%%%%%%%%%%%%%%%%%%%%%%%%%%%%%%%%%%%%%%%%%%%%%%%%%%%%%%%%%%% 
    }]{
      \includegraphics[scale=\figscale,width=0.45\figwidth]{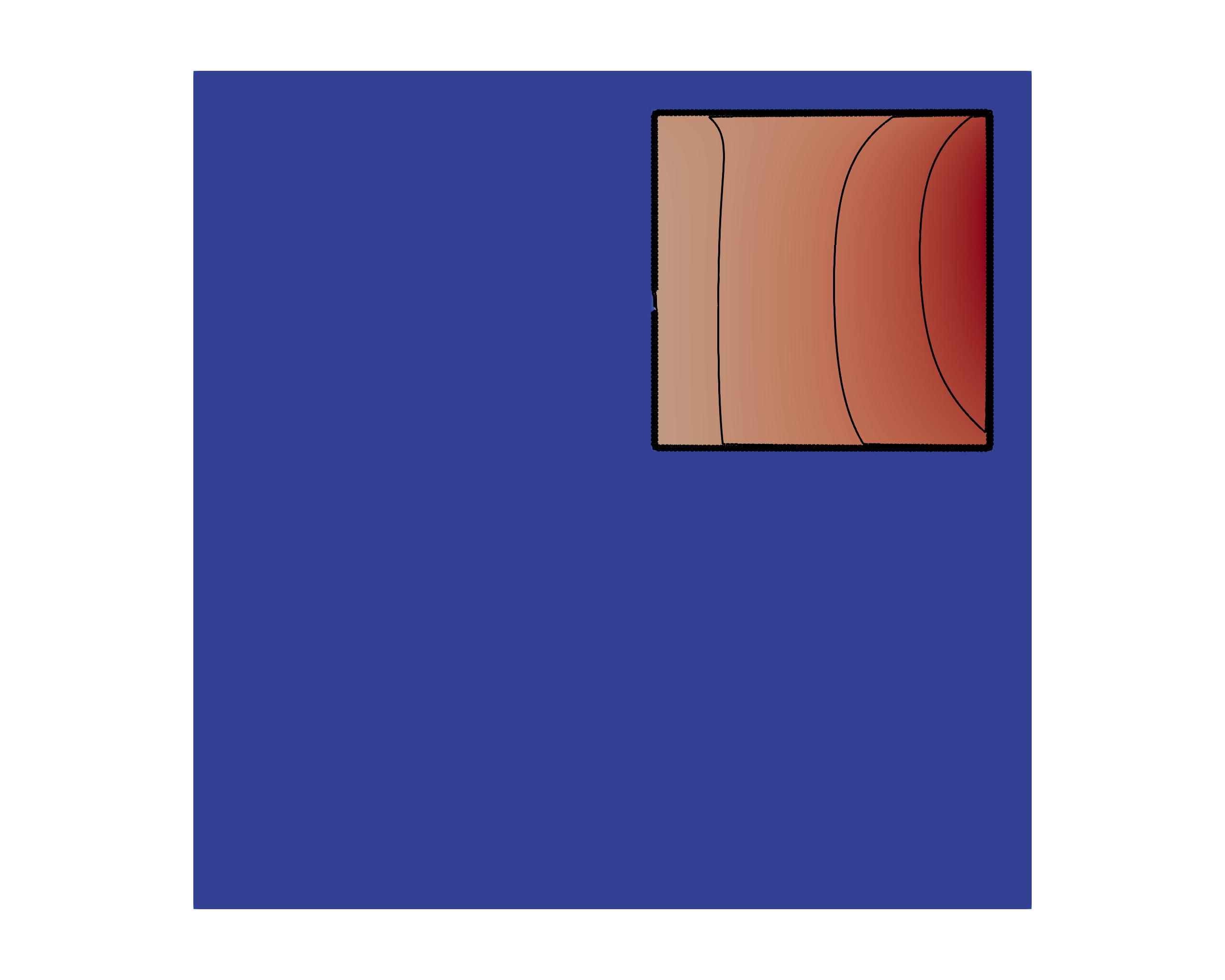}
    }
    %%%%%%%%%%%%%%%%%%%%%%%%%%%%%%%%%%%%%%%%%%%%%%%%%%%%%%%%%%%%%%%%%%%% 
     \hfill
     \subfloat[{
      %%%%%%%%%%%%%%%%%%%%%%%%%%%%%%%%%%%%%%%%%%%%%%%%%%%%%%%%%%%%%%% 
        {
          $p=100$
        }
        %%%%%%%%%%%%%%%%%%%%%%%%%%%%%%%%%%%%%%%%%%%%%%%%%%%%%%%%%%%%%%%% 
    }]{
       \includegraphics[scale=\figscale,width=0.45\figwidth]{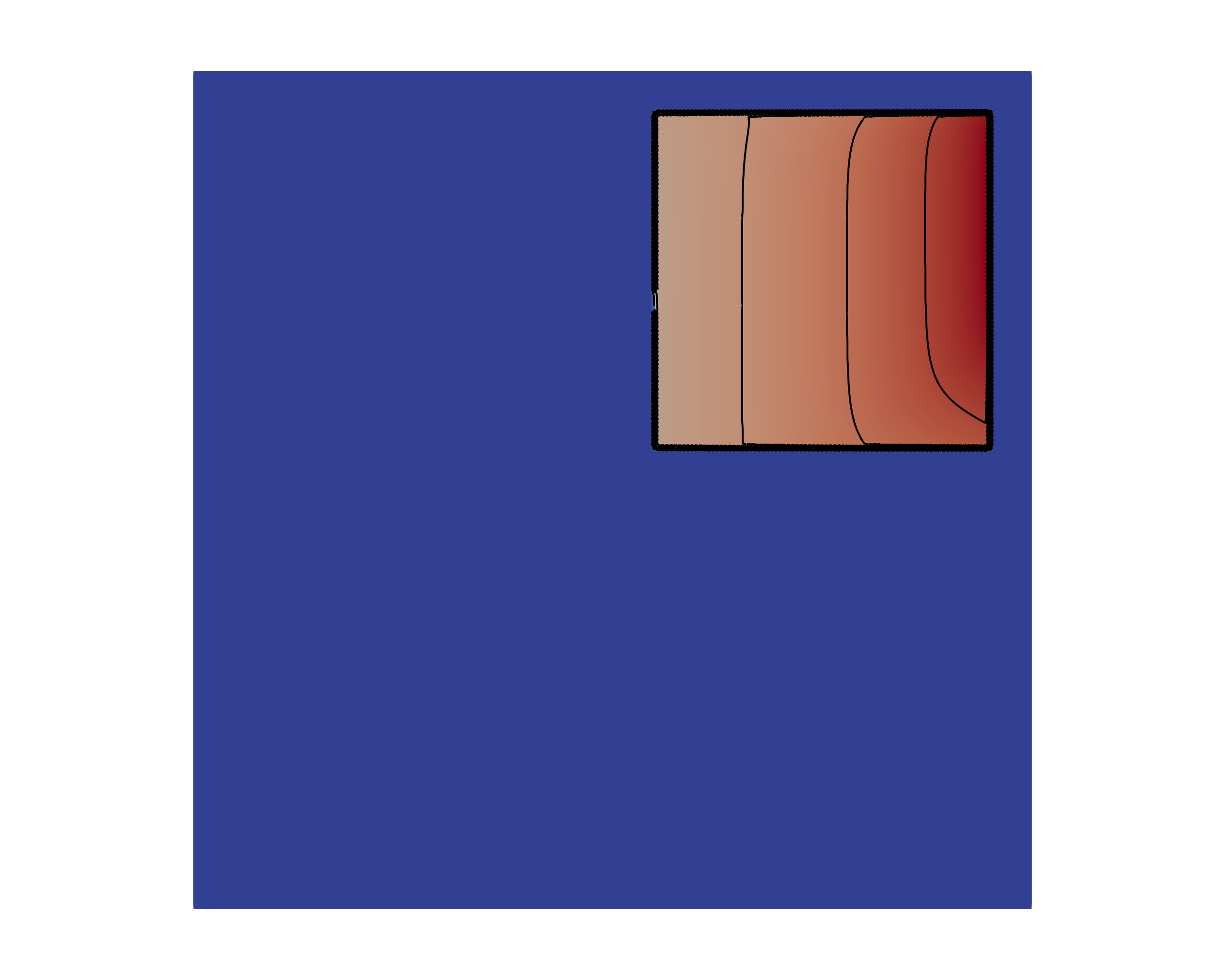}
    }
    %%%%%%%%%%%%%%%%%%%%%%%%%%%%%%%%%%%%%%%%%%%%%%%%%%%%%%%%%%%%%%%%%%%% 
  \end{center}
\end{figure}

In this test we examine the scalar problem with prescribed boundary
data given by the conic solution
\begin{equation}
  \label{eq:cone}
  u(x,y) = (x^2+y^2)^{1/2}.
\end{equation}
Notice this function is Eikonal. We study the measure
$\sigma_p^\mO(\mathcal T)$ with $\mO$ taken such that $0\notin
\mO$. Results are shown in Figure \ref{fig:scalareikonal} where we
plot the domain $\mO$ and $\sigma_p^\mO(\mathcal T)$ for
$p=2,4,10,20,100$. Notice that, in contrast to the Aronsson case, the
mass is not sent toward the boundary as $p\to\infty$. In this case we conjecture that $\si^\mO$ is always absolutely continuous with respect to the Lebesgue measure $\mL^n\LL_\mO$.

\ms

\subsection{Vectorial case - Eikonal solution}
\label{sec:vec}

\begin{figure}[h!]
  \caption[] {\label{fig:veceikonal} A test to characterise
    $\sigma^\mO_\infty$ for the boundary data given by
    (\ref{eq:cone}). We show the domain $\mO$ and the measure
    $\sigma_p^\mO(\mathcal T)$ for increasing values of $p$. In each
    of the subfigures B-F we plot 10 contours that are equally spaced
    between the minimum and maximum values of $\sigma_p^\mO$. Notice
    that as $p$ increases these contours remain equally spaced
    indicating mass is not pushed to the boundary.  }
  \begin{center}
    \subfloat[{}]
             {
               \hspace{.2cm}
               \begin{tikzpicture}[scale=2]
                 \newcommand\Square[1]{+(-#1,-#1) rectangle +(#1,#1)}
                 \draw [very thin, gray] (-1,-1) grid (1,1);  
                 %\draw [red] (2,3) +(-2pt,-2pt) rectangle +(2pt,2pt) ;
                 \draw [fill=red] (.0,.0) \Square{pi/4};
                 \node at (0.5,0.5) {$\mO$};
                 %\draw [fill=blue] (.875,.875) circle (.025);
                 %\node [above] at (.875,.875) {$\sigma_\infty^\mO = \mathcal L^2$};
                 \begin{scope}[thick,font=\scriptsize]
                   \draw [->] (-1.25,0) -- (1.25,0) node [above left]  {$x$};
                   \draw [->] (0,-1.25) -- (0,1.25) node [below right] {$y$};
                   \foreach \n in {-1,1}{%
                     \draw (\n,-3pt) -- (\n,3pt)   node [above] {$\n$};
                     \draw (-3pt,\n) -- (3pt,\n)   node [right] {$\n$};
                   }
                 \end{scope} 
                 %\path [draw=none,fill=gray,semitransparent] (+1,-1) circle (3);
               \end{tikzpicture}
             }
     \hfill
     \subfloat[{
        {
          $p=2$
        }
    }]{
      \includegraphics[scale=\figscale,width=0.45\figwidth]{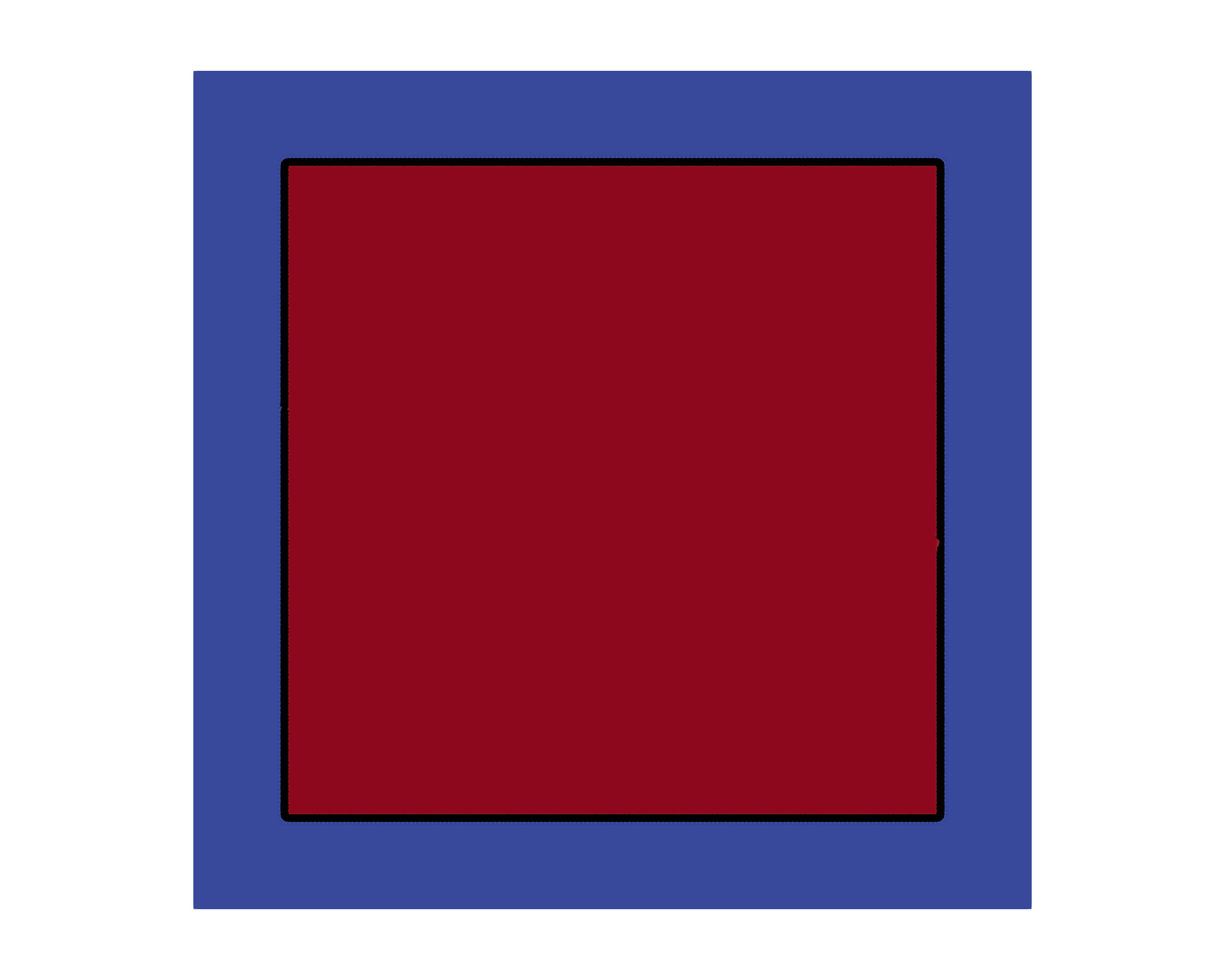}
    }
    \\
    \subfloat[{
        {
          $p=4$
        }  
    }]{
      \includegraphics[scale=\figscale,width=0.45\figwidth]{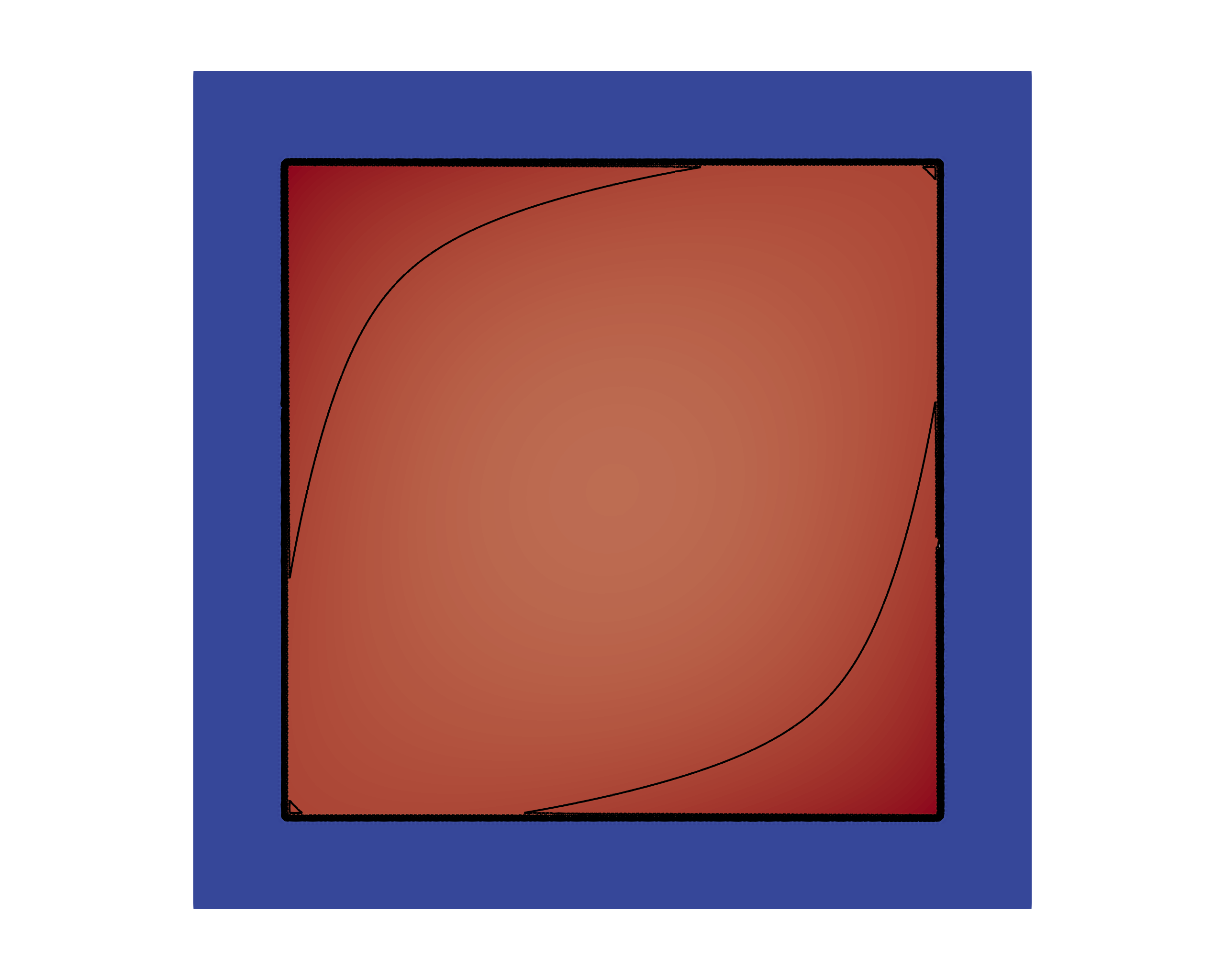}
    }
     \hfill
    \subfloat[{
        %%%%%%%%%%%%%%%%%%%%%%%%%%%%%%%%%%%%%%%%%%%%%%%%%%%%%%%%%%%%%%% 
        {
          $p=10$
        }
        %%%%%%%%%%%%%%%%%%%%%%%%%%%%%%%%%%%%%%%%%%%%%%%%%%%%%%%%%%%%%%%% 
    }]{
      \includegraphics[scale=\figscale,width=0.45\figwidth]{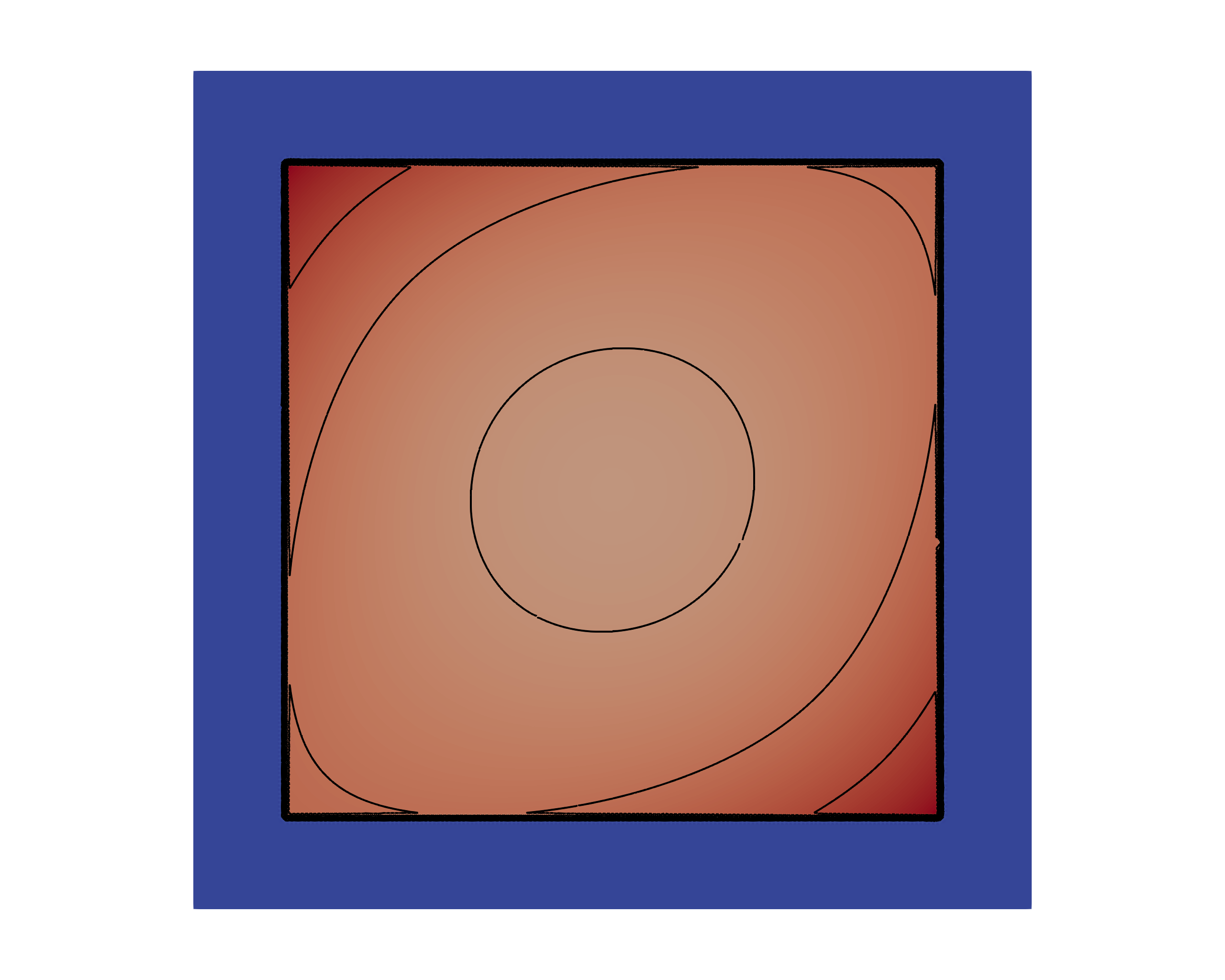}
        }
        \\
      \subfloat[{
      %%%%%%%%%%%%%%%%%%%%%%%%%%%%%%%%%%%%%%%%%%%%%%%%%%%%%%%%%%%%%%% 
        {
          $p=20$
        }
        %%%%%%%%%%%%%%%%%%%%%%%%%%%%%%%%%%%%%%%%%%%%%%%%%%%%%%%%%%%%%%%% 
    }]{
      \includegraphics[scale=\figscale,width=0.45\figwidth]{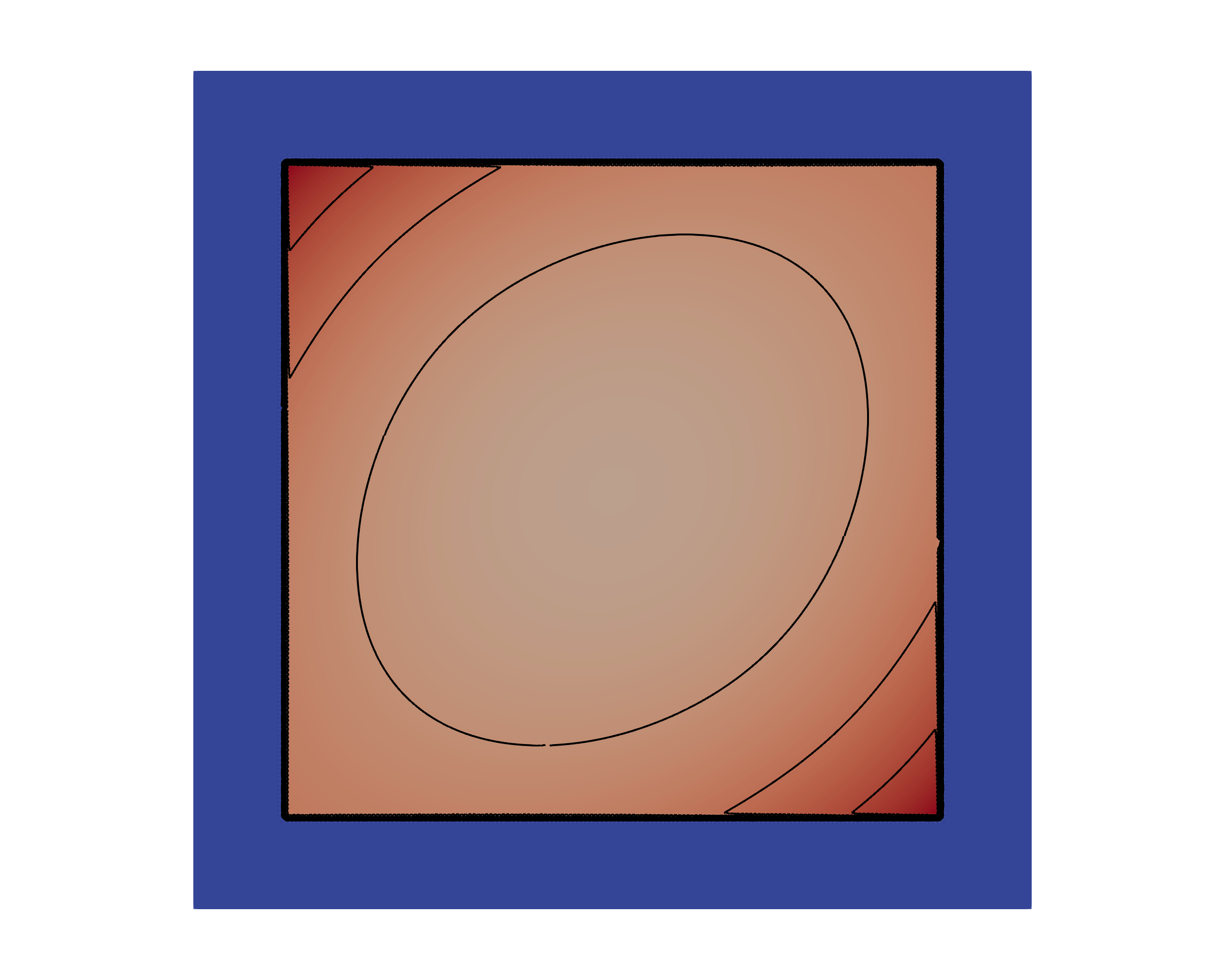}
    }
    %%%%%%%%%%%%%%%%%%%%%%%%%%%%%%%%%%%%%%%%%%%%%%%%%%%%%%%%%%%%%%%%%%%% 
     \hfill
     \subfloat[{
      %%%%%%%%%%%%%%%%%%%%%%%%%%%%%%%%%%%%%%%%%%%%%%%%%%%%%%%%%%%%%%% 
        {
          $p=100$
        }
        %%%%%%%%%%%%%%%%%%%%%%%%%%%%%%%%%%%%%%%%%%%%%%%%%%%%%%%%%%%%%%%% 
    }]{
       \includegraphics[scale=\figscale,width=0.45\figwidth]{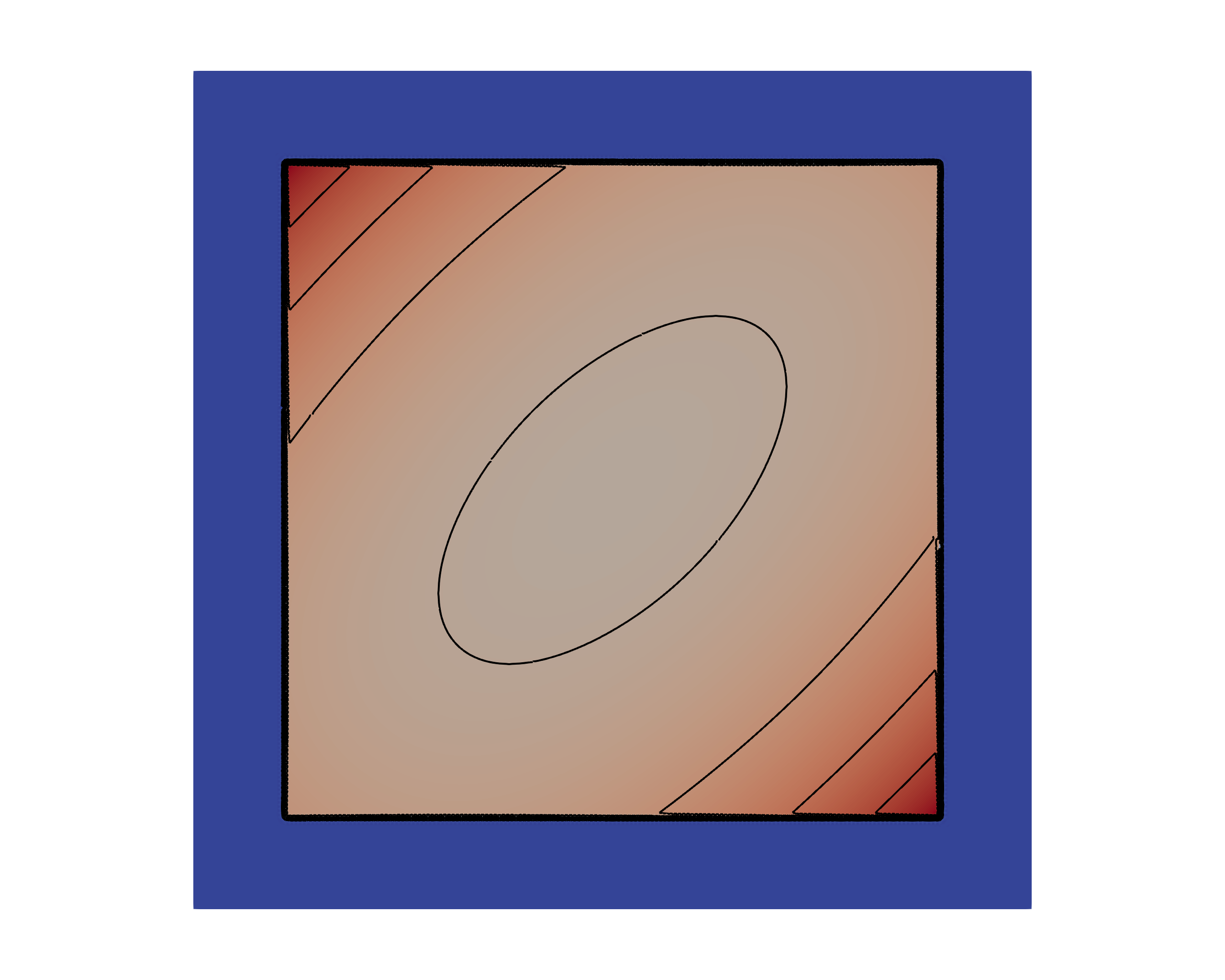}
    }
    %%%%%%%%%%%%%%%%%%%%%%%%%%%%%%%%%%%%%%%%%%%%%%%%%%%%%%%%%%%%%%%%%%%% 
  \end{center}
\end{figure}
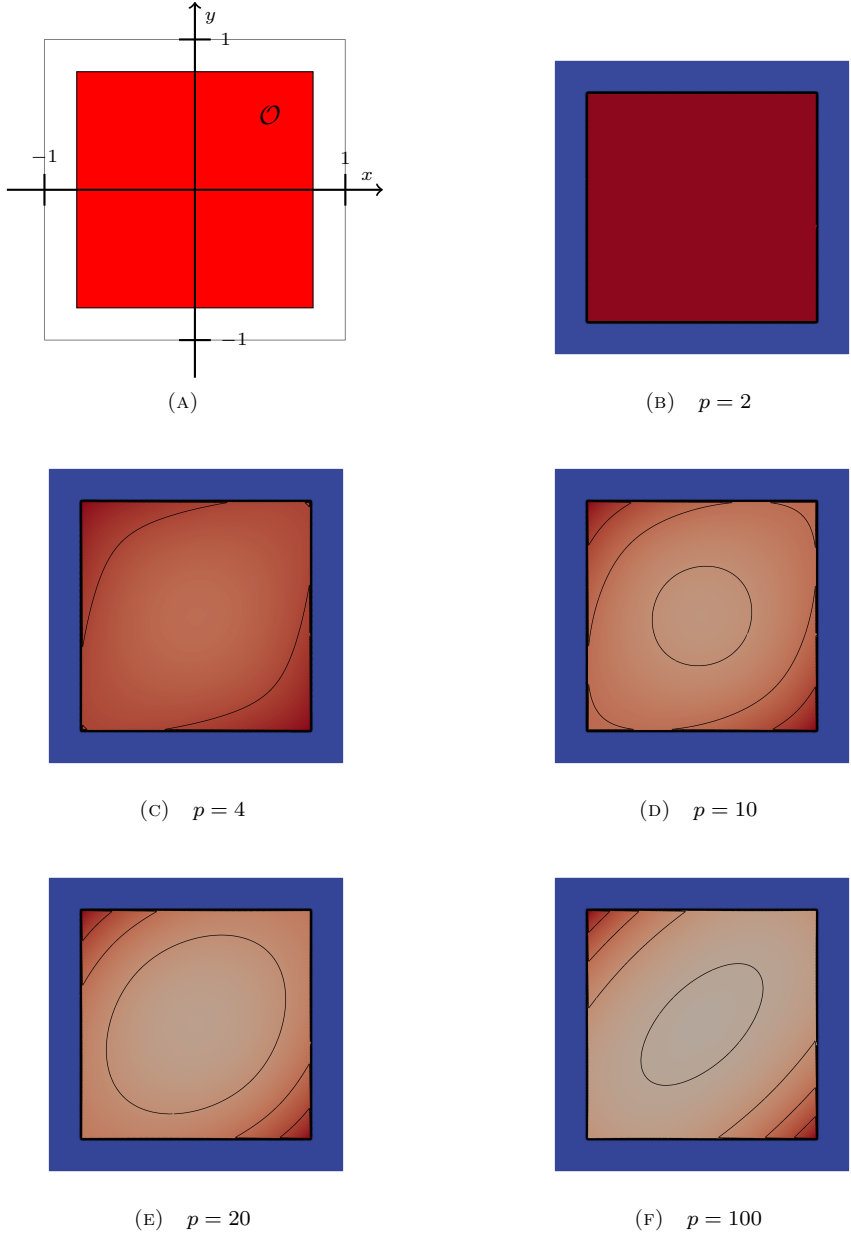

In this test we examine the vectorial problem with prescribed boundary
data given by
\begin{equation}
  \label{eq:cone}
  u(x,y) = e^{i x} - e^{i y},
\end{equation}
where we use the notation $e^{i t} = (\cos t, \sin t)$. Notice this
function is Eikonal. We study the measure $\sigma_p^\mO(\mathcal T)$
with $\mO$ taken such that $0\notin \mO$. Results are shown in Figure
\ref{fig:veceikonal} where we plot the domain $\mO$ and
$\sigma_p^\mO(\mathcal T)$ for $p=2,4,10,20,100$. Notice that mass is
not sent towards the boundary as $p\to\infty$.  We conjecture that again $\si^\mO$ is absolutely continuous with respect to the Lebesgue measure $\mL^n\LL_\mO$.

\ms

\subsection{Vectorial case - Mixed boundary conditions}
 
\begin{figure}[h!]
  \caption[] {\label{fig:mixed1} A test to characterise
    $\sigma^\mO_\infty$ for the boundary data given by (\ref{eq:lamb})
    with $\lambda = \frac 1 2$. We show the domain $\mO$ and the
    measure $\sigma_p^\mO(\mathcal T)$ for increasing values of
    $p$. In each of the subfigures B-F we plot 10 contours that are
    equally spaced between the minimum and maximum values of
    $\sigma_p^\mO$. Notice that as $p$ increases these contours become
    concentrated at the four points $\{x_i\}_{i=0}^3$.  In this case we conjecture that $\si^\mO$ is given by the sum of four Dirac masses at the points $\{x_i\}_{i=0}^3$ and otherwise it is absolutely continuous with respect to $\mL^n\LL_\mO$, rescaled so that it is a probability on $\overline \mO$.}
  \begin{center}
    \subfloat[{}]
             {
               \hspace{.2cm}
               \begin{tikzpicture}[scale=2]
                 \newcommand\Square[1]{+(-#1,-#1) rectangle +(#1,#1)}
                 \draw [very thin, gray] (-1,-1) grid (1,1);  
                 %\draw [red] (2,3) +(-2pt,-2pt) rectangle +(2pt,2pt) ;
                 \draw [fill=red] (.0,.0) \Square{1};
                 \node at (0.5,0.5) {$\mO$};
                 \draw [fill=blue] (-1,1) circle (.025);
                 \node [left] at (-1,1) {{$x_2$}};
                 \draw [fill=blue] (0,-1) circle (.025);
                 \node [below right] at (0,-1) {{$x_1$}};
                 \draw [fill=blue] (-1,-1) circle (.025);
                 \node [left] at (-1,-1) {{$x_0$}};
                 \draw [fill=blue] (0,1) circle (.025);
                 \node [above left] at (0,1) {{$x_3$}};
                 \begin{scope}[thick,font=\scriptsize]
                   \draw [->] (-1.25,0) -- (1.25,0) node [above left]  {$x$};
                   \draw [->] (0,-1.25) -- (0,1.25) node [below right] {$y$};
                   \foreach \n in {-1,1}{%
                     \draw (\n,-3pt) -- (\n,3pt)   node [above] {$\n$};
                     \draw (-3pt,\n) -- (3pt,\n)   node [right] {$\n$};
                   }
                 \end{scope} 
                 %\path [draw=none,fill=gray,semitransparent] (+1,-1) circle (3);
               \end{tikzpicture}
             }
     \hfill
     \subfloat[{
        {
          $p=2$
        }
    }]{
      \includegraphics[scale=\figscale,width=0.45\figwidth]{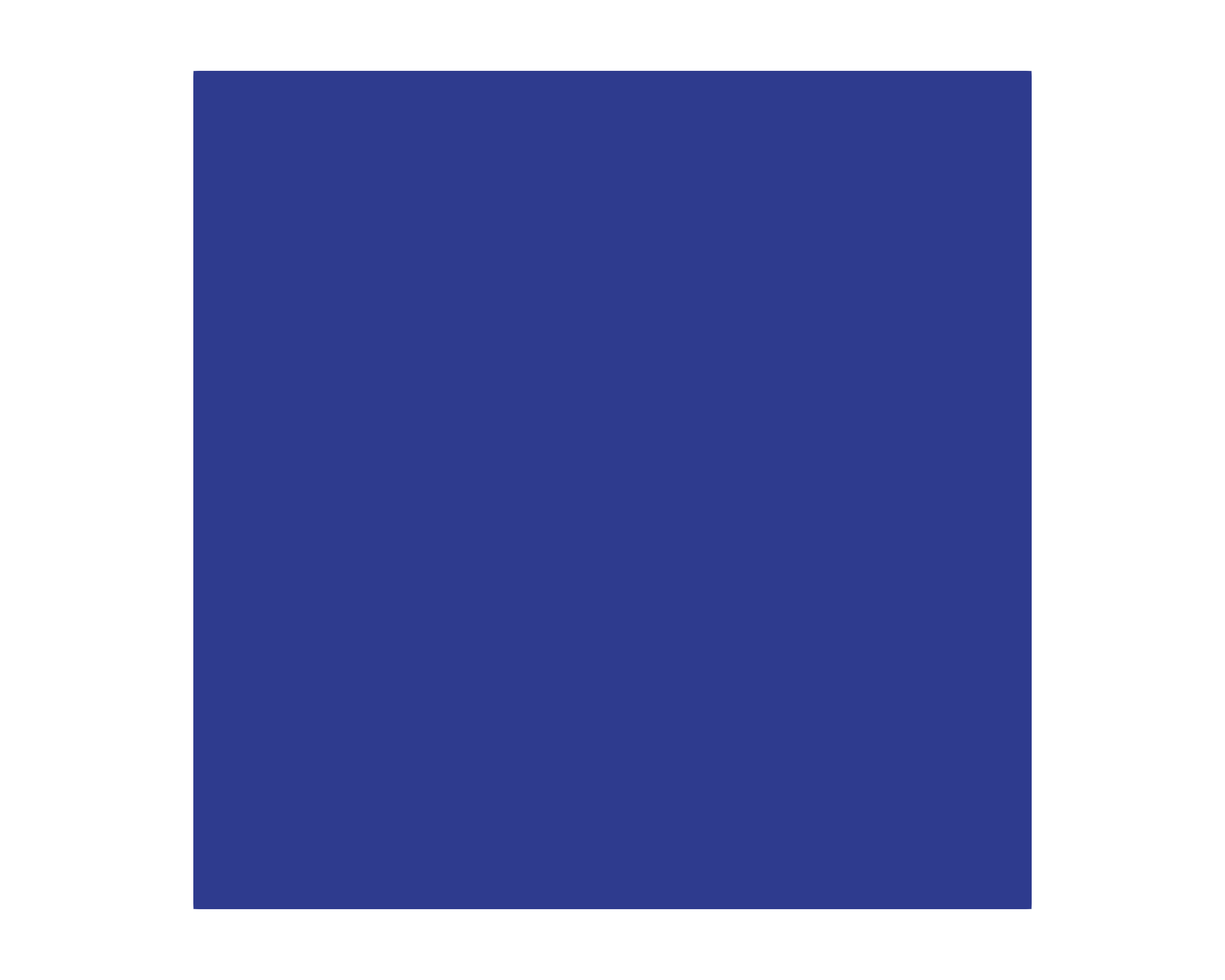}
    }
    \\
    \subfloat[{
        {
          $p=4$
        }  
    }]{
      \includegraphics[scale=\figscale,width=0.45\figwidth]{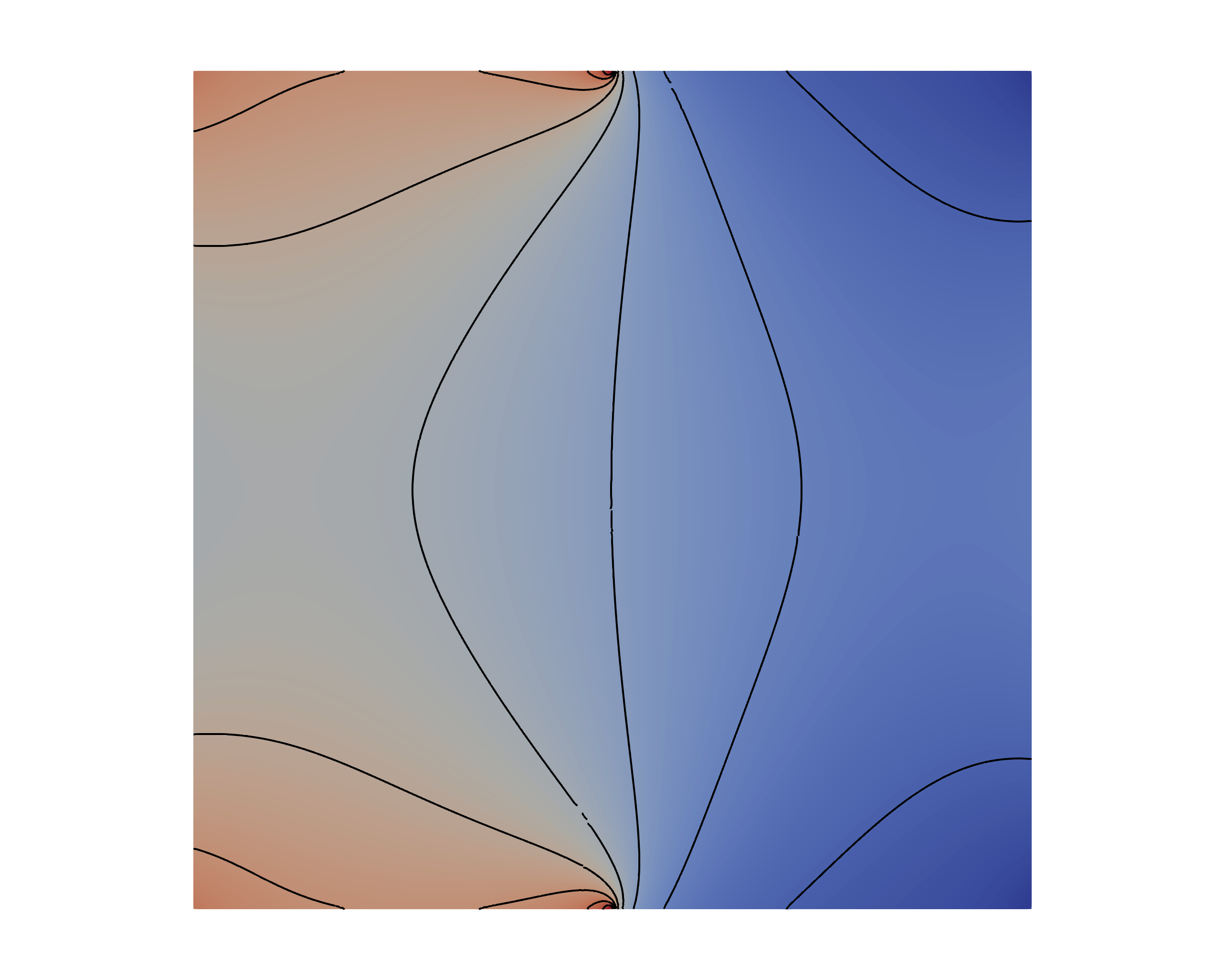}
    }
     \hfill
    \subfloat[{
        %%%%%%%%%%%%%%%%%%%%%%%%%%%%%%%%%%%%%%%%%%%%%%%%%%%%%%%%%%%%%%% 
        {
          $p=10$
        }
        %%%%%%%%%%%%%%%%%%%%%%%%%%%%%%%%%%%%%%%%%%%%%%%%%%%%%%%%%%%%%%%% 
    }]{
      \includegraphics[scale=\figscale,width=0.45\figwidth]{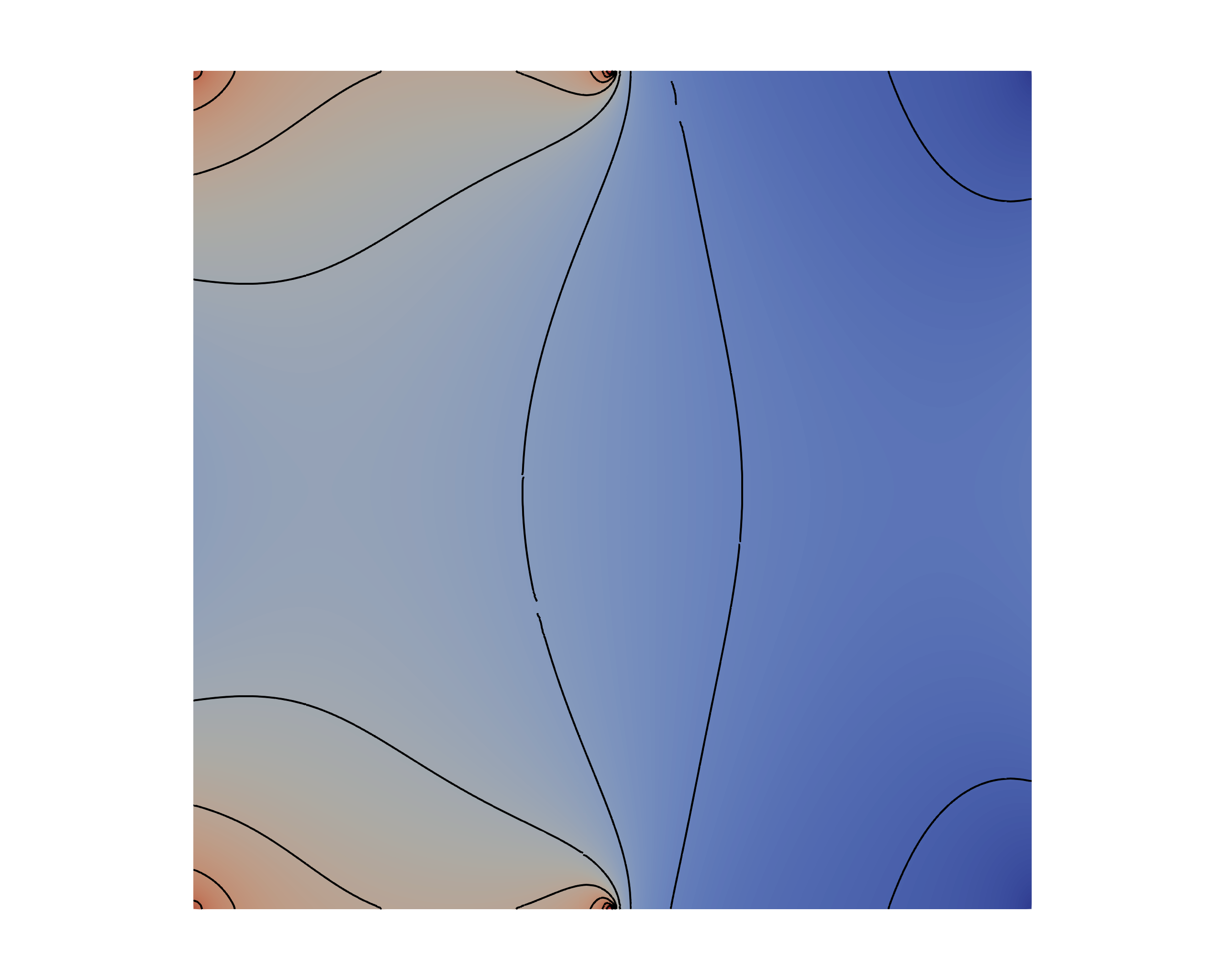}
        }
        \\
      \subfloat[{
      %%%%%%%%%%%%%%%%%%%%%%%%%%%%%%%%%%%%%%%%%%%%%%%%%%%%%%%%%%%%%%% 
        {
          $p=20$
        }
        %%%%%%%%%%%%%%%%%%%%%%%%%%%%%%%%%%%%%%%%%%%%%%%%%%%%%%%%%%%%%%%% 
    }]{
      \includegraphics[scale=\figscale,width=0.45\figwidth]{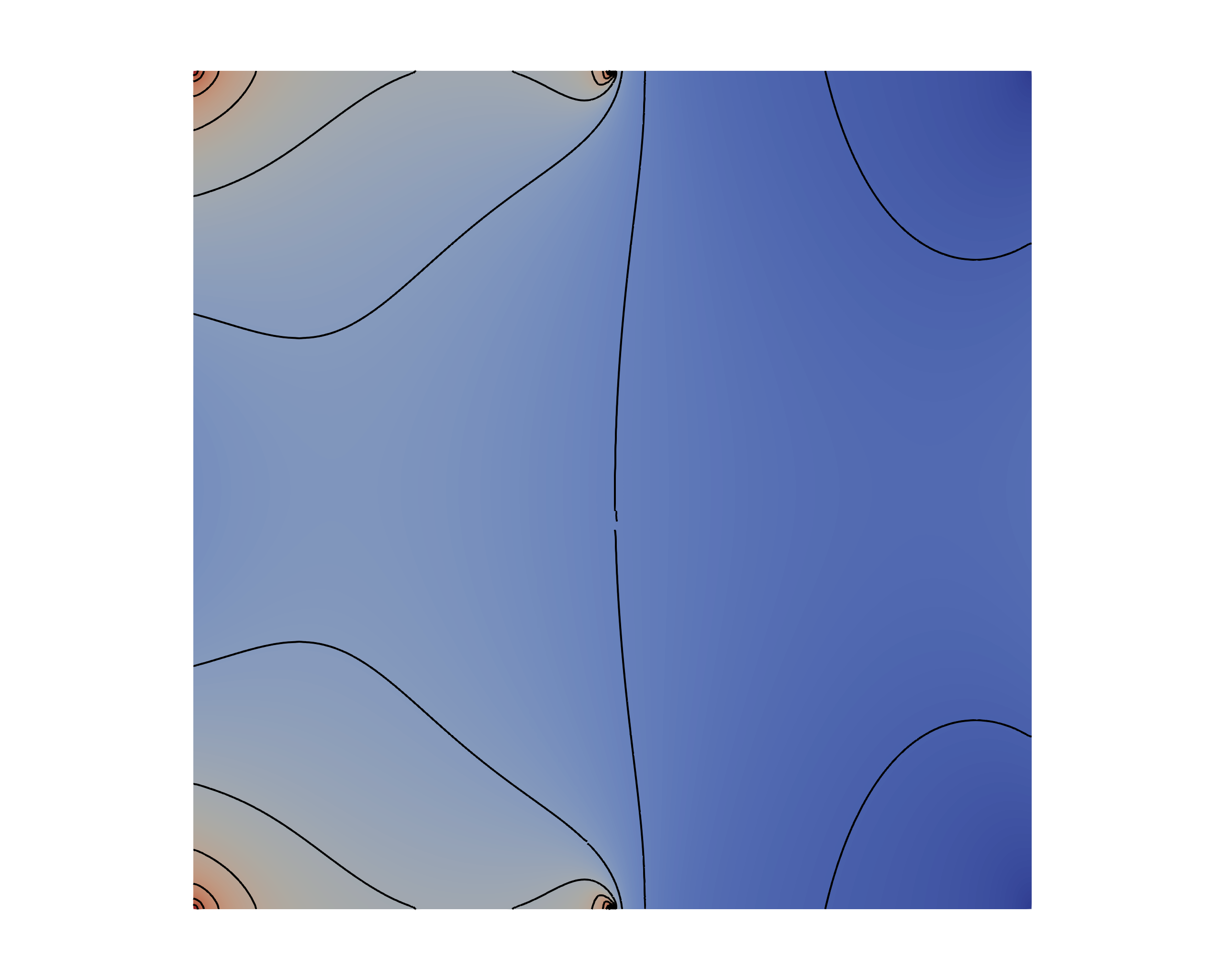}
    }
    %%%%%%%%%%%%%%%%%%%%%%%%%%%%%%%%%%%%%%%%%%%%%%%%%%%%%%%%%%%%%%%%%%%% 
     \hfill
     \subfloat[{
      %%%%%%%%%%%%%%%%%%%%%%%%%%%%%%%%%%%%%%%%%%%%%%%%%%%%%%%%%%%%%%% 
        {
          $p=100$
        }
        %%%%%%%%%%%%%%%%%%%%%%%%%%%%%%%%%%%%%%%%%%%%%%%%%%%%%%%%%%%%%%%% 
    }]{
       \includegraphics[scale=\figscale,width=0.45\figwidth]{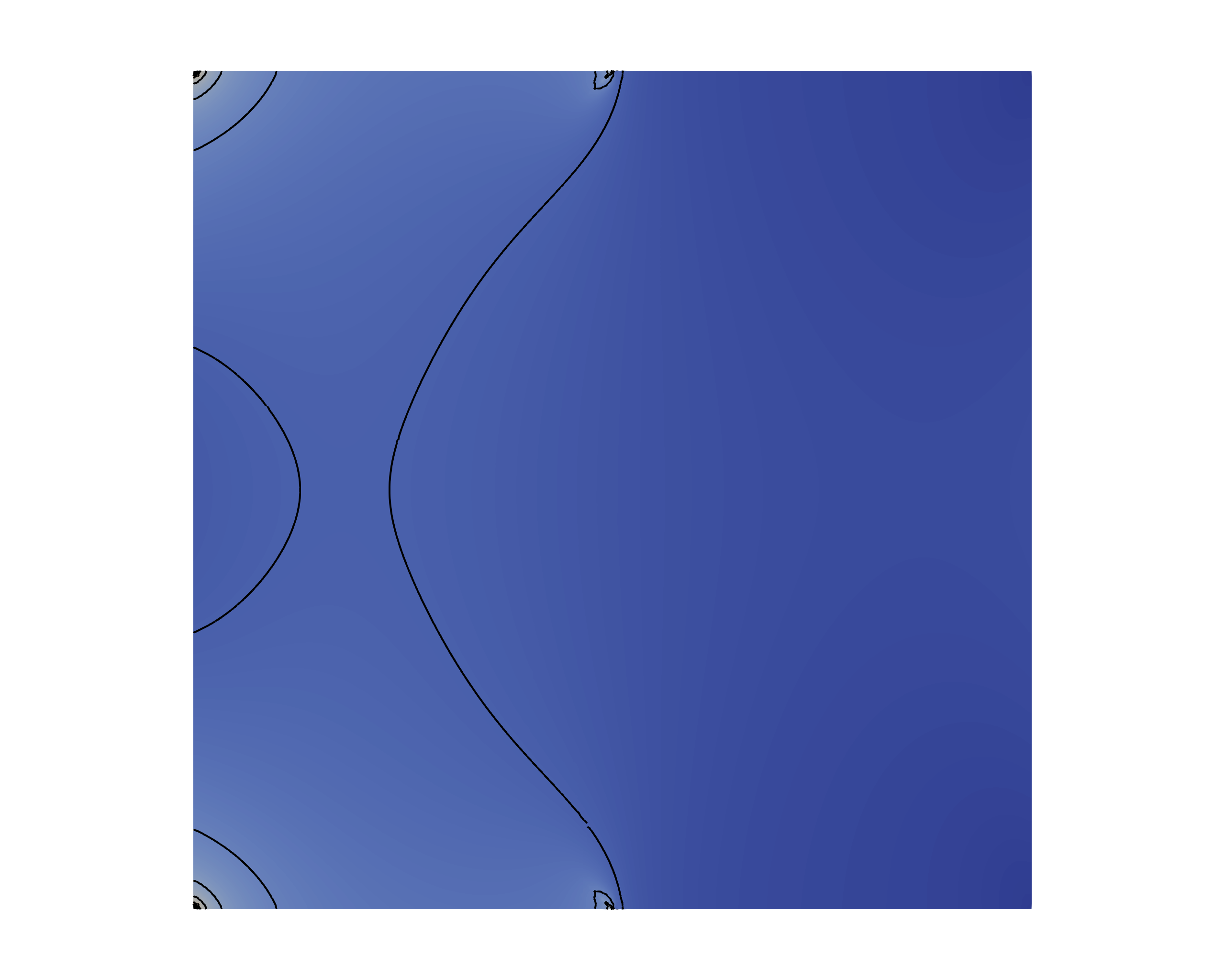}
    }
    %%%%%%%%%%%%%%%%%%%%%%%%%%%%%%%%%%%%%%%%%%%%%%%%%%%%%%%%%%%%%%%%%%%% 
  \end{center}
\end{figure}

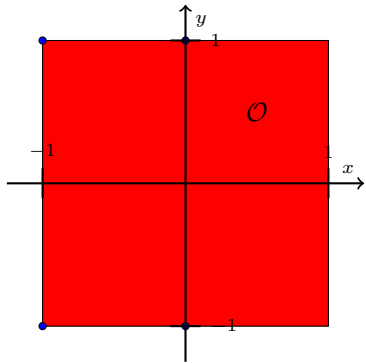
\begin{figure}[h!]
  \caption[]
  {\label{fig:mixed2}
    As Figure \ref{fig:mixed1} with $\lambda = -\frac 1 2$. Qualitatively the same behaviour is observed here as well.
  }
  \begin{center}
    \subfloat[{}]
             {
               \hspace{.2cm}
               \begin{tikzpicture}[scale=1.9]
                 \newcommand\Square[1]{+(-#1,-#1) rectangle +(#1,#1)}
                 \draw [very thin, gray] (-1,-1) grid (1,1);  
                 %\draw [red] (2,3) +(-2pt,-2pt) rectangle +(2pt,2pt) ;
                 \draw [fill=red] (.0,.0) \Square{1};
                 \node at (0.5,0.5) {$\mO$};
                 \draw [fill=blue] (0.,1.) circle (.025);
                 \draw [fill=blue] (0.,-1.) circle (.025);
                 \draw [fill=blue] (-1.,1.) circle (.025);
                 \draw [fill=blue] (-1.,-1.) circle (.025);
                 %\node [above] at (.875,.875) {$\sigma_\infty^\mO = \mathcal L^2$};
                 \begin{scope}[thick,font=\scriptsize]
                   \draw [->] (-1.25,0) -- (1.25,0) node [above left]  {$x$};
                   \draw [->] (0,-1.25) -- (0,1.25) node [below right] {$y$};
                   \foreach \n in {-1,1}{%
                     \draw (\n,-3pt) -- (\n,3pt)   node [above] {$\n$};
                     \draw (-3pt,\n) -- (3pt,\n)   node [right] {$\n$};
                   }
                 \end{scope} 
                 %\path [draw=none,fill=gray,semitransparent] (+1,-1) circle (3);
               \end{tikzpicture}
             }
     \hfill
     \subfloat[{
        {
          $p=2$
        }
    }]{
      \includegraphics[scale=\figscale,width=0.45\figwidth]{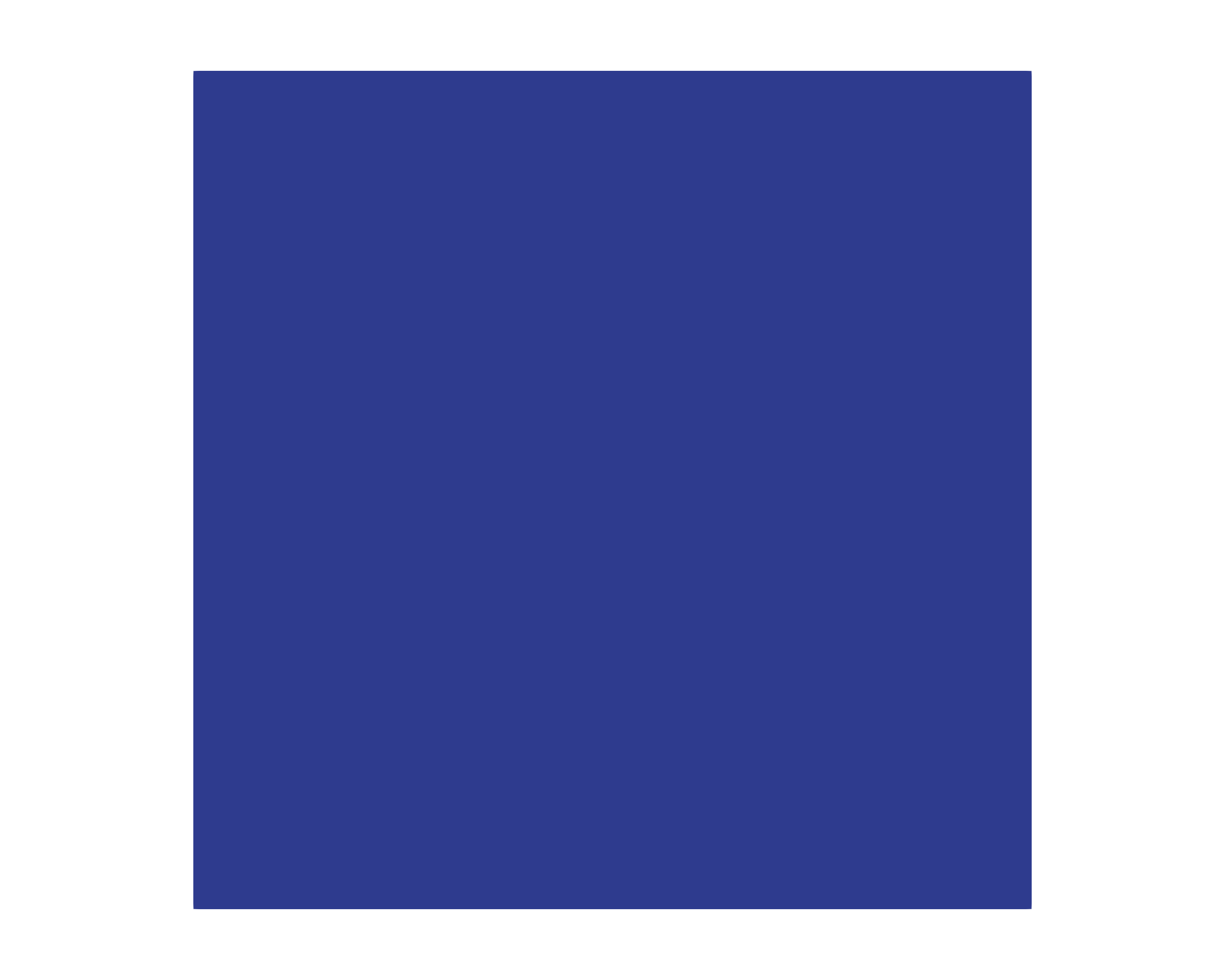}
    }
    \\
    \subfloat[{
        {
          $p=4$
        }  
    }]{
      \includegraphics[scale=\figscale,width=0.45\figwidth]{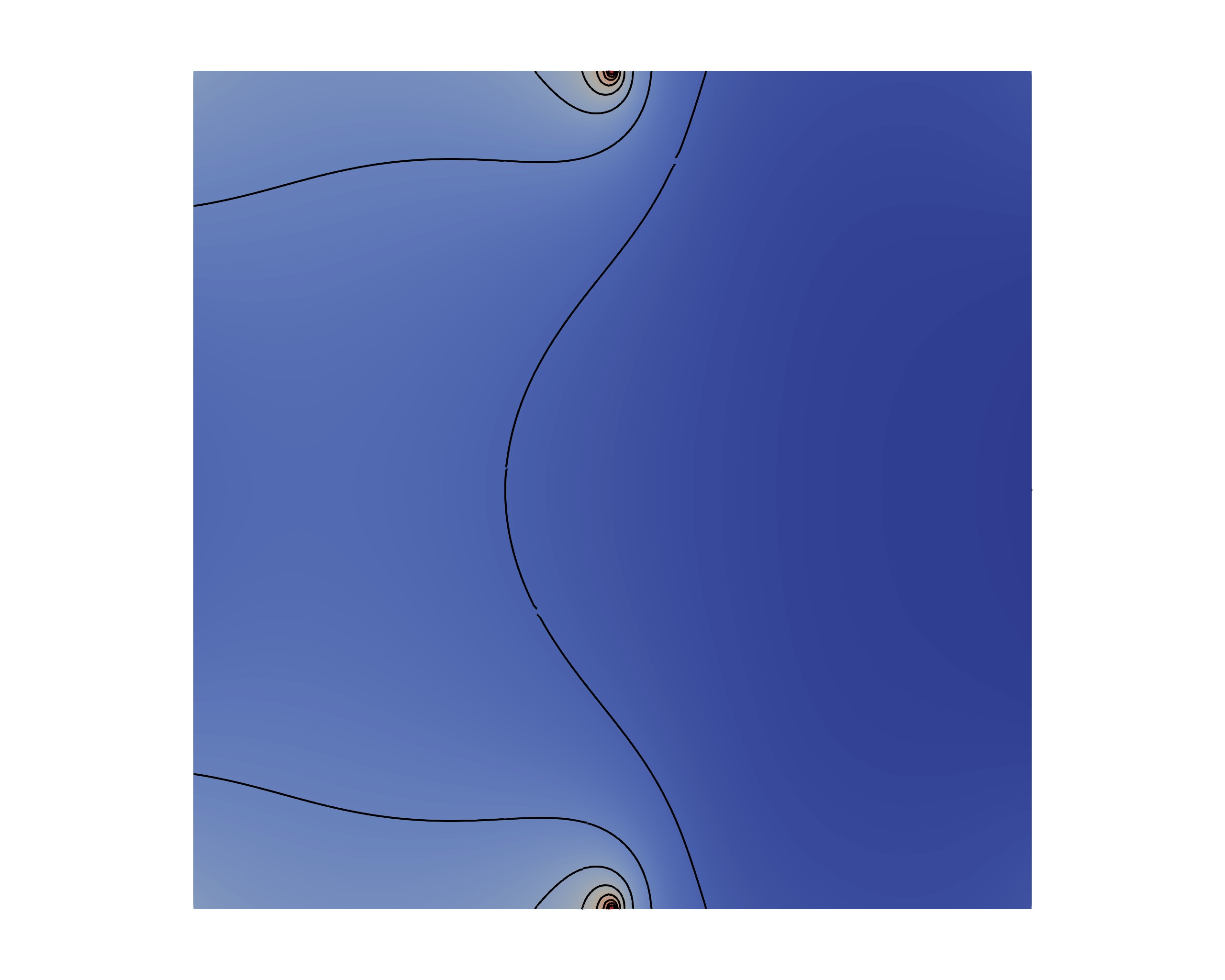}
    }
     \hfill
    \subfloat[{
        %%%%%%%%%%%%%%%%%%%%%%%%%%%%%%%%%%%%%%%%%%%%%%%%%%%%%%%%%%%%%%% 
        {
          $p=10$
        }
        %%%%%%%%%%%%%%%%%%%%%%%%%%%%%%%%%%%%%%%%%%%%%%%%%%%%%%%%%%%%%%%% 
    }]{
      \includegraphics[scale=\figscale,width=0.45\figwidth]{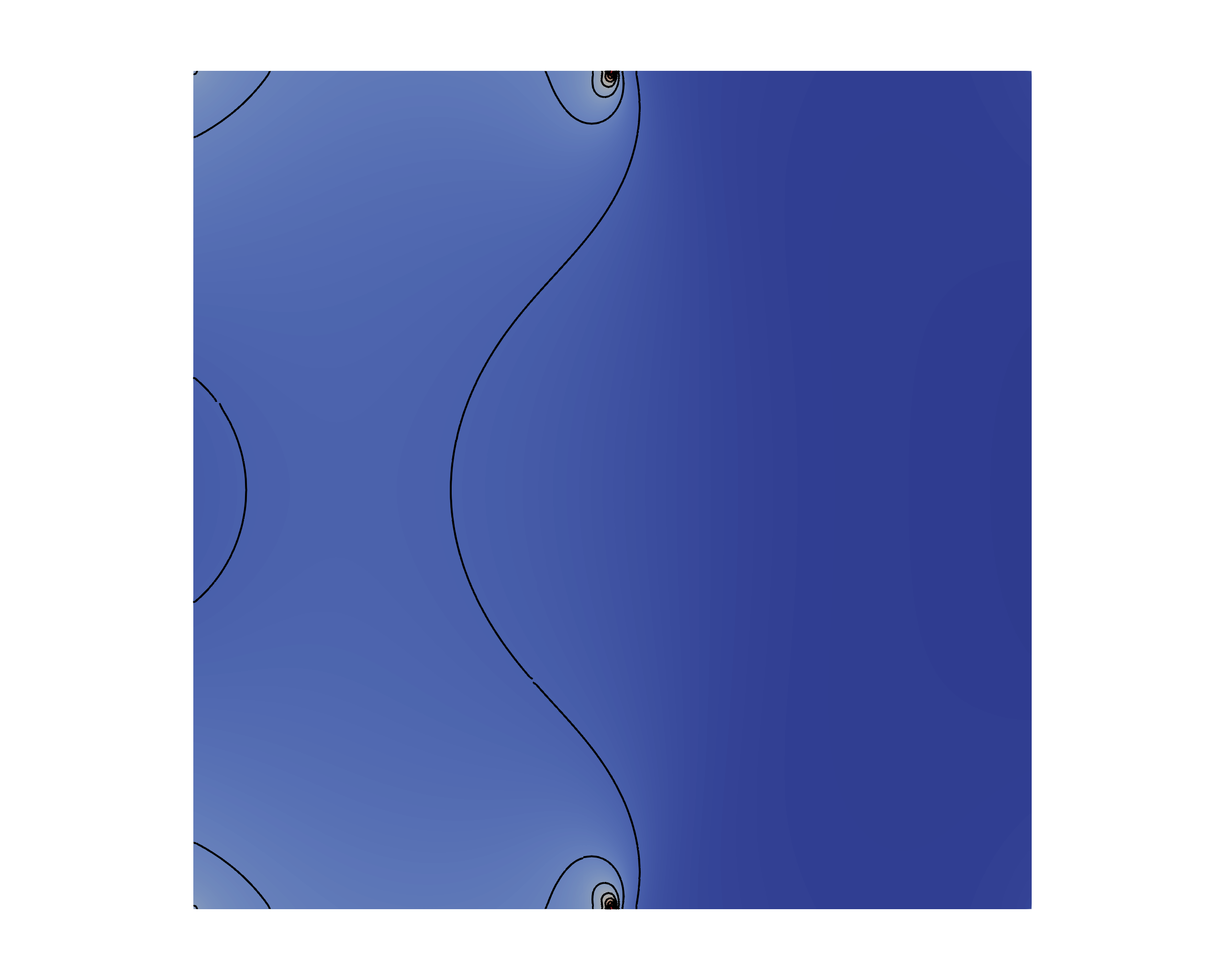}
        }
        \\
      \subfloat[{
      %%%%%%%%%%%%%%%%%%%%%%%%%%%%%%%%%%%%%%%%%%%%%%%%%%%%%%%%%%%%%%% 
        {
          $p=20$
        }
        %%%%%%%%%%%%%%%%%%%%%%%%%%%%%%%%%%%%%%%%%%%%%%%%%%%%%%%%%%%%%%%% 
    }]{
      \includegraphics[scale=\figscale,width=0.45\figwidth]{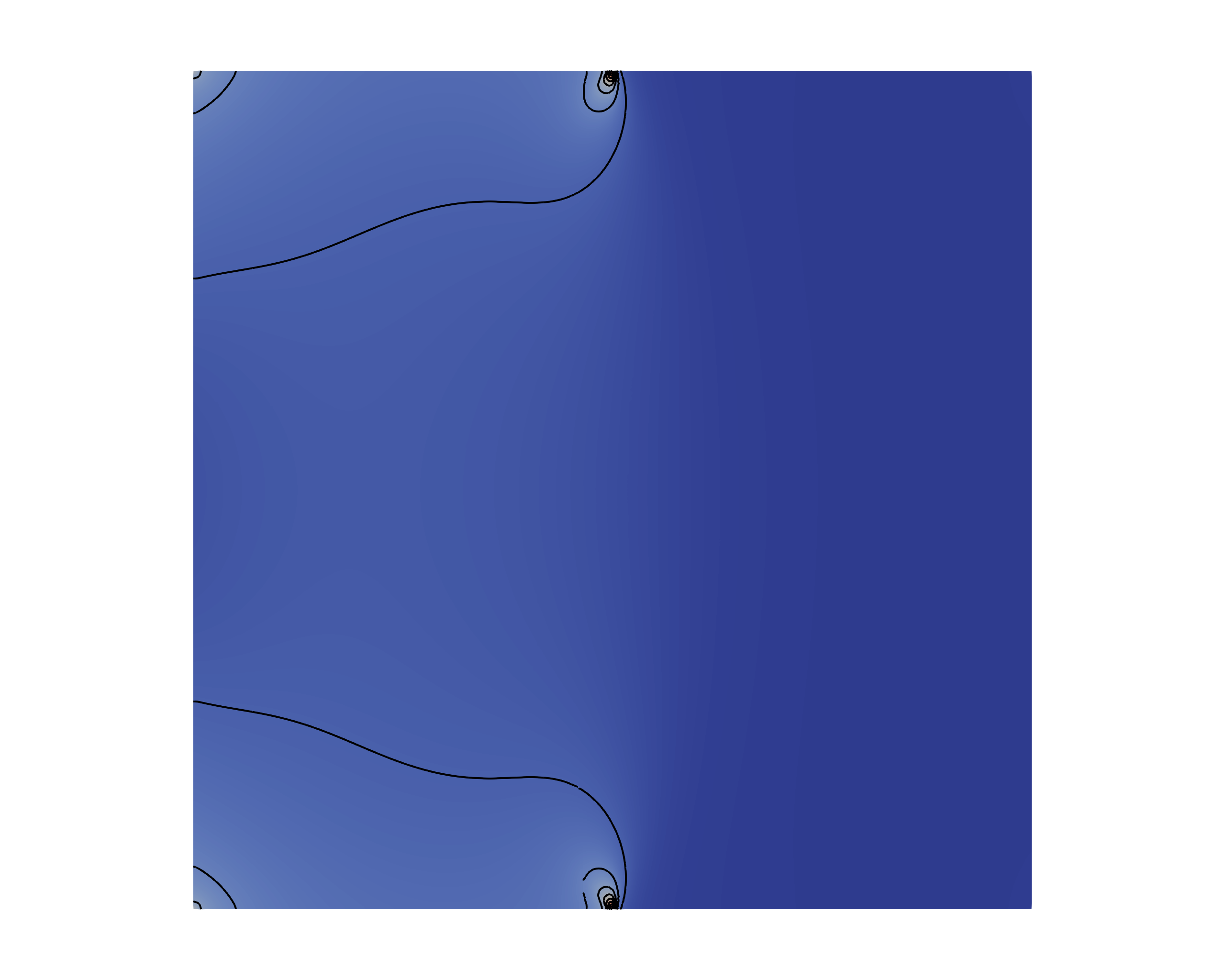}
    }
    %%%%%%%%%%%%%%%%%%%%%%%%%%%%%%%%%%%%%%%%%%%%%%%%%%%%%%%%%%%%%%%%%%%% 
     \hfill
     \subfloat[{
      %%%%%%%%%%%%%%%%%%%%%%%%%%%%%%%%%%%%%%%%%%%%%%%%%%%%%%%%%%%%%%% 
        {
          $p=100$
        }
        %%%%%%%%%%%%%%%%%%%%%%%%%%%%%%%%%%%%%%%%%%%%%%%%%%%%%%%%%%%%%%%% 
    }]{
       \includegraphics[scale=\figscale,width=0.45\figwidth]{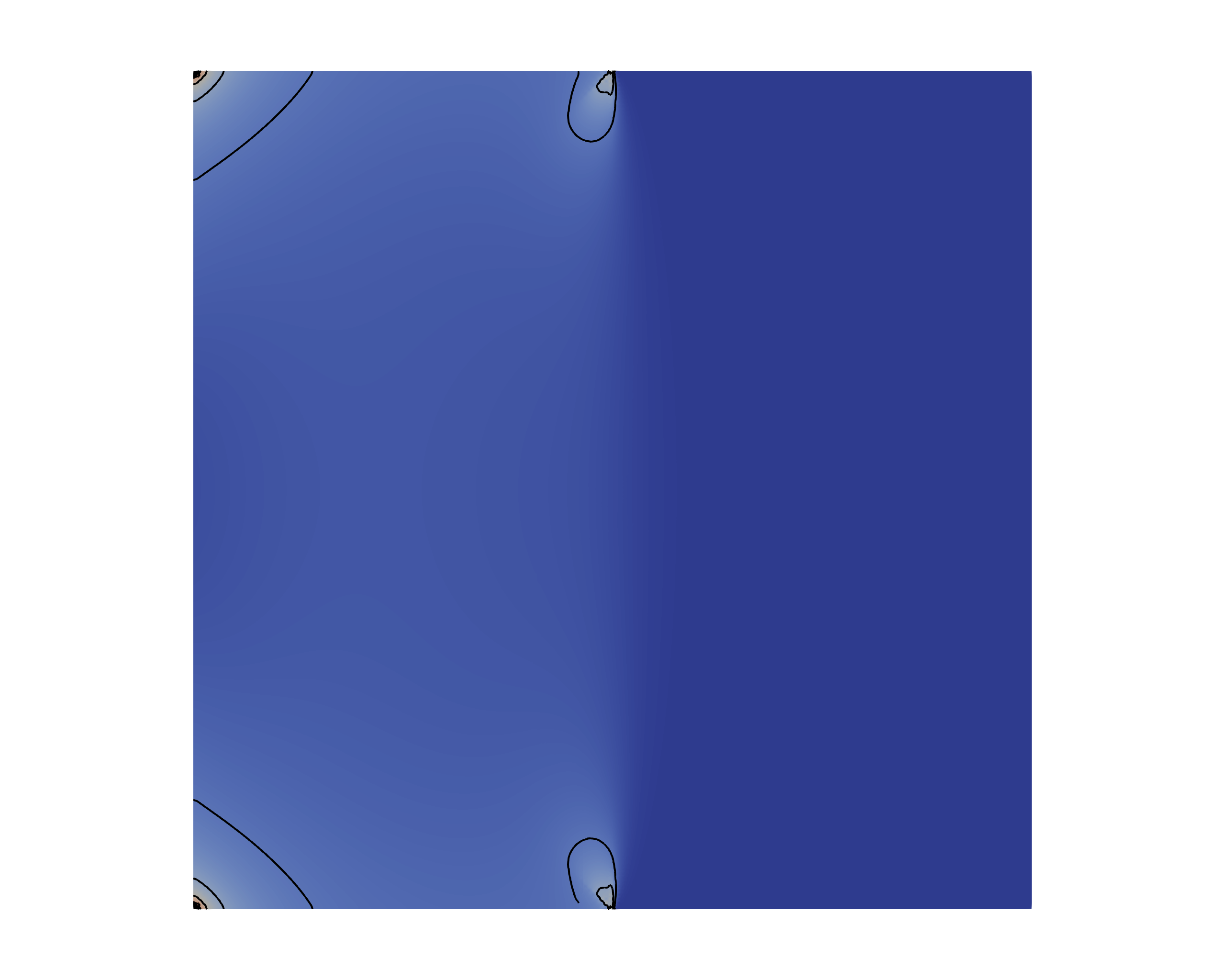}
    }
    %%%%%%%%%%%%%%%%%%%%%%%%%%%%%%%%%%%%%%%%%%%%%%%%%%%%%%%%%%%%%%%%%%%% 
  \end{center}
\end{figure}

In this test we examine the vectorial problem with prescribed boundary
data given by
 \begin{equation}
   \label{eq:lamb}
   u(x,y) =
   \begin{cases}
     (x,y) ,\ \ \ \text{ if } x \leq 0,
     \\
     (\lambda x, y),\ \, \text{ if } x > 0,
   \end{cases}
 \end{equation}
 for $\lambda \in \{\pm\frac 12\}$. We study the measure
 $\sigma_p^\mO(\mathcal T)$ with $\mO = \Omega$. Results are shown in
 Figure \ref{fig:mixed1} and \ref{fig:mixed2} where we plot the domain $\mO$ and
 $\sigma_p^\mO(\mathcal T)$ for $p=2,4,10,20,100$. Specific conjectures regarding the behaviour of the measure $\si^\mO$ is provided at the captions, since the particular behaviour depends on the shape of the domain and its position in relation to the level sets of the gradient. This case appears to combine both concentration and absolutely continuous parts.

\ms

\subsection{Vectorial case - Orientation preserving diffeomorphism}
 
\begin{figure}[h!]
  \caption[] {\label{fig:circ} A test to characterise
    $\sigma^\mO_\infty$ for the boundary data given by (\ref{eq:diff})
    with $\lambda = \frac 1 2$. We show the domain $\mO$ and the
    measure $\sigma_p^\mO(\mathcal T)$ for increasing values of
    $p$. In each of the subfigures B-F we plot 10 contours that are
    equally spaced between the minimum and maximum values of
    $\sigma_p^\mO$. Notice that as $p$ increases these contours are
    relatively evenly spaced. Again the conjecture is that $\si^\mO$ is absolutely continuous with respect to $\mL^n\LL_\mO$.}
  \begin{center}
    \subfloat[{}]
             {
               \hspace{.2cm}
   \begin{tikzpicture}[scale=2]
     \draw [very thin, gray] (-1,-1) grid (1,1);  
     \fill[red,even odd rule] (0,0) circle (.8) (0,0) circle (.2);
     \draw (0, 0) circle (.2);
     \draw (0, 0) circle (.8);
     \begin{scope}[thick,font=\scriptsize]
       \draw [->] (-1.25,0) -- (1.25,0) node [above left]  {$x$};
       \draw [->] (0,-1.25) -- (0,1.25) node [below right] {$y$};
       \foreach \n in {-1,1}{%
         \draw (\n,-3pt) -- (\n,3pt)   node [above] {$\n$};
         \draw (-3pt,\n) -- (3pt,\n)   node [right] {$\n$};
       }
     \end{scope}
   \end{tikzpicture}
             }
     \hfill
     \subfloat[{
        {
          $p=2$
        }
    }]{
      \includegraphics[scale=\figscale,width=0.45\figwidth]{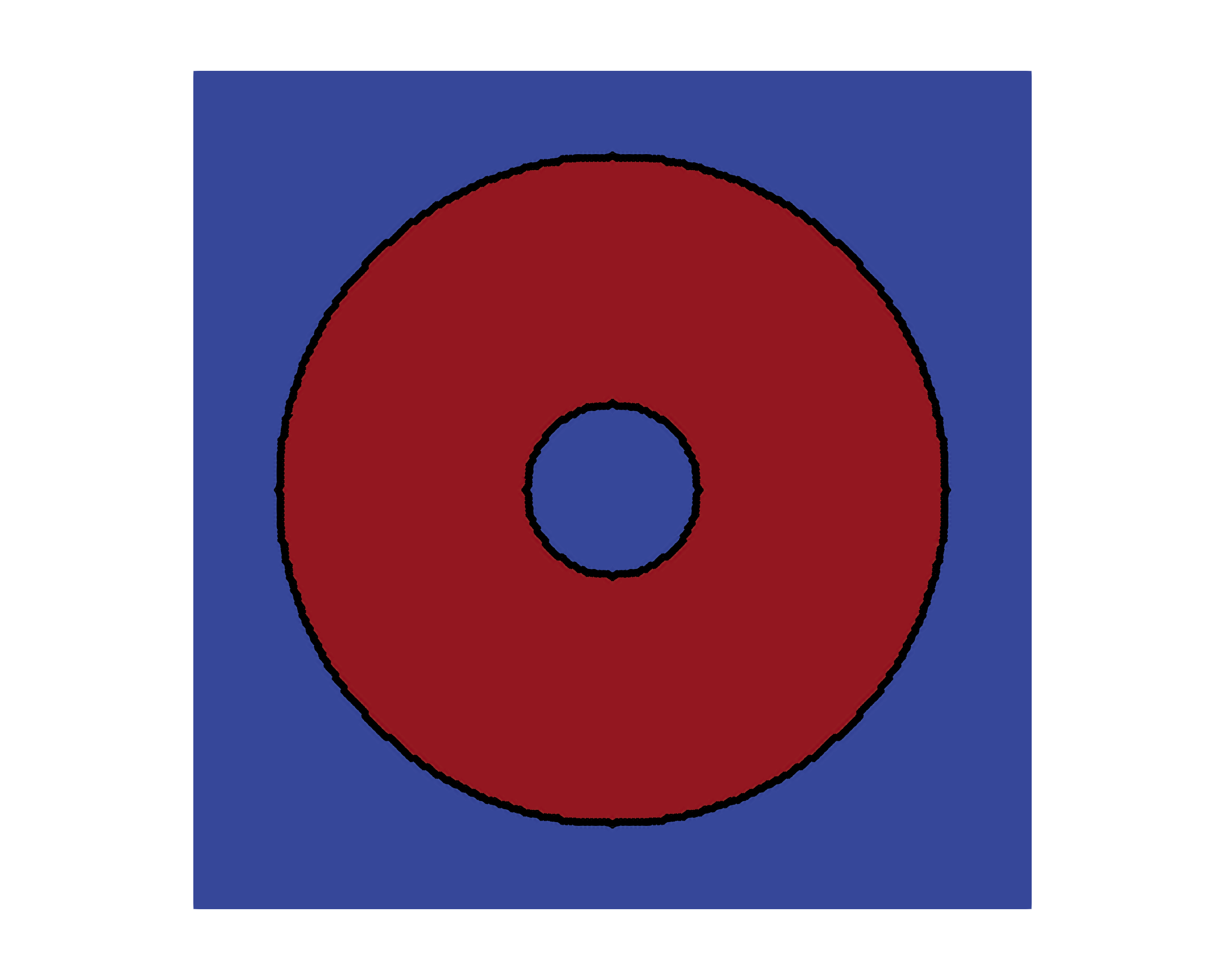}
    }
    \\
    \subfloat[{
        {
          $p=4$
        }  
    }]{
      \includegraphics[scale=\figscale,width=0.45\figwidth]{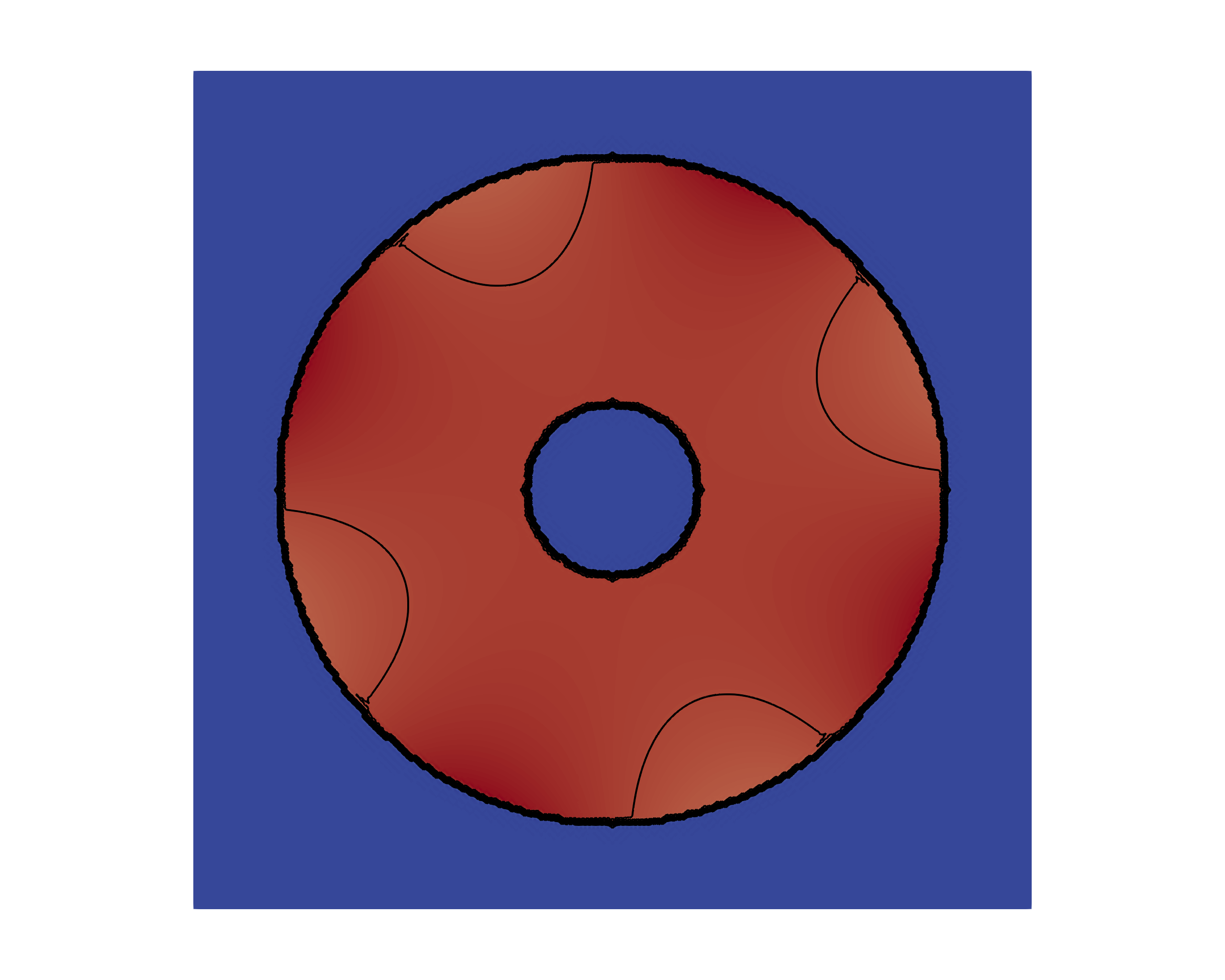}
    }
     \hfill
    \subfloat[{
        %%%%%%%%%%%%%%%%%%%%%%%%%%%%%%%%%%%%%%%%%%%%%%%%%%%%%%%%%%%%%%% 
        {
          $p=10$
        }
        %%%%%%%%%%%%%%%%%%%%%%%%%%%%%%%%%%%%%%%%%%%%%%%%%%%%%%%%%%%%%%%% 
    }]{
      \includegraphics[scale=\figscale,width=0.45\figwidth]{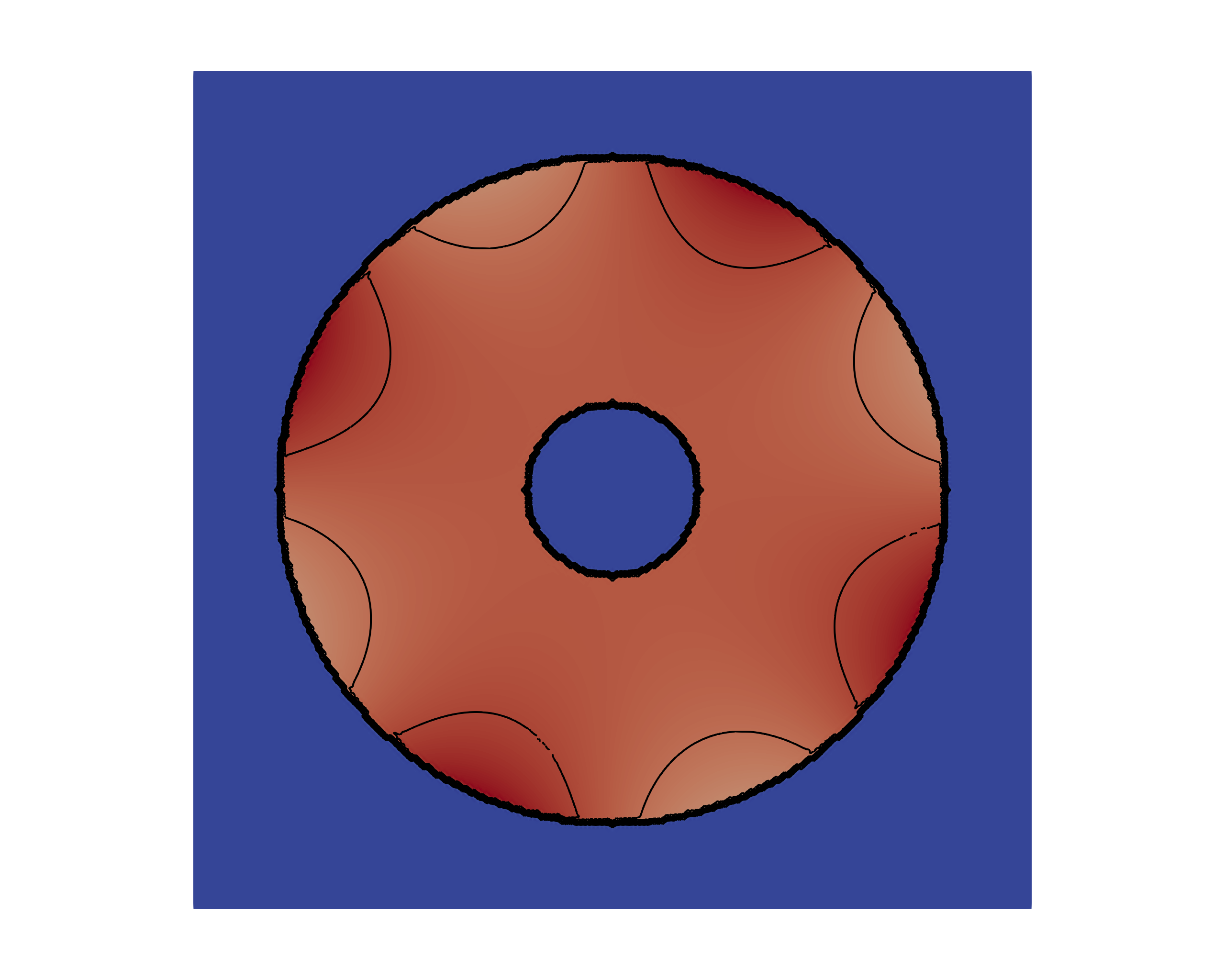}
        }
        \\
      \subfloat[{
      %%%%%%%%%%%%%%%%%%%%%%%%%%%%%%%%%%%%%%%%%%%%%%%%%%%%%%%%%%%%%%% 
        {
          $p=20$
        }
        %%%%%%%%%%%%%%%%%%%%%%%%%%%%%%%%%%%%%%%%%%%%%%%%%%%%%%%%%%%%%%%% 
    }]{
      \includegraphics[scale=\figscale,width=0.45\figwidth]{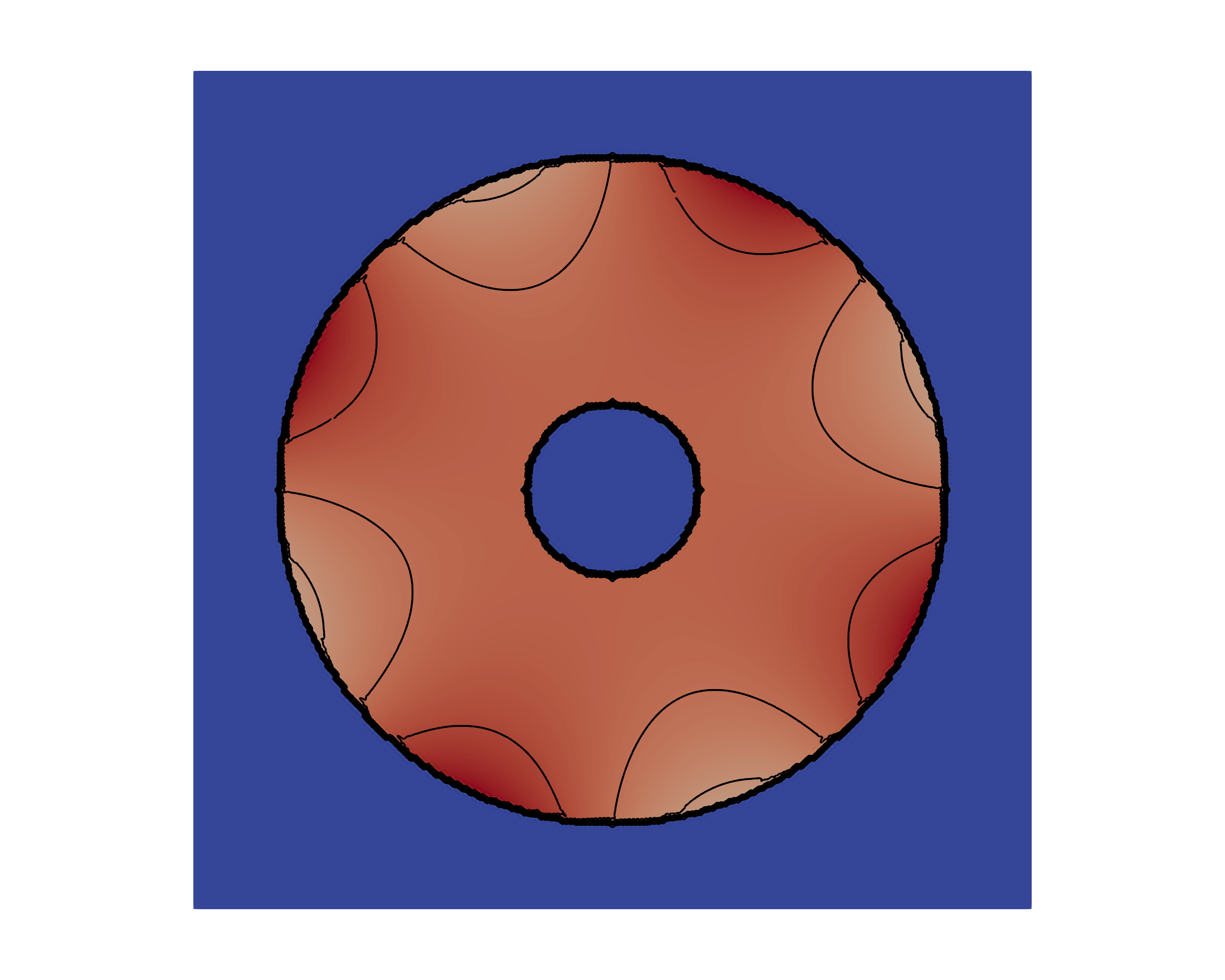}
    }
    %%%%%%%%%%%%%%%%%%%%%%%%%%%%%%%%%%%%%%%%%%%%%%%%%%%%%%%%%%%%%%%%%%%% 
     \hfill
     \subfloat[{
      %%%%%%%%%%%%%%%%%%%%%%%%%%%%%%%%%%%%%%%%%%%%%%%%%%%%%%%%%%%%%%% 
        {
          $p=100$
        }
        %%%%%%%%%%%%%%%%%%%%%%%%%%%%%%%%%%%%%%%%%%%%%%%%%%%%%%%%%%%%%%%% 
    }]{
       \includegraphics[scale=\figscale,width=0.45\figwidth]{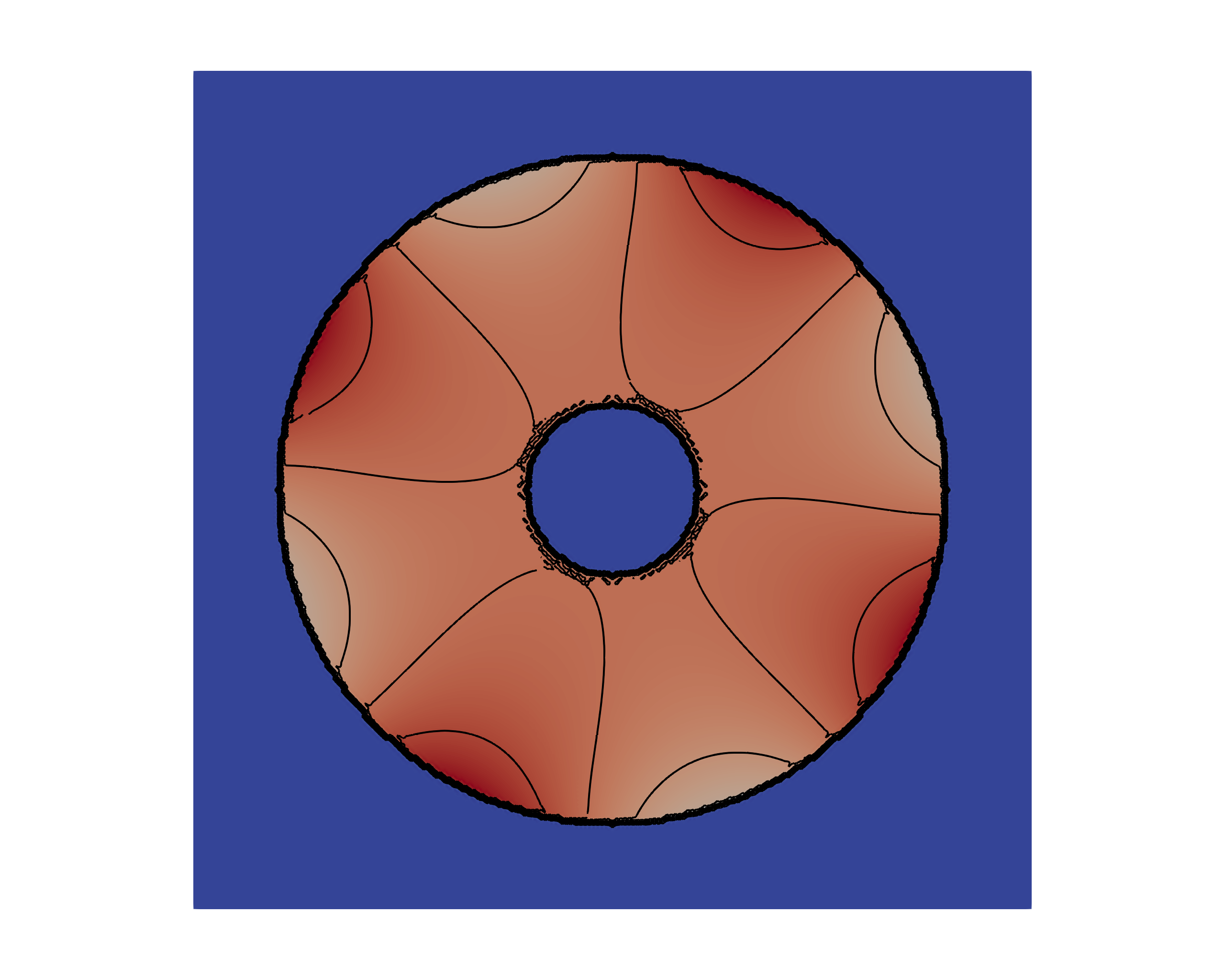}
    }
    %%%%%%%%%%%%%%%%%%%%%%%%%%%%%%%%%%%%%%%%%%%%%%%%%%%%%%%%%%%%%%%%%%%% 
  \end{center}
\end{figure}

In this test we examine the vectorial problem with prescribed boundary
data given by
\begin{equation}
  \label{eq:diff}
  u(x,y) =
  e^{(\log{|x|} S)} x,
\end{equation}
where $S$ is the orthogonal, skew-symmetric matrix
\begin{equation}
  S = 
  \begin{bmatrix}
    0 & -1 \\ 1 & 0
  \end{bmatrix}.
\end{equation}
This is the explicit example given in \cite[Lemma 3.1]{KS}. Results
are shown in Figure \ref{fig:circ} where we plot the domain $\mO$ and
$\sigma_p^\mO(\mathcal T)$ for $p=2,4,10,20,100$. In this case we also conjecture that  $\si^\mO$ must be absolutely continuous with respect to the Lebesgue measure $\mL^n\LL_\mO$. Given that absolute continuity appears in all cases we have eikonal data, it appears that this must be a general fact regardless of any particular additional structures.

%We show when the
%boundary of $\mO$ is not aligned with the level sets of $\norm{\D u}$.

\clearpage

\clearpage

\clearpage

\clearpage

\newpage

\bibliographystyle{amsplain}

\end{document}